\input amssym.def
\input amssym
\documentstyle[12pt]{article}
\title{{\bf Framed vertex operator algebras, codes and the moonshine module}
\footnotetext{The first author is supported by NSF grants DMS-9303374,
DMS-9700923 and a
research grant from the Committee on Research, UC Santa Cruz.}\footnotetext{
The second author is supported by NSF grant DMS-9623038 and the
University of Michigan Faculty Recognition 
Grant (1993--96).}\footnotetext{The third author is supported by a
research fellowship of the DFG, grant Ho~1842/1-1.}}
\author{Chongying Dong,\ Robert L.~Griess Jr.\ and\ Gerald H\"ohn}
\date{}

\makeatletter

\begin{document}

\newtheorem{thm}{Theorem}[section]
\newtheorem{prop}[thm]{Proposition}
\newtheorem{lem}[thm]{Lemma}
\newtheorem{rem}[thm]{Remark}
\newtheorem{cor}[thm]{Corollary}
\newtheorem{conj}[thm]{Conjecture}
\newtheorem{de}[thm]{Definition}
\newtheorem{notation}[thm]{Notation}
\pagestyle{plain}
\maketitle

\def \Z{\Bbb Z}
\def \F{\Bbb F}
\def \C{\Bbb C}
\def \R{\Bbb R}
\def \Q{\Bbb Q}
\def \N{\Bbb N}
\def \D{{\cal D}}
\def \wt{{\rm wt}}
\def \tr{{\rm tr}}
\def \sp{{\rm span}}
\def \Res{{\rm Res}}
\def \Res{{\rm QRes}}
\def \End{{\rm End}}
\def \E{{\rm End}}
\def \Ind {{\rm Ind}}
\def \Irr {{\rm Irr}}
\def \Aut{{\rm Aut}}
\def \Hom{{\rm Hom}}
\def \mod{{\rm mod}}
\def \ann{{\rm Ann}}
\def \rank{{\rm rank}}
\def \<{\langle} 
\def \>{\rangle} 
\def \t{\tau }
\def \a{\alpha }
\def \e{\epsilon }
\def \l{\lambda }
\def \L{\Lambda }
\def \g{\gamma}
\def \b{\beta }
\def \om{\omega }
\def \o{\omega }
\def \k{\kappa}
\def \c{\chi}
\def \ch{\chi}
\def \cg{\chi_g}
\def \ag{\alpha_g}
\def \ah{\alpha_h}
\def \ph{\psi_h}
\def \pf{{\bf Proof. }}
\def \voa{{vertex operator algebra\ }}
\def \svoa{{super vertex operator algebra\ }}
\def \qed{\mbox{$\square$}}
\def \lc{L_C}
\def \tlc{\widetilde{L}_C}
\def \tv{\widetilde{V}_L}
\def \vlc{V_{L_C}}
\def\tvlc{\widetilde{V}_{L_C}}
\def\vtlc{V_{\widetilde{L}_C}}
\def\tvtlc{\widetilde{V}_{\widetilde{L}_C}}
\def\ha{\frac{1}{2}}
\def\se{\frac{1}{16}}

\def\Ve{V^{0}}
\def\M{\rm \ {\sl  M}\llap{{\sl I\kern.80em}}\ }
\def\xx{\em}


{\bf Abstract.} For a simple vertex operator algebra whose Virasoro element
is a sum of commutative Virasoro elements of central charge~$\frac{1}{2}$,
two codes are introduced and studied. It is proved that such vertex operator
algebras are rational. For lattice vertex
operator algebras and related ones, decompositions into
direct sums of irreducible modules for the product of the Virasoro algebras 
of central charge~$\frac{1}{2}$ are explicitly described. As an application, 
the decomposition of the moonshine vertex operator algebra is obtained for a 
distinguished system of 48 Virasoro algebras.

\section{Introduction}

Vertex operator algebras (VOAs) have been studied by mathematicians for  
more than a decade, but still very little is known 
about the general structure of VOAs.  
Most of the examples so far come 
{}from an auxiliary mathematical structure like affine Kac-Moody algebras,
Virasoro algebras, integral lattices or are modifications of these
(like orbifolds and simple current extensions).  We use the definition
of VOA as in [FLM], Section~8.10. 

\smallskip 

In this paper we  develop a general structure theory for a class 
of VOAs containing a subVOA of the same rank and relatively simple form, namely
a tensor product of simple Virasoro VOAs of central charge $\frac{1}{2}$.   
We call this the class of {\sl framed VOAs\/}, abbreviated FVOAs. It contains
important examples of VOAs. We show how
VOAs constructed from certain integral lattices can be described as framed
VOAs.   In the case that the lattice itself comes from a binary code, 
this can be done even more explicitly.
As an application of the general structure theory we describe  
VOAs of small central charge as FVOAs, especially
the moonshine VOA $V^{\natural}$ of central charge~$24$.

\smallskip

The modules of a VOA together with the intertwining operators can be put 
together into a larger structure which is a called an intertwining 
algebra~[Hu1], [Hu2]. In the case where the fusion algebra of the VOA is
the group algebra of an abelian group $G$, like for lattice VOAs,
this specializes to an {\sl abelian\/} intertwining algebra~[DL1];
also see [Mo]. 
The description of VOAs containing a fixed VOA with abelian intertwining 
algebra is relatively simple: They correspond to the subgroups
$H\leq G$ such that all the conformal weights of the VOA-modules
indexed by $H$ are integral~[H3].
The Virasoro VOA of rank $\ha$ gives one of the easiest examples
of {\sl non abelian\/} intertwining algebras. Section~2 can be considered 
as a study of the extension problem for the tensor products of this 
Virasoro VOAs.

\smallskip

It is our hope that the ideas used in this work
can  be  extended to structure theories for
VOAs based on  other classes of rational subVOAs with
nonabelian intertwining algebras, like the VOAs
belonging to the discrete series
representations of the Virasoro algebra~[W].

\smallskip

We continue with a more detailed description of the results in this paper.

\medskip 

The Virasoro algebra of central charge $\frac{1}{2}$  has just three
irreducible highest weight unitary representations, 
with highest weights
$h=0$, ${1 \over 2}$, ${1 \over 16}$, and  the one with $h=0$ carries the
structure of a {\sl simple\/} VOA whose irreducible modules are exactly these 
irreducible highest weight unitary representations.
The relevant fusion rules here (Theorem~\ref{t2.1}) are relatively simple-looking.
A tensor product of $r$ such VOAs, 
denoted $T_r$, has irreducible representations in bijection with $r$-tuples
$(h_1,\dots ,h_r)$ such that each $h_i \in \{0, {1 \over 2}, {1 \over 16} \}$.

We are interested in the case of a VOA 
$V$ containing a subVOA isomorphic to $T_r$.
Such a subVOA arises from a {\sl Virasoro frame\/}, a set 
of elements $\omega_1$, $\ldots$, $\omega_r$ 
such that for each $i$, the vertex operator components of $\omega_i$ along 
with the vacuum element span a copy of the simple Virasoro VOA
of central charge $1 \over 2$ and
such that these subVOAs are mutually commutative
and $\omega_1 + \cdots + \omega_r$ is the Virasoro element of $V$.  
We abbreviate VF for Virasoro frame. 
Such elements may be characterized internally up to a factor $2$
as the unique indecomposable idempotents in the weight $2$ 
subalgebra of $T_r$ with respect to the algebra product 
$u_1v$ induced from the VOA structure on $T_r$. 

\smallskip

It was shown in~[DMZ] that the moonshine VOA $V^{\natural}$ is a FVOA
with $r=48$.  Partial results on decompositions of  $V^{\natural}$
into a direct sum of irreducible $T_{48}$-modules were obtained
in~[DMZ] and [H1]. These results were fundamental in proving that
$V^{\natural}$ is holomorphic~[D3]. 
In fact, the desire to understand $V^\natural$ was one of the
original motivations for us to study FVOAs. 

\medskip

In Section~2, we describe how the set of $r$-tuples which occur lead to two
linear codes ${\cal C}$, ${\cal D} \le \Bbb F_2^r$ where ${\cal D}$ is contained in
the annihilator code ${\cal C^{\perp}}$.
For self-dual FVOAs we give a proof that they are equal: 
$\cal C = \cal D ^\perp$.
Associated to these codes are normal 2-subgroups $G_{\cal D} \le G_{\cal C}$
of the subgroup $G$ of the automorphism group $Aut(V)$
of $V$ which stabilizes the VF (as a set).  The group $G$ is finite.
We get an accounting of all subVOAs of $V$ which contain $\Ve$, the
subVOA of $G_{\cal D}$-invariants.  
We obtain a general result (Theorem~\ref{tr1}) that FVOAs are rational,
establishing the existence of a new broad class of rational VOAs.  The
rationality of  FVOAs  is a very  important aspect of  their  representation
theory. In particular, a FVOA has only finitely many irreducible modules.

\smallskip

In Section 3, we describe the Virasoro decompositions of the lattice VOAs 
$V_{D_1^d}$, and closely related VOAs, with respect to a natural subVOA $T_{2d}$.

\smallskip

In Section 4, we study the familiar situation of the twisted or untwisted
lattice associated to a binary doubly-even code of length $d \in 8 \Bbb Z$ 
and the twisted and untwisted VOA associated to a lattice.
A {\sl marking\/} of the code is a partition of its coordinates into $2$-sets.
A marking determines a $D_1^d$ sublattice in the associated lattices
and a VF in the associated VOAs.
We give an explicit description of the coset decomposition of the lattices 
under the $D_1^d$ sublattice, a $\Z_4$-code,
and the decomposition of the VOA as a module for the subVOA generated by
the VF. As a corollary, we give 
information about various multiplicities of the decompositions under 
this subVOA using the symmetrized marked weight enumerator of the marked code
or the symmetrized weight enumerator of the $\Bbb Z_4$-code.

\smallskip

Finally, Section 5 is devoted to applications.  Two examples are
discussed in detail.  Example~I is about the Hamming code of length~$8$, the
root lattice $E_8$ and the VOA $V_{E_8}$. Here, $r=16$, and we find
at least $5$ different VFs.
Example~II is about the Golay code, the Leech lattice
and the moonshine module, $V^{\natural}$, where $r=48$.
For every  VF inside $V^{\natural}$, the code ${\cal C}$ has 
dimension at most~41. 
There is a special marking for which this
bound of~$41$ is achieved, and for this marking the complete
decomposition polynomial is explicitly given.
The $D_1$-frames inside the Leech lattice which arise
{}from a marking of the Golay code are characterized 
by properties of the corresponding $\Z_4$-codes. 

\smallskip   

Appendix~\ref{hamming-code} contains a few special results about orbits
on markings of the Hamming code, 
Appendix~\ref{golay-code} the stabilizer in $M_{24}$ of the above
special marking for the Golay code, and Appendix~\ref{moon-code} the structure
of the automorphism group of the above code of dimension~$41$. 
Appendix~\ref{append-lift} classifies automorphisms of a lattice VOA which 
correspond to $-1$ on the lattice.

\medskip

In [M1]--[M3],  there is a new treatment of the moonshine VOA and there
is some overlap with results of this article.  In particular, the vertex
operator subalgebra similar to
our $V^{0}$ (see section~2) and its representation theory have been
independently investigated in [M3].

\bigskip

The authors are grateful to E.~Betten for discussions on markings of the
Golay code.   They also thank A.~Feingold and G.~Mason for discussions.

\bigskip


{}{\bf{}Notation and terminology}
\bigbreak
\halign{#\hfil&\quad#\hfil\cr

${\bf 1}$ & The vacuum element of a VOA. \cr
${ Aut}(V)$ & The automorphism group of the VOA $V$.\cr 
$B_V$ & The conformal block on the torus of the VOA $V$.\cr
${\cal B}_2^n$  & The FVOA $(M(0,0)\oplus M(\ha,\ha))^{\otimes n}$
with binary code \cr
& ${\cal C}({\cal B}_2^n)=\{(0,0),(1,1)\}^n$ of length $2n$. \cr
$({\cal B}_2^n)_0$  & The subVOA of ${\cal B}_2^n$ belonging to the
subcode of ${\cal C}({\cal B}_2^n)$ consisting \cr 
& of codewords of weights divisible by $4$. \cr
$c$ & An element of $\F_2^n$. \cr
$C$ & A linear binary code, often self-annihilating and doubly-even.\cr
$C^{\perp}$ & The annihilator code of $C$. \cr
${\cal C}={\cal C}(V)$& The binary code determined by the $T_r$-module structure
of $\Ve$.  \cr
$\C [L]$ & The complex group algebra of the group $L$. \cr
$\C \{ L \}$ & The twisted complex group algebra of the lattice $L$; \cr
 & it is the group algebra $\C [ \hat L ]$ modulo the ideal generated by
$\kappa + 1$. \cr 
$Co_0$ & The Conway group which is $Aut(\Lambda)$, a finite group of order  \cr
 & $2^{22}3^95^47^211.13.23$; its quotient by the center $\{ \pm 1 \}$ is a finite  \cr
 & simple group. \cr
$d$ & The length of a binary code $C$, usually divisible by $8$. \cr
$d_4^n$  & The marked binary code $\{(0,0,0,0),(1,1,1,1)\}^n$ of length $4n$. \cr 
$(d_4^n)_0$  & The subcode of $d_4^n$ consisting of codewords of 
weights divisible by $8$. \cr
$D_n$ & The index $2$ sublattice of $\Z^n$ consisting of vectors  whose\cr
& coordinate sum is even (the ``checkerboard lattice"). \cr
${\cal D}={\cal D}(V) $ &  The binary code of the $I\subseteq 
\{1,\ldots,r\}$ with
$V^I\not = 0$. \cr
$\delta_2^n$  & The marked Kleinian or $\F_4$-code $\{(0,0),(1,1)\}^n$ of length 
$2n$. \cr
$(\delta_2^n)_0$  & The subcode of $\delta_2^n$ consisting of codewords of 
weights divisible by $4$. \cr
$\delta(c)$ & The number of $k$ with $c(k)=(c_{2k-1},c_{2k})\in \{(0,1),(1,0)\}$.\cr
$\Delta(L)$ & The $\Z_4$-code associated to a lattice $L$ with fixed
$D_1$-frame. \cr
$E_8$ & The root lattice of the  Lie group $E_8(\C)$.\cr
$\epsilon$ & A vector with components $+$ or $-$. \cr
FVOA & Abbreviation for framed vertex operator algebra. \cr
$\cal F$ & The set $\{M(0),M(\frac{1}{2}),M(\frac{1}{16})\}$.\cr
$G$ & The subgroup of ${ Aut}(V)$ fixing a VF of $V$.  \cr 
${G_{\cal C}} $& The normal subgroup of $G$ acting trivially on $T_r$.\cr
${G_{\cal D}} $& The normal subgroup of $G$ acting trivially on $\Ve$.\cr
${{\cal G}={\cal G}_{24}} $& The Golay code of length $24$. \cr
$\gamma$ &  An element of $\Z_4^n$. \cr
$\gamma_{\epsilon_k}^a$ & A map $\F_2^2\longrightarrow \Z_4^2$.\cr 
$\Gamma_{\epsilon}^a$ & A map $\F_2^{2n}\longrightarrow \Z_4^{2n}$.\cr
$h$, $h_i$ & weights of elements or modules of a VOA,
 usually $h_i \in \{0,\frac{1}{2},\frac{1}{16}\}$. \cr
${H_8} $& The Hamming code of length $8$. \cr
${\cal H}={\cal H}_6$ & The hexacode of length $6$, a code over 
$\F_4=\{0,1,\omega,\bar\omega\}$  or over \cr
 & the Kleinian fourgroup $\Z_2\times \Z_2=\{0,a,b,c\}$. \cr
$I$ & A subset of $\{1,\ldots,r\}$. \cr
$I+J$ & The symmetric difference, for subsets of $\{1,\ldots,r\}$. \cr
$\kappa$ & A central element of order 2 in the group $\hat L$. \cr 
$L$ & An integral lattice, often self-dual and even. \cr
$\hat L$ & A central extension of $L$ by a central subgroup $\kappa$. \cr
$L^*$ & The dual lattice of $L$. \cr
$L_C$ & The even lattice constructed from a doubly-even code $C$. \cr
$\tlc$ & The ``twisted'' even lattice constructed from a doubly-even code $C$. \cr
$L(n)$, $L^i(n)$ & The generator of a Virasoro algebra given by the expansion \cr
&  $Y(\omega , z)= \sum_{n \in \Z}L(n)z^{-n-2}$,  
resp.~$Y(\omega_i , z)= \sum_{n \in \Z}L^i(n)z^{-n-2}$.  \cr 
 $\Lambda$ & The Leech lattice. \cr
$\M$ & The Monster simple group. \cr
$M_{24}$ & The simple Mathieu group of order $2^{10}.3^3.5.7.11.23=244,823,040$.\cr
$M(0)$, $M(\frac{1}{2})$, $M(\frac{1}{16})$ & The irreducible modules for the 
   Virasoro algebra with central \cr & charge $\frac{1}{2}$. \cr
$M(1)$ &  The canonical irreducible module for Heisenberg algebras. \cr
$M(h_1,\ldots,h_r)$ & The irreducible $T_r$-module of highest weight 
  $(h_1,\ldots,h_r)$. \cr
$m_h(V)=m_{h_1,\ldots,h_r}$  & The multiplicity of the $T_r$-module 
$M(h_1,\ldots,h_r)$ in the FVOA $V$; \cr
 & we think of this as a function of 
$(h_1,\ldots,h_r) \in \{0, {1 \over 2}, {1 \over 16} \}^r$.  \cr
${\cal M}$ & A marking of a binary code. \cr
$n$ & A natural number. \cr
$N_{\mu_k\epsilon_k}^{ab}$ & A map $\F_2^2\longrightarrow \C[{\cal F}^4]$.\cr
${\bf N}_{\mu,\epsilon}^{ab}$ & A map $\F_2^{2n}\longrightarrow \C[{\cal F}^{4n}]$. \cr
$\mu$ & A vector with components $+$ or $-$. \cr
$P_V(a,b,c)$ & The decomposition polynomial of a FVOA $V$. \cr
$r$ & The number of elements in a VF.\cr
$R_{\mu_k}^a$ & A map $\Z_4\longrightarrow \C[{\cal F}^2]$.\cr
${\bf R}_{\mu}^a$ & A map $\Z_4^n\longrightarrow \C[{\cal F}^{2n}]$. \cr
${\rm smwe}_C(x,y,z)$ & The symmetrized marked weight enumerator of a binary code\cr & $C$
with marking ${\cal M}$. \cr
${\rm swe}_{\Delta}(A,B,C)$ & The symmetrized weight enumerator of a $\Z_4$-code $\Delta$. \cr 
$Sym_r$, $Sym_{\Omega}$& The symmetric group on a set of $r$ objects, usually the index set \cr
& $\{1,\ldots,r\}$, resp.~the symmetric group on the set $\Omega$. \cr
$\Sigma_2^n$  & The $\Z_4$-code $\{(0,0),(2,2)\}^n$ of length $2n$. \cr
$(\Sigma_2^n)_0$  & The subcode of $\Sigma_2^n$ consisting of codewords of 
weights divisible by $4$. \cr
$T$ & A faithful module of dimension $2^m$ for an extraspecial group of  \cr
 & order $2^{1+2m}$, for some $m$, or for a finite quotient of some $\hat L$.  \cr
$T_r=M(0)^{\otimes r}$ & The tensor product of $r$ simple Virasoro VOAs of rank $\frac{1}{2}$. \cr
$V$ & An arbitrary VOA, often holomorphic = self-dual. \cr
$V(c)$ & The submodule of the FVOA $V$ isomorphic
 to  $M(\frac{c}{2})$.\cr
$V_L$ & The VOA constructed from an even lattice $L$. \cr
$V_L^T$ & The $\Z_2$-twisted module of the lattice VOA $V_L$. \cr
$\tv$ & The ``twisted'' VOA constructed from an even lattice $L$. \cr
VF  & Abbreviation for Virasoro frame. \cr
VOA & Abbreviation for vertex operator algebra. \cr
$V^I$ & The sum of irreducible $T_r$-submodules of $V$ isomorphic to \cr  
 & $M(h_1,\ldots,h_r)$ with $h_i=\frac{1}{16}$ if and only if $i\in I$.\cr
$\Ve=V^{\emptyset}$ & This is $V^I$, for $I=0$. \cr
$V^{\natural}$ & The moonshine VOA, or moonshine module. \cr
$W(R)$ & The Weyl group of type $R$, a root system.  \cr 
$Y(\,.\, ,z)$ &  A vertex operator.\cr
$\Xi_1$, $\Xi_3$ &  Two $D_8^*/D_8$-codes of length $1$ and $3$.\cr
$\omega$, $\omega_i$ & Virasoro elements of rank $r$, $\frac{1}{2}$, 
respectively. \cr
$\Omega$ & The ``all ones vector'' $(1,1,\dots,1)$ in $\F_2^n$.\cr
}


\section{Framed vertex operator algebras}\label{one}
\setcounter{equation}{0}

Recall that the Virasoro algebra of central charge $\frac{1}{2}$ has 
three irreducible unitary representations $M(h)$ of highest
weights $h=0$, $\frac{1}{2}$, $\frac{1}{16}$ (cf.~[FQS], [GKO], [KR]).
Moreover, $M(0)$ can be made into a simple vertex operator algebra with 
central charge $\frac{1}{2}$ (cf.~[FZ]). 

In [DMZ], a class of simple vertex operator algebras $(V,Y,{\bf 1},\o)$ containing
an even number of commuting Virasoro algebras of rank $\frac{1}{2}$ were defined. 

\begin{de}{\rm  Let $r$ be any natural number. A simple vertex operator algebra 
$V$ is called a {\em framed vertex operator algebra\/} (FVOA) if the following 
conditions are satisfied: There exist $\omega_i\in V$ for $i=1$,~$\ldots$,~$r$ such that 
(a) each $\o_i$ generates a copy of the simple Virasoro vertex operator 
algebra of central charge $\frac{1}{2}$ and 
the component operators $L^i(n)$ of 
$Y(\omega_i,z)=\sum_{n\in\Z}L^i(n)z^{-n-2}$ 
satisfy $[L^i(m),L^i(n)]=(m-n)L^i(m+n)+\frac{m^3-m}{24}\delta_{m,-n};$
(b) the $r$ Virasoro algebras are mutually commutative; 
and (c) $\omega=\omega_1+\cdots+\omega_{r}$.
The set $\{\o_1,\ldots,\o_r\}$ is called a {\em Virasoro frame\/} (VF).}
\end{de}

In this paper we assume that $V$ is a FVOA.
It follows that $V$ is a unitary representation for each of
the $r$ Virasoro algebras of central charge $\frac{1}{2}$. 

In~[DMZ] it is also assumed that $V_0$ is one-dimensional.
This assumption is now a consequence of the simplicity of $V$: 
\begin{rem}\label{l2.2}   
A FVOA is truncated below from zero: 
$V=\bigoplus_{n\geq 0}V_n$ and $V_0$ is one dimensional: $V_0=\C\,{\bf 1}$. 
\end{rem}

\pf Let $Y(\omega_i, z)= \sum_{n \in \Z}L^i(n)z^{-n-2}$.
Since $V$ is a unitary representation for the Virasoro algebra
generated by the components for $Y(\o,z)=\sum_{n\in\Z}L(n)z^{-n-2}$ 
as $L(n)=\sum_{i=1}^rL^i(n)$ all weights of $V$ are nonnegative
that is, $V=\bigoplus_{n\geq 0}V_n$. 

Then each nonzero vector $v\in V_0$ is a highest weight vector for the 
$r$ Virasoro algebras with highest weight $(0,\ldots,0)$. 
The highest weight module for the $i$-th Virasoro algebra generated by
$v$ is necessarily isomorphic to $M(0)$. From the construction of
$M(0)$ we see immediately that $L^i(0)v=0$.
So $L(-1)v=\sum_{i}L^i(-1)v=0$, i.e.~$v$ is a vacuum-like vector (see~[L1]). 
It is proved in~[L1] that a simple vertex
operator algebra has at most one vacuum-like vector up to a scalar. 
Since ${\bf 1}$ is a vacuum like vector, we conclude that 
$V_0=\C\,{\bf 1}$. \qed

The following theorem can be found in~[DMZ]:
\begin{thm}\label{t2.1} (1) The VOA $M(0)$ has exactly three 
irreducible $M(0)$-modules, namely $M(h)$,  with
$h=0$, $\frac{1}{2}$, $\frac{1}{16}$, and any module is completely reducible. 

(2) The nontrivial fusion rules among these
modules are given by: $M(\frac{1}{2})\times M(\frac{1}{2})=M(0)$, 
$M(\frac{1}{2})\times M(\frac{1}{16})=M(\frac{1}{16})$ and 
$M(\frac{1}{16})\times M(\frac{1}{16})=M(0)+M(\frac{1}{2})$. 

(3) Any module for the tensor product vertex operator algebra
$T_{r}=M(0)^{\otimes r}$, 
where $r$ a positive integer, is a direct sum of irreducible modules 
$M(h_1,\ldots,h_r):=M(h_1)\otimes\cdots\otimes M(h_r)$ 
with $h_i\in\{0,\frac{1}{2},\frac{1}{16}\}$. 

(4) As $T_{r}$-modules, 
$$V=\bigoplus_{h_i\in\{0,\frac{1}{2},\frac{1}{16}\}}
m_{h_1,\ldots, h_{r}}M(h_1,\ldots,h_{r})$$
where the nonnegative integer 
$m_{h_1,\ldots,h_r}$ is the multiplicity of $M(h_1,\ldots,h_r)$ in $V$. 
In particular, all the multiplicities are finite and $m_{h_1,\ldots,h_r}$
is at most $1$ if all $h_i$ are different from $\frac{1}{16}$. 
\end{thm}

\smallskip
Let $I$ be a subset of $\{1,\ldots,r\}$. Define
$V^I$ as the sum of all irreducible submodules isomorphic to
$M(h_1, \ldots , h_r)$ such that $h_i=\frac{1}{16}$ if and only if 
$i \in I$. Then 
$$V=\bigoplus_{I\subseteq \{1,\ldots,r\}}V^I.$$

Here and elsewhere we identify a subset of $\{1,2,\ldots,r\}$ 
with its characteristic function, an integer vector of 
zeros and ones. We further identify such vectors with their image under the 
reduction modulo $2$, i.e.~we consider them as binary codewords in $\F_2^r$.
Interpretation should be clear from context, e.g. we think
of the codeword $c$ as an $r$-tuple of integers in the expression 
$\frac{1}{2}c$.  

For each  $c\in \F_2^r$ let $V(c)$ be the sum of the irreducible submodules 
isomorphic to $M(\frac{1}{2}c_1,\ldots,\frac{1}{2}c_{r})$. 
Then $\Ve=\bigoplus_{c\in \F_2^r}V(c)$. Recall the important fact mentioned in Theorem~\ref{t2.1}~(4) 
that for $c \in {\cal C}$ the 
$T_r$-module $M(\frac{1}{2}c_1,\ldots,\frac{1}{2}c_{r})$ 
has multiplicity $1$ in $V$. So, $V(c)=0$ or is isomorphic to
$M(\frac{1}{2}c_1,\ldots,\frac{1}{2}c_{r})$. 
\smallskip 

We can now define two important binary codes ${\cal C}={\cal C}(V)$ 
and ${\cal D}={\cal D}(V)$. 
\begin{de} \rm For every FVOA $V$,  let:
\begin{equation}\label{2.1}
{\cal C}={\cal C}(V)=\{c\in \F_2^r \mid V(c)\ne 0\}, \hbox{ \ and\ } 
{\cal D}={\cal D}(V)=\{I\in \F_2^r \mid V^I \ne 0\}.
\end{equation}
\end{de}

The vector of all multiplicities $m_{h_1,\ldots,h_r}$ will be
denoted by $m_h(V)$. Note that the codes ${\cal C}$ and 
${\cal D}$ are completely determined by $m_h(V)$. 

The following Proposition generalizes Proposition~5.1 of [DMZ]
and Theorem~4.2.1 of [H1]. In particular, it shows 
${\cal C}$ and ${\cal D}$ are linear binary codes.

\begin{prop}\label{p2.3} (1) $\Ve=V^{\emptyset}$ is a simple vertex operator 
algebra and the $V^I$ are irreducible $\Ve$-modules. Moreover
$V^I$ and $V^J$ are inequivalent if $I\ne J$. 

(2) For any $I$ and $J$ and $0\ne v\in V^J$, $\sp\{u_nv\mid u\in V^I\}=
V^{I+J}$ where $I+J$ is the symmetric difference of $I$ and $J$. 
Moreover, ${\cal D}$ is an abelian group under symmetric difference.

(3) There is one to one correspondence between the subgroups ${\cal D}_0$
of ${\cal D}$ and the vertex operator subalgebras which contain
$\Ve$ via ${\cal D}_0\mapsto V^{{\cal D}_0}$, 
where we define $V^S :=\oplus_{I\in S}V^I$ for any subset $S$ of ${\cal D}_0$. 
Moreover, $V^{I+{\cal D}_0}$ is an irreducible $V^{{\cal D}_0}$-module
for $I\in{\cal D}$ and $V^{I+{\cal D}_0}$ and $V^{J+{\cal D}_0}$
are nonisomorphic if the two cosets are different.  

(4) Let $I \subseteq \{1,\ldots,r\}$ be given and suppose
that $(h_1,\ldots,h_r)$ and $(h_1',\ldots,h'_r)$ are $r$-tuples 
with $h_i$, $h_i'\in\{0,\ha,\se\}$ such that $h_i={1 \over 16}$ 
(resp.~$h_i'={1 \over 16}$) if and only if $i\in I$.
If both $m_{h_1,\ldots,h_r}$ and $m_{h_1',\ldots,h'_r}$  are nonzero 
then $m_{h_1,\ldots,h_r}=m_{h_1',\ldots,h'_r}$.  
That is, all
irreducible modules inside $V^I$ for $T_r$ have the same multiplicities.

(5) The binary code ${\cal C}$ is linear and 
$\{u_nv\mid u\in V(c)\}=V(c+d)$ for any
$c$, $d \in {\cal C}$ and $0\ne v\in V(d)$. 

(6) Moreover, there is a one to one correspondence between vertex operator
subalgebras of $\Ve$ which contain $T_{r}$ and the subgroups
of ${\cal C}$, and $V$ is completely reducible 
for such vertex operator subalgebras whose irreducible modules in
$\Ve$ are indexed by the corresponding cosets in ${\cal C}$.
\end{prop}

\pf Let $v\in V^J$ be nonzero. It follows from Proposition~2.4 of [DM]
or Lemma 6.1.1 of [L2] and the simplicity of $V$ 
that $V=\sp\{u_nv\mid u\in V,n\in\Z\}$. 

\smallskip   
 
{}From the fusion
rules given in Theorem~\ref{t2.1}~(2) and Proposition 2.10 of [DMZ] we see
that $u_nv\in V^{I+J}$ exactly for $u\in V^I$. 
In particular, $\sp\{u_nv\mid u\in \Ve,\,\hbox{$n\in\Z$}\}=V^J$. 
So, $V^J$ can be generated by any nonzero vector and
$V^J$ is a irreducible $\Ve$-module. Since
$V^I$ and $V^J$ are inequivalent $T_{r}$-modules 
if $I\ne J$ they are certainly 
inequivalent $\Ve$-modules. By Proposition~11.9 of [DL1], 
we know that $Y(u,z)v\ne 0$ if 
$u$ and $v$ are not $0$. Thus $V^{I+ J}\ne 0$ if
neither $V^I$ or $V^J$ are 0. This shows that ${\cal D}$ is
a group. So, we finish the proof of (1) and (2).

\smallskip

For (3), we first observe that for a subgroup ${\cal D}_0$ of ${\cal D}$,
(2) implies that $V^{{\cal D}_0}$ is a subVOA which contains $\Ve$.
On the other hand, since $V=V^{{\cal D}}$, $V$ is a completely reducible
$\Ve$-module. Also $V^I$ and $V^J$ are inequivalent $\Ve$-modules if
$I$ and $J$ are different. Let $U$ be any vertex operator subalgebra of $V$
 which contains
$\Ve$. Then $U$ is a direct sum of certain $V^I$. Let ${\cal D}_0$ be
the set of $I\in {\cal D}$ such that $V^I \leq  U$. 
Then $0 \in {\cal D}_0$. Also from (2) if 
$I$, $J\in{\cal D}_0$ then $I+ J\in {\cal D}_0$.
Thus ${\cal D}_0$ is a subgroup of ${\cal D}$. In order to see the simplicity
of $U$, we take a vector $v\in V^I$ for some $I\in {\cal D}_0$. Then
$\sp\{u_nv\mid u\in V^J,\, n\in \Z\}=V^{I+J}$ for any $J\in {\cal D}_0$.
It is obvious that $\{I+ J\mid J\in {\cal D}_0\}={\cal D}_0$. Thus
$U$ is simple. The proof of the irreducibility of $V^{I+ {\cal D}_0}$
is similar to that of simplicity of $V^{{\cal D}_0}$. Inequivalence of
$V^{I+ {\cal D}_0}$ and $V^{J+ {\cal D}_0}$ is clear as they are
inequivalent $T_r$-modules.

The proofs of (5) and (6) are similar to those of (2) and (3). 

For (4) we set $p=m_{h_1,\ldots,h_r}$ and $q=m_{h_1',\ldots,h_r'}$.
Let $W_1$, $\ldots$, $W_p$ be submodules of $V$ isomorphic to $M(h_1,\ldots,h_r)$
such that $\sum_{i=1}^pW_i$ is a direct sum.  Let
$d=(d_1,\ldots,d_{r})\in{\cal C}$ such that $V(d)\times
M(h_1,\ldots,h_r)=M(h_1',\ldots,h_r')$. Set $W_i'=\sp\{u_nW_i \mid u\in
V(d),\,n\in\Z\}$ for $i=1$, $\ldots$, $p$. Then $W_i'$ is isomorphic to
$M(h_1',\ldots,h_r')$ for all $i$. Note that 
\begin{eqnarray*}
& &\ \ \ \ \sp\{u_nW_i' \mid u\in
V(d),n\in\Z\}\\
& &=\sp\{u_nv_mW_i \mid u,v\in V(d),\,m,n\in\Z\}\\
& &=\sp\{u_nW_i \mid u\in T_r,n\in\Z\}=W_i
\end{eqnarray*}
(cf.~Proposition~4.1 of [DM]). 
Thus $\sum_{i=1}^pW_i'$ must be a
direct sum in $V$. This shows that $p\leq q$. Similarly, $p\geq q$. \qed

\begin{rem}\label{radd} {\rm We can also define framed vertex 
operator superalgebras. 
The analogue of Proposition~\ref{p2.3} 
still holds. In particular we have the binary codes ${\cal C}$ and ${\cal D}$.}
\end{rem}

\medskip

\begin{de} \rm Let $G$ be the subgroup of $Aut(V)$ consisting of automorphisms
which stabilize the Virasoro frame $\{\o_i\}$. Namely,
\begin{equation}\label{g1}
G=\left\{g\in Aut(V)\mid g\{\o_1,\ldots,\o_r\}=\{\o_1,\ldots,\o_r\}\,\right\}.
\end{equation}
The two subgroups $G_{\cal C}$ and $G_{\cal D}$ are defined by:
\begin{eqnarray*}
G_{\cal C}&=&\{g\in G \mid  g|_{T_r}=1\}, \\
G_{\cal D}&=&\{g\in G \mid  g|_{\Ve}=1\}.
\end{eqnarray*}
Finally, we define the automorphism group $Aut(m_h(V))$ as the subgroup of the
group $Sym_r$ of permutations of $\{1,\ldots,r\}$ which fixes the multiplicity 
function $m_h(V)$, i.e.~which consists of the permutations
$\sigma\in Sym_r$ such that $m_{h_1,\ldots,h_r}=m_{h_{\sigma(1)},\ldots,h_{\sigma(r)}}$.
\end{de}

It is easy to see that both $G_{\cal D}$ and and $G_{\cal C}$ are normal 
subgroups of $G$ and $G_{\cal D}$ is a subgroup of $G_{\cal C}$.

\smallskip
Following Miyamoto [M1], we
define for $i=1$, $\ldots$, $r$ an involution $\tau_i$ on $V$ which acts on 
$V^I$ as $-1$ if $i\in I$ and as $1$ otherwise. 
The group generated by all $\tau_i$ is a subgroup of the group 
of all automorphisms of $V$ and is isomorphic to the dual group 
$\hat{\cal D}$ of ${\cal D}$.

We  define another group $F_{\cal C}$ which is a subgroup of $Aut(\Ve)$ 
and is generated by $\sigma_i$ which acts on $M(h_1,\ldots,h_r)$ 
by $-1$ if $h_i=\frac{1}{2}$ and $1$ otherwise. 
The group $F_{\cal C}$ is isomorphic to the dual group 
$\hat{\cal C}$ of ${\cal C}$.

\smallskip

\begin{thm}\label{tg1} 

(1) The subgroup $G_{\cal D}$ is isomorphic to the dual group 
$\hat{\cal D}$ of ${\cal D}$.

(2) $G_{\cal C}/G_{\cal D}$ is isomorphic to a subgroup of the dual group 
$\hat{\cal C}$ of ${\cal C}$. 

(3) $G/G_{\cal C}$ is isomorphic to a subgroup of $Aut(m_h(V))\leq Sym_r$.
In particular, $G$ is a finite group.

(4) For any $g\in G$ and a $T_r$-submodule $W$ of $V$ isomorphic
to $M(h_1,\ldots,h_r)$ then $gW$ is isomorphic to $M(h_{\mu_g^{-1}(1)},\ldots,
h_{\mu_g^{-1}(r)})$ where $\mu_g\in Sym_r$ such that $g\omega_i=\omega_{\mu_g(i)}$ 
for all $i$. 

(5) If the eigenvalues of $g \in G_{\cal C}$ on $V^I$ are $i$ and $-i$, 
then $i$ and $-i$ have the same multiplicity. 
\end{thm}

\pf (1) Let $g\in G$ such that $g\mid_{\Ve}=1$.
Recall from Proposition~\ref{p2.3} that $V=\bigoplus_{I\in {\cal D}}V^I$. 
Since each $V^I$ is an irreducible $\Ve$-module we have
$V^I=\sp\{v_nu\mid v\in \Ve,n\in\Z\}$ for any nonzero vector 
$u\in V^I$. Note that $g$ preserves each homogeneous subspace $V^I_n,$
which is finite-dimensional. 
Take $u \in V^I$ to be an eigenvector of $g$
with eigenvalue $x_I$ and let $v \in V^0$. Then $g(v_nu)=v_ngu=x_Iv_nu$. Thus
$g$ acts on $V^I$ as the constant $x_I$. For any
$0\ne u\in V^I$ and $0\ne v\in V^J$ 
we have $0\ne Y(u,z)v\in V^{I+J}[[z,z^{-1}]]$. 
Since $x_{I+J}Y(u,z)v=gY(u,z)v=x_Ix_JY(u,z)v$, we see that $x_Ix_J=x_{I+J}$. 
In particular 
$x\in \hat{\cal D}$ and $x$ takes values in $\{ \pm 1 \}$.
Clearly,  each $g\in \hat{\cal D}$ acts on $\Ve$
trivially since $\hat{\cal D}$ is generated by the $\tau_i$. This proves (1).

For (2) we take $g\in G_{\cal C}$. A similar argument as in the first 
paragraph shows that $g|_{V(c)}$ is a constant $y_c=\pm 1$ and $y_{c+d}=y_cy_d$. 
In other words
we have defined an element $y$ of $\hat {\cal C}$ which maps $c\in {\cal C}$ 
to $y_c$. One can easily see that this gives a group homomorphism from
$G_{\cal C}$ to $\hat{\cal C}$ with kernel $G_{\cal D}$. 

For (3) let $g\in G$. Then there exists a unique
$\mu_g\in Sym_r$ such that $g\o_i=\o_{\mu_g(i)}$. Clearly we have
$\mu_{g_1g_2}=\mu_{g_1}\mu_{g_2}$ for $g_1$, $g_2\in G$. 
It is obvious that the kernel of the map $g\mapsto \mu_g$ is $G_{\cal C}$. 

In order to prove (4), we take a highest weight vector $v$ of $W$.
Then $L^i(0)v=h_iv$ for $i=1$, $\ldots$, $r$. 
So $L^i(0)gv=gL^{\mu_g^{-1}(i)}(0)v=h_{\mu_g^{-1}}v$ and
$gv$ is a highest weight vector with highest weight
$(h_{\mu_g^{-1}(1)},\ldots,h_{\mu_g^{-1}(r)})$. 
That is, $gW$ is isomorphic to $M(h_{\mu_g^{-1}(1)},\ldots,
h_{\mu_g^{-1}(r)})$. 

Finally,  we turn to (5). 
We first mention how a general $g\in G_{\cal C}$ acts on $V^I$ for $I\in {\cal D}$. 
Note that $g^2=1$ on $\Ve$ by the proof of (2), that is, $g^2\in G_{\cal D}$. So
$g^2=\pm 1$ on each $V^I$. This implies that $g$ is diagonalizable on
$V^I$ whose eigenvalues are $\pm 1$ if $g^2=1$ on $V^I$ and are $\pm i$ if
$g^2=-1$ on $V^I$. 

In the second case,
let $V^I=W_1\oplus \cdots\oplus  W_p\oplus M_1\oplus\cdots \oplus M_q$ 
where all $W_j$, $M_k$ are irreducible $T_r$-modules and $g=i$ on each $W_j$ 
and $g=-i$ on each  $M_k$. On $V^0$, 
$g$ is not $1$, since otherwise $g$ is in $G_{\cal D}$ and $g$ would have only 
$\pm 1$ for eigenvalues, by (1).  
Take  an irreducible $T_r$-submodule $U$ of $V^0$ so that $g|_U=-1$.  
Set $W_j'=\{ \ u_nW_j  \mid  u\in U, n\in\Z \}.$ Then $g=-i$ on each $W_j'$. 

{}Claim: $\sum_{j=1}^p W_j' $ is a direct sum.

Using associativity, we see that
\begin{eqnarray*}
& & \sp\{u_nW'_j \mid  u\in U,  n\in\Z\}= 
\sp\{u_mv_nW_j \mid  u,v\in U,\  m,n\in\Z\}\\
& &\hspace{2 cm}=\sp\{v_nW_j \mid  v\in T_r,\  n\in\Z \}=W_j.
\end{eqnarray*} 
This proves the claim. Thus $p\leq q$.   
Similarly, $q\leq p$. So they must be equal. This finishes the proof.
\qed

The results in Proposition~\ref{p2.3} (2) and (3) resp.~(5) and (6) 
can be interpreted by the ``quantum Galois theory'' developed in 
[DM] and [DLM2]. For example, Proposition~\ref{p2.3} (2) and (3)
is now a special case of Theorems~1 and~3 of [DM] applied
for the group $G_{\cal D}$:

\begin{rem}\label{r2.4} 
{\rm Note that $\Ve$ is the space of $G_{\cal D}$-invariants.
There is a one to one correspondence between the subgroups of $G_{\cal D}$
and vertex operator subalgebras of $V$ containing $\Ve$ via
$H\mapsto V^H$. In fact, $V^H=\bigoplus_{I\in H'}V^I$ where 
$H'=\{I\in {\cal D} \mid H|_{V^I}=1\}$. Under the identification of
$G_{\cal D}$ with $\hat{\cal D}$, the subcode $H'$ of ${\cal D}$
corresponds to the common kernel of the functionals in $H$.} 
\end{rem}

Next we prove that a FVOA is always rational. 
Recall the definition of rationality and regularity as defined in [DLM1].
A vertex operator algebra is called {\em rational\/} if any admissible module is 
a direct sum of irreducible admissible modules and 
a rational vertex operator algebra is {\em regular\/} if any weak module
is a direct sum of ordinary irreducible modules. 
(The reader is refereed to 
[DLM3] for the definitions of weak module, admissible module, and ordinary
module.)

It is proved in [DLM3] that if $V$ is a rational vertex operator algebra then $V$
has only finitely many irreducible admissible modules and each
is an ordinary irreducible module. 

We need two lemmas.
\begin{lem}\label{r1} Let $V$ be a FVOA such that ${\cal D}(V)=0$. Then any nonzero
weak $V$-module $W$ contains an ordinary irreducible module.
\end{lem}

\pf Since ${\cal D}(V)=0$ we have the decomposition $V=\bigoplus_{c\in {\cal C}}V(c)$. 
Since $T_r$ is regular (see Proposition~3.3 of [DLM1]), 
$W$ is a direct sum of ordinary irreducible $T_r$-modules.
Let $M$ be an irreducible $T_r$-submodule of $W$. Then 
$$N:=\sp\{u_nM\mid u\in V,\,n\in\Z\}$$ 
is an ordinary $V$-module as 
each $ \sp\{u_nM\mid u\in V(c),\,n\in\Z\}$ is an ordinary irreducible
$T_{r}$-module and ${\cal C}$ is a finite set. 
For an ordinary $V$-module $X$ we define
$m(X)$ to be the sum of the multiplicities $m_{h_1,\ldots,h_r}$ of all
modules $M(h_1,\ldots,h_r)$ in $X$, i.e., the $T_r$-composition length.
Let $K$ be a $V$-submodule of $N$ such that $m(K)$ is the smallest
among all nonzero $V$-submodules of $N$. Then $K$ is an irreducible ordinary
$V$-submodule of $N$ and of $W$. 
\qed

\begin{lem}\label{r2} Any FVOA $V$ with ${\cal D}(V)=0$ is rational.
\end{lem}

\pf We must show that any admissible $V$-module is a direct sum of irreducible
ones.
Let $W$ be an admissible $V$-module and $M$ the sum of all irreducible
$V$-submodules. We prove that $W=M$. Otherwise by Lemma~\ref{r1} the quotient
module $W/M$ has an irreducible submodule $W'/M$ where $W'$ is a submodule
of $W$ which contains $M$. Let $U$ be an irreducible $T_r$ submodule
of $W'$ such that $U\cap M=0$ and set $X:=\sp\{v_nU\mid v\in V,\,n\in\Z\}$. Then
$X$ is a submodule of $W'$ and $W'=M+X$. Note that 
$U[c]:=\sp\{v_nU\mid v\in V(c),\,n\in\Z\}$ for each $c\in {\cal C}$ is an 
irreducible $T_r$-module. Then either $U[c]\cap M=0$ or
$U[c]\cap M=U[c]$. If the latter happens, then $Y(v,z)(U+M/M)=0$ in the
quotient module $W/M$, which is impossible by Proposition~11.9 of~[DL].
Thus $U[c]\cap M=0$ for all $c\in {\cal C}$ and $W'=M\oplus X$.  By
Lemma~\ref{r1}, $X$ has an irreducible $V$-submodule $Y$ and certainly
$M\oplus Y$ strictly contains $M$. This is a contradiction. \qed

\begin{thm}\label{tr1} Any FVOA $V$ is rational.
\end{thm}

\pf Let $W$ be an admissible $V$-module. Then $W$ is a direct sum of
irreducible $\Ve$-modules by Lemma~\ref{r2}. Let $M$ be an irreducible 
$\Ve$-module. It is enough to show that $M$ is contained in an irreducible
$V$-submodule of $W$. 
First note that there exists a subset $I$ of $\{1,\ldots,r\}$ such that
for every irreducible $T_r$-module $M(h_1,\ldots,h_r)$ inside $M$ we have 
$h_i=\frac{1}{16}$ if and only if 
$i\in I$. Let $X$ be the $V$-submodule generated by $M$.
Then $X=\sum_{J\in{\cal D}}X[J]\leq W$ where  
$X[J]=\sp\{u_nM\mid u\in V^J,\,n\in \Z\}$
is a $\Ve$-module. We will show that $X$ is an irreducible $V$-module.

By the fusion rules, we know that for every irreducible $T_r$-submodule 
of $V$ which is 
isomorphic to 
$M(h_1,\ldots,h_r)$ has $h_k=\frac{1}{16}$ if and only if $k\in I+J$.
The $X[J]$ for $J\in{\cal D}$
are nonisomorphic $\Ve$-modules as they are nonisomorphic $T_r$-modules.
Thus $X=\bigoplus_{J\in{\cal D}}X[J]$.

Let $Y$ be a nonzero $V$-submodule of $X$. Then $Y=\bigoplus_{J\in{\cal D}}Y[J]$
where $Y[J]=Y\cap X[J]$ is a $\Ve$-module. 
If $Y[J]\ne 0$ then
$\sp\{v_nY_J\mid v\in V^J,\,n\in \Z\}\ne 0$. Otherwise  use
the associativity of vertex operators to obtain
$$0\!=\!\sp\{u_mv_nY[J] \mid u,v\in V^J,\, m,n\in \Z\}\!=\!
                   \sp\{v_nY[J]\mid v\in \Ve,\, n\in \Z\}\!=\!Y[J].$$
By associativity again we see that 
$\sp\{v_nY[J]\mid v\in V^J,\, n\in \Z\}$ is a nonzero $\Ve$-submodule of
$M$. Since $M$ is irreducible it follows immediately  that  
$\sp\{v_nY[J]\mid v\in V^J, \, n\in\Z\}=M$. So $M$ is a subspace of $Y$.
Since $X$ is generated by $M$ as a $V$-module we immediately have
$X=Y$. This shows that $X$ is indeed an irreducible $V$-module. 

It should be pointed out that each $X[J]$ in fact is an irreducible 
$V^0$-module. Let $0\ne u\in X[J]$. Since
$X=\sp\{v_nu\mid v\in V,\, n\in\Z\}$ we see that 
$\sp\{v_nu\mid v\in V^0,\, n\in\Z\}=X[J]$. \qed

\begin{cor} 
Let $V$ be a FVOA. Then

(1) $V$ has only finitely many irreducible admissible modules and every 
irreducible admissible $V$-module is an ordinary irreducible $V$-module. 

(2) $V$ is regular, that is, any weak $V$-module is a direct sum of ordinary 
irreducible $V$-modules.
\end{cor}

\pf We have already mentioned that (1) is true for all rational vertex
operator algebra (see [DLM3]). So,  (1) is an immediate consequence of 
Theorem~\ref{tr1}. In [DLM1] we proved that (2) is true
for any rational vertex operator algebra which has a regular
vertex operator subalgebra with the same Virasoro element. 
Note that $T_r$ is such a vertex operator subalgebra of $V$. \qed 

Theorem~\ref{tr1} is very useful. We will see in the later sections that the 
FVOAs $V_{\L}^+$ and $V^{\natural}$ are rational vertex operator algebras.
Theorem~\ref{tr1} 
simplifies the original proofs of the rationality of $V_{\L}^+$ in~[D3]
and $V^{\natural}$ in~[DLM1]. Most important, 
we do {\em not\/} use the self-dual property of $V^{\natural}$ 
(i.e., $V^{\natural}$ is the only irreducible module for itself) as proved in~[D3]. 

\medskip

It is a interesting problem to find suitable invariants for a FVOA $V$.
Two invariants of $V$ are the binary codes
codes ${\cal C}$ and ${\cal D}$ of length $r$ as defined before. They
cannot be arbitrary but must satisfy the following conditions:

\begin{prop}\label{codes}
(1) The code ${\cal C}$ is even, i.e.~the weight~${\rm wt}(c)=\sum_{i=1}^r
c_i\in\Z_+$ of every codeword $c\in{\cal C}$ is divisible by $2$.

(2) The weights of all codewords $d\in{\cal D}$ are divisible by $8$.

(3) The binary code ${\cal D}$ is a subcode of the annihilator code \hfill\break
${\cal C}^{\perp}=
\{ d=(d_i)\in \F_2^r\mid (d,c)=\sum_i d_i c_i=0~\hbox{for all $c=(c_i) 
\in{\cal C}$}\}$.
\end{prop}

\pf Let $W$ be a $T_r$-submodule isomorphic to $M(h_1,\ldots,h_r)$. Then the
weight of a highest weight vector of $W$ is $h_1+h_2+\cdots+h_r$ which is 
necessarily an integer as $V$ is $\Z$-graded.  
The parts (1) and (2) now follow immediately.
To see (3), note that for $c\in{\cal C}$ and 
$M\leq V^I$ isomorphic to $M(g_1,\ldots,g_r)$ 
one has from the fusion rules given in Theorem~\ref{t2.1} (2) that
$$M' = \sp\{u_nM(g_1,\ldots,g_r)\mid u\in V(c),\, n\in\Z\}\leq V^I$$
is isomorphic to $M(h_1,\ldots,h_r)$ 
with $h_i=g_i=\frac{1}{16}$ if $i\in I$ and $h_i=0$ 
(resp.~$h_i=\frac{1}{2}$) if $c_i+2g_i=0$ 
in ${\F}_2$ (resp.~$c_i+2g_i=1$).
Since the conformal weights $g_1+\cdots +g_r$ and $h_1+\cdots+h_r$
of $M(g_1,\ldots,g_r)$ and $M(h_1,\ldots,h_r)$ are both integral we see that
$\#(\{i\in \{1,\ldots,r\}\mid c_i=1\}\setminus I)$ is an even  integer.
Thus $\#\{i\in I\mid c_i=1\}$ is also even as ${\rm wt}(c)$ is even. 
This implies that $(d,c)=0$, as required,
where $d\in \cal D$ is the codeword belonging to $I\subseteq\{1,\ldots,r\}$. 
\qed

\smallskip
Here are a few remarks on the action of $Aut(V)$ on $V_2$, which is an 
action preserving the algebra product $a*b := a_3b$ coming from the VOA 
structure.  

\begin{rem}{\rm 
(1) If $V$ is a VOA and is generated as a VOA by $V_2$, then $Aut(V)$ 
acts faithfully on $V_2$.  This happens in the case  $V=V_L^+$, where $L$ 
is a lattice spanned by its vectors $x$ such that $(x,x)=4$.  

(2) If $V$ is a FVOA, the kernel of the action of $Aut(V)$ on $V_2$ is 
contained in the intersection of the groups $G_{\cal C}$, as we vary over 
all frames.  Hence, this kernel is a finite 2-group, of nilpotence class 
at most two and order dividing 
$2^r$,  where $r = \rank(V)$. }
\end{rem}

The framed vertex operator algebras with ${\cal D}=0$ can be completely 
understood in an easy way. 
\begin{prop}\label{C-voas}
For every even linear code $C\leq {\F}_2^r$ there is up to isomorphism
exactly one FVOA $V_{C}$ such that the associated binary codes are
${\cal C}=C$ and ${\cal D}=0$.
\end{prop}

\pf Let $V_{\rm Fermi}=M(0)\oplus M(\frac{1}{2})$ be the \svoa as described in
[KW]. The (graded) tensor product $V_{{\rm Fermi}}^{\otimes r}$ is a \svoa 
whose code ${\cal C}$ is the complete code $\F_2^r$ (see Remark~\ref{radd}).
It has the property, that the even vertex operator subalgebra is the 
vertex operator algebra 
associated to the level $1$ irreducible highest weight representation
for the affine Kac-Moody algebra $D_{r/2}$
if $r$ is even and $B_{(r-1)/2}$ if $r$ odd (see [H1], chapter~2). The code 
${\cal C}$ for this \voa is the even subcode of $\F_2^r$.
Proposition~\ref{p2.3}~(6) gives, for every even code $C\le \F_2^r$, 
a FVOA $V$ such that ${\cal C}(V)=C$ and ${\cal D}(V)=0$. 
The uniqueness of the FVOA with code ${\cal C}(V)=C$
up to isomorphism follows from a general result on the uniqueness of 
simple current extensions of vertex operator algebras~[H3].  \qed

This proposition is also proved in a different way by Miyamoto 
in [M2], [M3].

\medskip

Recall that a holomorphic (or self-dual) VOA is a VOA $V$ whose only 
irreducible module is $V$ itself. In the case of holomorphic FVOAs, 
we can show that the subcode 
${\cal D}\leq {\cal C}^{\perp}$ is in fact equal to 
${\cal C}^{\perp}$.

We need some basic facts from~[Z]
and~[DLM4] about the ``conformal block on the torus'' $B_V$~[Z] of a VOA 
$V$. To apply Zhu's modular invariance theorems one has to assume that 
$V$ is rational and satisfies the $C_2$ condition.\footnote{The 
condition that $V$ be a 
direct sum of highest weight representations for the Virasoro algebra
was also required in~[Z], but was removed in~[DLM4].} 
It was proved in~[DLM4] that the moonshine VOA 
satisfies the $C_2$ condition. The same proof in fact works for any
FVOA. We also know from Theorem~\ref{tr1} that a FVOA is rational.

Applying Zhu's result to a FVOA $V$ yields that $B_V$ is a finite dimensional
complex vector space with a canonical base $T_{M_i}$ indexed by the 
inequivalent irreducible $V$-modules $M_i$ and that $B_V$ carries a natural
${ SL}_2({\Z})$-module structure $\rho_V:{ SL}_2({\Z})\longrightarrow {GL}(B_V)$.

Let $V$ and $W$ be two rational VOAs satisfying the $C_2$ condition.
The following two properties of the conformal block follow directly from 
the definition:

(B1)\qquad $B_{V\otimes W}=B_V\otimes B_W$ as ${ SL}_2({\Z})$-modules
and $T_{M_i\otimes M_j}=T_{M_i}\otimes T_{M_j}$.

(B2)\qquad If $W$ is a subVOA of $V$ with the same Virasoro
element then there is a natural 
${ SL}_2({\Z})$-module map $\iota^*:B_V\longrightarrow B_W$.

We also need the following well-known result:

(B3)\qquad For the vertex operator algebra $M(0)$, the action of
$S=\left(0\ -1 \atop 1\ \ 0 \right)\in {SL}_2({\Z})$ on $B_{M(0)}$ in 
the canonical basis $\{T_{M(0)},T_{M(\ha)},T_{M(\frac{1}{16})}\}$ is given
by the matrix
\begin{equation}\label{april}
\left(\matrix{ 1/2 & 1/2 & 1/\sqrt{2} \cr 
1/2 & 1/2 & -1/\sqrt{2} \cr
1/\sqrt{2} & -1/\sqrt{2} & 0}\right)_.
\end{equation}

Here is a result about binary codes used in the proof of
Theorem~\ref{codeholo} below: 
\begin{lem}\label{mayl}
Let $\mu^{\otimes n}$ the $n$-fold tensor product of the matrix
$\mu=\left( 1\quad 1 \atop 1\ -1\right)$ 
considered as a linear endomorphism of the vector space  
${\C} [{\F}_2^n]\cong {\C} [{\F}_2]^{\otimes n}$ on the canonical
base $\{e_v\mid v\in {\F}_2^n\}$. For a subset $X\subseteq {\F}_2^n$ denote by
$\chi_X=\sum_{v\in X}e_v$ the characteristic function of $X$.
Then the following relation between a linear code $C$ and its 
annihilator $C^\perp$ holds:
$$\chi_{C^\perp}=\frac{1}{|C|}\cdot \mu^{\otimes n}(\chi_C).$$
\end{lem}

\begin{rem} {\rm $\mu^{\otimes n}$ is a Hadamard matrix of size $2^n$ and
the corresponding linear map is called the Hadamard transform.}
\end{rem}

\pf For every ${\Z}$-module $R$ and function 
$f:{\F}_2^n\longrightarrow R$
the following relation holds (cf.~Ch.~5, Lemma~2 of [MaS])
\begin{equation}\label{transform}
|C|\sum_{v \in C^\perp}f(v)=\sum_{u\in C}\sum_{v\in {\F}_2^n} (-1)^{(u,v)}f(v).
\end{equation}
Now let $R$ be the abelian group ${\C} [{\F}_2^n]$ and 
define $f$ by $f(v)=e_v$ for all $v\in {\F}_2^n$. The left hand side of
(\ref{transform}) is $|C|\cdot \chi_{C^\perp}$.
Expansion of the right side gives: 
$$\sum_{u\in C}\sum_{v_1,\ldots, v_n\in {\F}_2} \prod_{i=1}^n (-1)^{u_iv_i} e_{v_1}
\otimes\cdots\otimes e_{v_n}=\sum_{u\in C} \mu^{\otimes n}(e_u)= 
\mu^{\otimes n}(\chi_C).\quad\qed $$

\bigskip
\begin{thm}\label{codeholo}
For a holomorphic FVOA the binary codes ${\cal C}$ and ${\cal D}$
satisfy ${\cal D}={\cal C^\perp}$.
\end{thm}

\pf The vector of multiplicities $m_h(V)$ can be regarded as an
element in the vector space ${\C}[{\cal F}^r]\cong
{\C}[{\cal F}]^{\otimes r} $ where ${\cal F}=\{M(0),M(\ha),M(\frac{1}{16})\}$. 
Define two linear maps $\pi$, $\theta : {\C}[{\cal F}]
\longrightarrow {\C} [{\F}_2]={\C}\,e_0\oplus {\C}\,e_1 $ by
$$\pi(M(0))=e_0,\qquad \pi(M(\ha))=e_0, \qquad \pi(M(\frac{1}{16}))=\sqrt{2} e_1, $$
and
$$\theta(M(0))=e_0,\qquad \theta(M(\ha))=e_1, \qquad \theta(M(\frac{1}{16}))=0. $$
Finally let $\sigma:{\C}[{\cal F}]\longrightarrow {\C}[{\cal F}]$
the linear map given by the matrix (\ref{april}) relative to the
basis ${\cal F}$. 
Now one has $\pi\circ\sigma=\mu\circ\theta$ and thus the following diagram 
commutes:
\begin{equation}\label{diag}
 \qquad\matrix{
{\C}[{\cal F}^r] &  {\textstyle \sigma^{\otimes r}}    \atop \longrightarrow  & {\C}[{\cal F}^r] \cr
\cr \downarrow\,{ \theta^{\otimes r}} & &
\downarrow\,{\pi^{\otimes r}} \cr \cr
 {\C} [{\F}_2^r] & {\textstyle \mu^{\otimes r} } \atop \longrightarrow  & {\C} [{\F}_2^r].
} \qquad \qquad \qquad \qquad
\end{equation}

By definition, the support of $\pi^{\otimes r}(m_h(V))$ is 
${\cal D}\leq {\F}_2^r$. From Lemma~\ref{mayl} 
$\mu^{\otimes r}\circ\theta^{\otimes r}(m_h(V))=
|{\cal C}| \cdot \chi_{\cal C^{\perp}}\in {\C}[{\F}_2^r]$. 
Note that the support of $\chi_{\cal C^{\perp}}$ is ${\cal C^{\perp}}$.
These facts together with (\ref{diag}) imply
the theorem if we can show that $\sigma^{\otimes r}(m_h(V))=m_h(V)$.

We identify ${\C}[{\cal F}^r]$ with the conformal block on
the torus of the VOA $T_r$ by identifying the canonical
bases: $M=T_M$. Using (B1) and (B3) we observe
that $\sigma^{\otimes r}=\rho_{T_r}(S)$, where $\rho_{T_r}$ is the 
representation $\rho_{T_r}: { SL}_2({\Z})\longrightarrow { GL}(B_{T_r})$
of degree~$3^r$. 

Define the shifted graded character $ch_V(\tau):=q^{-c/24}\sum_{n\geq 0}
(\dim V_n)q^n$ where $c$ is the central charge of $V.$  
Since $V$ is holomorphic, the conformal block 
$B_V$ is one dimensional. 
Then \hbox{$\rho_V(S)=1$} (the case $\rho_V(S)=-1$ is impossible since
$ch_V(i)>0$ where $i$ is the square root of
$-1$ in upper half plane; cf.~[H1], proof of Cor.~2.1.3).
Now we use (B2).
The generator $T_V$ of $B_V$ is mapped by ${\iota}^*$ to
$\sum m_{h_1,\ldots h_r}\,T_{M(h_1,\ldots,h_r)}=m_h(V)$. Since ${\iota}^*$ is 
${ SL}_2({\Z})$-equivariant we get
$\sigma^{\otimes r}(m_h(V))=\rho_{T_r}(S)(m_h(V))={\iota}^*(\rho_V(S)(T_V))=m_h(V)$.
\qed

\smallskip
The same kind of argument was used in the proof of Theorem~4.1.5 in~[H1].

\section{Vertex operator algebras $V_{D_1^d}$}
\setcounter{equation}{0}

Let $D_n=\{(x_1,\ldots,x_n)\in \Z^n\mid \sum_{i=i}^n x_i\hbox{\ even}\}
\leq \R^n$, $n\geq 1$, be the root lattice of type $D_n$, 
the ``checkerboard lattice''. 
In this section, we describe the Virasoro decomposition of modules and
twisted modules for the vertex operator algebra $V_{D_1^d}$.

We work in the setting of [FLM] and [DMZ]. In particular $L$ is an
even lattice with nondegenerate symmetric ${\Bbb Z}$-bilinear 
form $\langle\cdot,\cdot\rangle$; ${\frak h}=L\otimes_{\Bbb Z}{\Bbb C}$;
$\hat{\frak h}_{\Bbb Z}$ is the corresponding Heisenberg algebra; $M(1)$
is the associated irreducible induced module for $\hat{\frak h}_{\Bbb Z}$ 
such that the canonical central element of $\hat{\frak h}_{\Bbb Z}$ acts as 1; 
$(\hat L,{}^-)$ 
is the central extension of $L$ by $\langle
\kappa\mid \kappa^2=1\rangle$, a group of order 2, 
 with commutator map 
$c_0(\alpha,\beta)=\<\a,\b\>+2\Z$; $c(\cdot,\cdot)$ is the alternating bilinear
form given by 
$c(\alpha,\beta)=(-1)^{c_0(\alpha,\beta)}$ for $\alpha,\beta\in L$;
$\chi$ is a faithful linear 
character of $\langle\kappa\rangle$ such that 
$\chi(\kappa)=-1$; ${\Bbb C}\{L\}=\mbox{
Ind}_{\langle\kappa\rangle}^{\hat L}{\Bbb C}_{\chi}$ ($\simeq {\Bbb C}[L]$, linearly), 
where ${\Bbb C}_{\chi}$ is the one-dimensional $\langle\kappa\rangle$-module defined by $\chi$;
$\iota(a)=a\otimes 1\in {\Bbb C}\{L\}$ for $a\in \hat L$; $V_L=M(1)\otimes 
{\Bbb C}\{L\}$; $\mbox{\bf 1}=\iota(1)$; $\omega=\frac{1}{2}\sum_{r=1}^{d}\beta_r(-1)^2$
where $\{\beta_1,\ldots,\beta_d\}$ is an orthonormal basis of ${\frak h}$;
It was proved in~[B] and [FLM] that there is a linear map 
$$\begin{array}{lcr}
V_{L}&\to& (\mbox{End}\,V_{L})[[z,z^{-1}]],\hspace*{3.6 cm} \\
v&\mapsto& Y(v,z)=\displaystyle{\sum_{n\in\Z}v_nz^{-n-1}\ \ \ (v_n\in
\mbox{End}\,V_{L})}
\end{array}$$ such that
$V_{L}=(V_{L},Y,{\bf 1},\omega)$ is a simple vertex operator algebra.
Let $L^{*}=\{x\in\frak h\mid \<x,L\>\leq \Z\}$ be the dual lattice
of $L$. Then the irreducible modules of $V_L$ are the $V_{L+\gamma}$
(which are defined in [D1]) indexed by the elements 
of the quotient group $L^{*}/L$ (see [D1]). 
In fact, $V_L$ is a rational \voa (see [DLM1]). 

Let $\theta$ be the automorphism of $\hat L$ such that 
$\theta(a)=a^{-1}\kappa^{\<\bar a,\bar a\>/2}$. Then $\theta$ is   
a lift of the $-1$ automorphism of $L$.
We have an automorphism of $V_L$,
denoted again by $\theta$, such that $\theta (u\otimes \iota(a))
=\theta(u)\otimes\iota(\theta a)$ for $u\in M(1)$ and $a\in \hat L$.
(See Appendix~\ref{append-lift} for a fuller discussion.)  
Here the action of $\theta$ on $M(1)$ is given by $\theta(\a_1(n_1)\cdots
\a_k(n_k))=(-1)^k\a_1(n_1)\cdots\a_k(n_k)$. The 
$\theta$-invariants $V_L^+$ of $V_L$ form a simple vertex operator 
subalgebra and the $-1$-eigenspace $V_L^-$ is an irreducible $V_L^+$-module
(see Theorem~2 of [DM]). Clearly $V_L=V_L^+ \oplus V_L^-$.

\medskip
Now we take for $L$ the lattice 
$$D_1^d=\bigoplus_{i=1}^d\Z\alpha_i,\  \<\a_i,\a_j\>=4\delta_{i,j}.$$ 
Then $L$ is an even lattice
and the central extension $\hat L$ is a direct product of $D_1^d$ with 
$\<\kappa\>$ and $\C\{L\}$ is simply the group algebra $\C[L]$ 
with basis 
 $e^{\a}$ for $\alpha\in L$. It is clear that
$\theta(e^{\alpha})=e^{-\alpha}$ for $\alpha\in D_1^d$. 
We extend the action of $\theta$ from $V_{D_1^d}$ to $V_{(D_1^*)^d}=M(1)\otimes 
\C[L^*]$ such that $\theta (u\otimes e^{\alpha})=(\theta u)\otimes e^{-\alpha}$
for $u\in M(1)$ and $\alpha\in L^*$. 
One can easily verify that $\theta$ has order $2$ and 
$\theta Y(u,z)\theta^{-1}=Y(\theta u,z)$ for $u\in V_{D_1^d},$
where $Y(v,z)$ $(v\in V_{D_1^d})$ are the vertex operators on $V_{(D_1^*)^d}.$
For any $\theta$-invariant subspace $V$ of $V_{L^*}$ we use $V^{\pm}$ to denote
the $\pm$-eigenspaces. 

First we turn our attention to the case that $d=1$. 
Then $L=\Z\alpha\cong 2\Z=D_1$ where $\<\a,\a\>=4$. 
Note that the dual lattice $D_1^*$ is $\frac{1}{4}D_1$ and 
$\{0,1,\ha,-\ha\}$ is a system of coset 
representatives of $D_1^{*}/D_1$. 

Set 
\begin{eqnarray}\label{twostar}
& &\omega_{1} =\frac{1}{16}\alpha(-1)^2
+\frac{1}{4}(e^{\alpha}+e^{-\alpha}),\nonumber \\
& &\omega_{2}=\frac{1}{16}\alpha(-1)^2
- -\frac{1}{4}(e^{\alpha}+e^{-\alpha}).
\end{eqnarray}
Then $\omega_i\in V_{D_1}^+$. 
\begin{lem}\label{l3.1}   For $D_1\cong L = \Z \alpha$, $\langle \alpha, 
\alpha \rangle = 4$, we have: 

(1) $V_{D_1}$ is a FVOA with $r=2$.

(2) We have the following Virasoro decompositions of $V_{D_1}^+$ and 
$V_{D_1}^-:$
$$V_{D_1}^+ \cong M(0,0), 
\ \ V_{D_1}^- \cong M(\frac{1}{2}, \frac{1}{2})$$
with highest weight vectors ${\bf 1}$ and $\a(-1)$, respectively.

(3) The decompositions for $V_{D_1+1}^{\pm}$ are:
$$V_{D_1+1}^+ \cong M(\ha,0), \ \ 
V_{D_1+1}^-\cong M(0,\ha)$$
with highest weight vectors
$(e^{\frac{1}{2}\a}-e^{-\frac{1}{2}\a})$ and 
$(e^{\frac{1}{2}\a}+e^{-\frac{1}{2}\a})$, respectively. 

(4) For $V_{D_1+\frac{1}{2}}\oplus V_{D_1-\frac{1}{2}}$ 
we get, in both cases,
  $$(V_{D_1+\frac{1}{2}}\oplus V_{D_1-\frac{1}{2}})^{\pm} \cong 
M(\frac{1}{16},\frac{1}{16})$$
with highest weight vectors $e^{\frac{1}{4}\a}\pm e^{-\frac{1}{4}\a}$. In fact, both
$ V_{D_1+\frac{1}{2}}$ and $V_{D_1-\frac{1}{2}}$ 
are irreducible $V_{D_1}^+$-modules.
\end{lem}

\pf
It was proved in [DMZ] (see Theorem~6.3 there) that 
$Y(\o_1,z_1)=\sum_{n\in\Z}L^1(n)z^{-n-2}$ and 
$Y(\o_2,z_2)=\sum_{n\in\Z}L^2(n)z^{-n-2}$ give 
two commuting Virasoro algebras with central charge~$\ha$.  
We first show that $\a(-1)$ is highest weight
vector with highest weight  $(\frac{1}{2},\frac{1}{2})$
for the two Virasoro algebras. 
Since $\a(-1)\in V_{D_1}^-$ has the smallest weight in $V_{D_1}^-$ it
is immediate to see that $L^i(n)\a(-1)=0$ if $n>0$. 
It is a straightforward computation by using the definition of vertex operators
to show that $L^1(0)\alpha(-1)=L^2(0)\a(-1)={1\over 2}\a(-1)$. 

Clearly, ${\bf 1}\in(V_{D_1})^+$ is a highest weight vector for the 
Virasoro algebras with highest weight $(0,0)$. So $V_{D_1}$ contains
two highest weight modules for the two Virasoro algebras with highest
weights $(0,0)$ and $(\frac{1}{2},\frac{1}{2})$. 
Since $M(0,0)\oplus M(\frac{1}{2},\frac{1}{2})$ and
$V_{D_1}$ have the same graded dimension we conclude that
$V_{D_1} \cong M(0,0)\oplus M(\frac{1}{2},\frac{1}{2})$ 
and $V_{D_1}^+ \cong M(0,0)$, $V_{D_1}^- \cong 
M(\frac{1}{2},\frac{1}{2})$.
This proves (2) and shows also (1): $V_{D_1}$ is a FVOA
with $r=2.$ Additionally we see that $V_{D_1}$ is a unitary representation 
of the two Virasoro algebras.

By Theorem~\ref{t2.1} (3) we know that $V_{D_1+\l}$,  for $\l=0$, 
$\pm \frac{1}{2}$, $1$,  
is a direct sum of irreducible modules $M(h_1,h_2)$
with $h_i\in \{0,\frac{1}{2},\frac{1}{16} \}$. It is easy to
find all highest weight vectors in $V_{D_1+\l}$. 
Part (3) and (4) follow immediately then.  
\qed

We return to the lattice
$L=\bigoplus_{i=1}^d\Z\a_i$, $\<\a_i,\a_j\>=4\delta_{i,j}$, 
$L\cong D_1^d=(2\Z)^d$. 
We sometimes identify $L$ with $(2\Z)^d$.
The component $\Z\a_i$ gives two Virasoro elements $\o_{2i-1}$
and $\o_{2i}$, as in (\ref{twostar}), above.

\begin{de}\label{newde}\rm   
The VF associated to the FVOAs derived 
{}from the $D_1^d$-lattice is the set $\{\o_1,\ldots,\o_{2d}\}$.
\end{de}

\begin{cor}\label{c3.1} (1) The decomposition of $V_{D_1^d}^{\pm}$ into 
irreducible modules for $T_{2d}$ is given by
$$V_{D_1^d}^{\pm}\cong \bigoplus_{\hbox{\scriptsize $\begin{array}{c}
(h_{2i-1},h_{2i})\in \{(0,0),(\frac{1}{2},\frac{1}{2})\} \\
(-1)^{\#\{i\mid h_{2i}=0\}}=\pm 1
\end{array}$ } }
M(h_1,\ldots,h_{2d}).$$
In particular,
$V_{D_1^d}^{\pm}$ is a direct sum of $2^{d-1}$ irreducible modules for 
$T_{2d}$.

(2) Let $\gamma=(\gamma_i)\in (D_1^*)^d$ such that $\gamma_i\in\{0,1\}$.
Then we get the decomposition
$$(V_{D_1^d+\gamma})^{\pm}\cong \bigoplus_{(h_{2i-1},h_{2i}) \in \left\{
\hbox{ \scriptsize $\begin{array}{ll}
\{(0,0),(\frac{1}{2},\frac{1}{2})\} & \hbox{if\ } \gamma_i=0, \\
\{(0,\frac{1}{2}),(\frac{1}{2},0)\} & \hbox{if\ } \gamma_i=1 \\
\end{array}$ }\right.  \atop
(-1)^{\#\{i \mid h_{2i}=0\}}=\pm1}
M(h_1,\ldots,h_{2d}).$$

(3) Let $\gamma=(\gamma_i)\in (D_1^{*})^d$, such that 
$2\gamma\not\in D_1^d$, i.e.~there is at least one $i$ such that
 $\gamma_i=\pm\frac{1}{2}$.
Then $(V_{D_1^d+\gamma}\oplus V_{D_1^d-\gamma})^{\pm}$, 
$V_{D_1^d\pm\gamma}$ have the same decomposition:
$$\bigoplus_{(h_{2i-1},h_{2i})\in \left\{
\hbox{ \scriptsize $\begin{array}{ll}
\{(0,0),(\frac{1}{2},\frac{1}{2})\} & \hbox{if\ }\gamma_i=0, \\
\{(\frac{1}{2},0),(0,\frac{1}{2})\} & \hbox{if\ }\gamma_i=1, \\
\{(\frac{1}{16},\frac{1}{16})\}        
& \hbox{if\ }\gamma_i=\pm\frac{1}{2}
\end{array}$ }\right. }
M(h_1,\ldots,h_{2d}).$$
\end{cor}

\pf  Note that $V_{D_1^d}$ is isomorphic to the tensor product 
vertex operator algebra
$V_{D_1}\otimes\cdots\otimes V_{D_1}$ ($d$ factors) and 
that
 $V_{D_1^d+\gamma}$ is isomorphic to the tensor product module 
$V_{\Z\a_1+\gamma_1}\otimes\cdots\otimes V_{\Z\a_d+\gamma_d}$. Thus 
$$(V_{D_1^d+\gamma}\oplus V_{D_1^d-\gamma})^{\pm}=
\bigoplus_{\stackrel{\mu\in\{+,-\}^d}{\prod\mu_i=\pm}}
V_{D_1+\gamma_1}^{\mu_1}\otimes\cdots\otimes V_{D_1+\gamma_d}^{\mu_d}.$$
The results (1) and (2) now follow from Lemma~\ref{l3.1} immediately.  

For (3) it is clear that the decompositions for $V_{D_1^d\pm\gamma}$
hold by Lemma~\ref{l3.1}. It remains to show that $V_{D_1^d\pm\gamma}$
and $(V_{D_1^d+\gamma}\oplus V_{D_1^d-\gamma})^{\pm}$ are all isomorphic
$T_{2d}$-modules. 
Note from Lemma~\ref{l3.1} that $V_{D_1+h}$ and
$V_{D_1-h}$ are isomorphic $T_2$-modules for any $h\in \{0,1,\pm\ha\}$.
Thus $V_{D_1^d+\gamma}$ and $V_{D_1^d-\gamma}$ are isomorphic 
$T_{2d}$-modules. In fact, $\theta: V_{D_1^d+\gamma}\to V_{D_1^d-\gamma}$ 
is such an isomorphism. Thus, 
$(V_{D_1^d+\gamma}\oplus V_{D_1^d-\gamma})^{\pm}=
 \{v\pm \theta v\mid v\in V_{D_1^d+\gamma}\}$ are isomorphic to 
$V_{D_1^d+\gamma}$ as $T_{2d}$-modules.
\qed  

\bigskip
Next we discuss the twisted modules of $V_L$ for an arbitrary 
$d$-dimensional positive definite even lattice $L$.
Recall from [FLM] the definition of the 
twisted sectors associated to an even lattice
$L$. Let $K=\{\theta(a)a^{-1}\mid a\in \hat L\}$. Then $\bar K=2L$ (bar
is the quotient map $\hat L \rightarrow L$). Also
set $R:=\{\alpha\in L\mid \langle\alpha,L\rangle \le 2\Z \}$; then 
$R \geq  2L$.
Then the inverse image $\hat R$ of $R$ in $\hat L$ is the center of $\hat
L$ and $K$ is a subgroup of $\hat R$. It was proved in [FLM]
(Proposition 7.4.8) there are exactly $|R/2L|$ central characters
$\chi: \hat{R}/K\to \C^{\times}$ of $\hat L/K$ such that
$\chi(\k K)=-1$. For each such $\chi$, there is a unique (up to
equivalence) irreducible $\hat L/K$-module $T_{\chi}$ with central
character $\chi$ and every irreducible $\hat L/K$-module on which
$\k K$ acts as $-1$ is equivalent to one of these. In particular,
viewing $T_{\chi}$ as $\hat L$-module, $\theta a$ and $a$ agree 
as operators on $T_{\chi}$ for
$a\in \hat L$. Let $\hat {\frak h}[-1]$ be the twisted Heisenberg
algebra. As in Section 1.7 of [FLM] we also denote by $M(1)$
the unique irreducible $\hat{\frak h}[-1]$-module with
the canonical central element acting by $1.$
Define the twisted space $V_{L}^{T_{\chi}}=M(1)\otimes T_{\chi}$.  It was
shown in [FLM] and [DL2] that there is a linear map
\begin{eqnarray*}
V_{L}&\to& (\mbox{End}\,V_{L}^{T_{\chi}})[[z^{1/2},z^{-1/2}]],\hspace*{3.6 cm} \\
v&\mapsto& Y(v,z)=\displaystyle{\sum_{n\in\frac{1}{2}\Z}v_nz^{-n-1}}\ \
\end{eqnarray*}
such that $V_{L}^{T_{\chi}}$ is an irreducible
$\theta$-twisted module for $V_L$. Moreover, 
every irreducible $\theta$-twisted $V_L$-module
is isomorphic to  $V_{L}^{T_{\chi}}$ for some $\chi$.

Define a linear operator $\hat\theta_d$ on $V_L^{T_{\chi}}$ such that
$$\hat\theta(\alpha_1(-n_1)\cdots\a_k(-n_k)\otimes t)=(-1)^ke^{d\pi i/8}
\alpha_1(-n_1)\cdots\a_k(-n_k)\otimes t$$
for $\a_i\in \frak h$, $n_i\in\frac{1}{2}+\Z$ and $t\in T$.
Then $\hat\theta_d Y(u,z)(\hat \theta_d)^{-1}=Y(\theta u,z)$ for $u\in V_L$
(cf.~[FLM]).  
We have the decomposition $V_L^{T_{\chi}}=
(V_L^{T_{\chi}})^+\oplus(V_L^{T_{\chi}})^-$ where $(V_L^{T_{\chi}})^+$ 
and $(V_L^{T_{\chi}})^-$ are the $\hat\theta_d$-eigenspaces with 
eigenvalues $-e^{d\pi i/8}$ and $e^{d\pi i/8}$ respectively. Then 
both $(V_L^{T_{\chi}})^+$ and $(V_L^{T_{\chi}})^-$ are irreducible 
$V_L^+$-modules (cf.~Theorem~5.5 of~[DLi]). 

\medskip
As before, we now take $L=\Z\a\cong D_1$ with $\<\a,\a\>=4$. Then $K=2L$,
$R=L$ and $R/K \cong \Z_2$. 
Let $\chi_1$ be the trivial character of
$R/K$ and $\chi_{-1}$ the nontrivial character. Then both $T_{\chi_1}$
and $T_{\chi_{-1}}$ are one-dimensional $L$ modules and
$\a$ acts on $T_{\chi_{\pm 1}}$ as $\pm1$.

\begin{lem}\label{l3.3} We have the Virasoro decompositions:
$$ (1)\qquad (V_{D_1}^{T_{\chi_1}})^+\cong M(\frac{1}{16},\frac{1}{2}), 
\qquad 
             (V_{D_1}^{T_{\chi_1}})^-\cong M(\frac{1}{16},0).$$
$$ (2)\qquad (V_{D_1}^{T_{\chi_{-1}}})^+\cong M(\frac{1}{2},\frac{1}{16}),\qquad
             (V_{D_1}^{T_{\chi_{-1}}})^-\cong M(0,\frac{1}{16}).$$
\end{lem}

\pf Recall from~[DL2] that 
$$V_L^{T_{\chi_1}}=\sum_{n\in{1\over 2}\Z,
\ n\geq 0}(V_L^{T_{\chi_1}})_{\frac{1}{16}+n}$$
(see Proposition 6.3 and formula (6.28) of~[DL2]). Note that
$$(V_L^{T_{\chi_1}})^+=\sum_{n\in \Z,\ n\ge0}
(V_L^{T_{\chi_1}})_{\frac{1}{16}+\ha+n}$$ 
and that
$$(V_L^{T_{\chi_1}})^-=\sum_{n\in \Z,\ n\ge0}(V_L^{T_{\chi_1}})_{\frac{1}{16}+n}.$$
Since both
$(V_L^{T_{\chi_1}})^+$ and $(V_L^{T_{\chi_1}})^-$ are irreducible $V_L^+$-modules
we only need to calculate highest weights for nonzero highest weight 
vectors in these spaces. Note that $T_{\chi_1}$ is a space of highest weight
vectors of $(V_L^{T_{\chi_1}})^-$. One can easily verify that 
$L^1(0)=\frac{1}{16}$ and $L^2(0)=0$ on $T_{\chi_1}$. 
Thus $(V_L^{T_{\chi_1}})^-\cong M(\frac{1}{16},0)$.

Also observe
that $\a(-\frac{1}{2})\otimes T_{\chi_1}$ is a space of highest weight vectors
of $(V_L^{T_{\chi_1}})^+$. From Lemma~\ref{l3.1} we know that $\a(-1)\in 
V_L^-\cong M(\frac{1}{2},\frac{1}{2})$.  
Now use the fusion rule given in Theorem \ref{t2.1}
to conclude that $(V_L^{T_{\chi_1}})^+\cong M(\frac{1}{16},\frac{1}{2})$. 
Part (2) is proved similarly.
\qed

\medskip
As we did in the untwisted case, we now consider 
the twisted modules for the lattice
$L=\bigoplus_{i=1}^d\Z\a_i\cong D_1^d$, $\<\a_i,\a_j\>=4\delta_{i,j}$, 
where $d$ is now a positive integer divisible by $8$. 
Then $K=2L$, $R=L$ and $R/2L\cong  \Z_2^{d}$. Thus,  there are
$2^d$ irreducible characters for $R/2L$ which are denoted by $\chi_J$
(where $J$ is a subset of $\{1,\ldots,d\}$) sending $\a_j$ to $-1$ if
$j\in J$ and to $1$ otherwise. Then we have $\chi_J=\prod_{j}\chi_{x_j}$
where $\chi_{x_j}$ is a character of $\Z\a_j/\Z2\a_j$ 
and $x_j=\chi_J(a_j)$. Moreover,  
$T_{\chi_J}\cong T_{\chi_{x_1}}\otimes\cdots \otimes T_{\chi_{x_d}}$.
In particular, each $T_{\chi_J}$ is one dimensional. 

\begin{cor}\label{c3.4} We have the Virasoro decompositions:
$$(V_{D_1^d}^{T_{\chi_J}})^{\pm}=\bigoplus_{\hbox{ \scriptsize
$(h_{2i-i},h_{2i})\in \left\{
\begin{array}{c}
\{(\frac{1}{16},0),(\frac{1}{16},\frac{1}{2})\}\ \hbox{if}\ i\not\in J\\
\{(0,\frac{1}{16}),(\frac{1}{2},\frac{1}{16})\}\ \hbox{if}\ i\in J
\end{array}\right.$}  \atop
(-1)^{\#\{j\mid h_{j}=\frac{1}{2}\}}=\pm (-1)^{d/8}}
M(h_1,\ldots,h_{2d}).$$ 
\end{cor}

\pf Recall from the proof of Corollary~\ref{c3.1} that 
$V_{D_1^d}$ is isomorphic to the tensor product vertex operator algebra
$V_{D_1}\otimes\cdots\otimes V_{D_1}$. Note that $V_{D_1^d}^{T_{\chi_J}}$ 
is isomorphic to the tensor product
$V_{D_1}^{T_{\chi_{x_1}}}\otimes\cdots\otimes 
V_{D_1}^{T_{\chi_{x_d}}}$ and $\hat\theta_d$ is also a tensor product
$\hat\theta_1\otimes\cdots \otimes\hat\theta_d$. 
By Lemma \ref{l3.3},
$$V_{D_1^d}^{T_{\chi_J}}
=\bigoplus_{\hbox{ \scriptsize
$(h_{2i-i},h_{2i})\in \left\{
\begin{array}{c}
\{(\frac{1}{16},0),(\frac{1}{16},\frac{1}{2})\}\ \hbox{if}\ i\not\in J\\
\{(0,\frac{1}{16}),(\frac{1}{2},\frac{1}{16})\}\ \hbox{if}\ i\in J
\end{array}\right.$}}M(h_1,\ldots,h_{2d}).$$ 
Since $\hat\theta_d=(-1)^{\#\{j\mid h_{j}=\frac{1}{2}\}}(-1)^{d/8}=(-1)^{\#\{j\mid h_{j}=0\}}(-1)^{d/8}$
on $M(h_1,\ldots,h_{2d})$ we see that $M(h_1,\ldots,h_{2d})$ 
embeds in  $(V_L^{T_{\chi_J}})^{\pm}$ 
if and only if $(-1)^{\#\{j\mid h_{j}=\frac{1}{2}\}}(-1)^{d/8}=\pm 1$. The
proof is complete.
\qed 

\begin{rem} {\rm Note that $V_{D_1^d}^{T_{\chi_J}}$ is $\frac{1}{2}\Z$ graded if 
$d$ is divisible by $8$ (cf.~[DL2]). In fact $(V_{D_1^d}^{T_{\chi_J}})^+$ is then the
subspace of $V_{D_1^d}^{T_{\chi_J}}$ consisting of vectors of integral weights
while $(V_{D_1^d}^{T_{\chi_J}})^-$  is the subspace of $V_{D_1^d}^{T_{\chi_J}}$ 
consisting of vectors of non-integral weights.}
\end{rem}

\section{Vertex operator algebras associated to binary codes}
\setcounter{equation}{0}

Let $C$ be a doubly-even linear binary code of length 
$d\in 8\Z$ containing the all ones vector $\Omega=(1,\ldots,1)$.
As mentioned in Section 2, we can regard a vector of $\F_2^d$
as an element in $\Z^d$ in an obvious way.   
One can associate (cf.~[CS1]) to such a code the two even lattices
$$\lc=\{\frac{1}{\sqrt{2}}(c+x)\mid c\in C,\ x\in (2\Z)^d\} \qquad\qquad  $$
and
$$\tlc= \{\frac{1}{\sqrt{2}}(c+y)\mid c\in C,
\ y\in (2\Z)^d,\ \hbox{$4\mid\sum y_i$}\}\ \cup\qquad\qquad $$
$$\qquad\qquad
\{\frac{1}{\sqrt{2}}(c+y+
\hbox{$(\frac{1}{2},\ldots,\frac{1}{2})$})
\mid c\in C,
\ y\in (2\Z)^d,\ \hbox{$4\mid (1-(-1)^{d/8}+\sum y_i)$}\} $$
and for every {\em self-dual\/} even lattice there are two 
vertex operator algebras $V_L$ and $\tv=V_L^+\oplus (V_L^T)^+$ (see~[FLM], [DGM1]).

\begin{de} \rm
A {\em marking\/} for the code $C$ is a partition
${\cal M}=\{(i_1,i_2),\ldots,(i_{d-1},i_d)\}$ of the positions
$1$, $2$, $\ldots$, $d$ into ${d\over 2}$ pairs. 
\end{de}

A marking ${\cal M}=\{(i_1,i_2),\ldots,(i_{d-1},i_d)\}$ determines
the $D_1^d$ sublattice $\bigoplus_{l=1}^d \Z\alpha_l$ inside 
$L_C$ and $\tlc$, where $\alpha_{2k-1}=\sqrt{2}(e_{i_{2k-1}}
+e_{i_{2k}})$ and $\alpha_{2k}=\sqrt{2}(e_{i_{2k-1}}-e_{i_{2k}})$
for $k=1$, $\ldots$, ${d\over 2}$ using $\{e_i\}$ as the standard base of 
$\lc\otimes\Q=\Q^{\Omega}$. Let us simplify notation and arrange for the 
marking to be ${\cal M}=\{(1,2),(3,4),\ldots,(d-1,d)\}.$

{}From Definition \ref{newde}, we see that every such marking defines a 
system of $2d$ commuting Virasoro algebras inside the 
vertex operator algebras $\vlc$, $\vtlc\cong\tvlc$ and $\tvtlc$.
As the main theorem we describe explicitly the
decomposition into irreducible $T_{2d}$-modules in terms of the marked
code. The triality symmetry of $\tvtlc$ given in~[FLM] and~[DGM1]
is directly visible in this decomposition. (See also [G1].) 

In order to give the Virasoro decompositions in a readable way,  we 
need the next lemma which describes $\lc$ and $\tlc$ in the coordinate system 
spanned by the $\alpha_{i}$. We use the following notation.  
Let $\gamma_+^0$, $\gamma_-^0$, $\gamma_+^1$
and $\gamma_-^1$ be the maps $\F_2^2\longrightarrow 
(D_1^{*}/D_1)^2=\{0,\frac{1}{2},-\frac{1}{2},1\}^2$ defined by
the table
$$
\begin{array}{c||cccc}
 & (0,0) & (1,1) & (1,0) & (0,1) \\ \hline \hline
\gamma_+^0 & (0,0) & (1,0) & (\ha,\ha) & (\ha,-\ha) \\
\gamma_-^0 & (1,1) & (0,1) & (-\ha,-\ha) & (-\ha,\ha) \\ \hline 
\gamma_+^1 & (\ha,0) & (-\ha,0) & (1,\ha) & (1,-\ha) \\
\gamma_-^1 & (-\ha,1) & (\ha,1) & (0,-\ha) & (0,\ha) 
\end{array}
$$
and write $c(k)$ ($k=1$, $\ldots$, ${d\over 2}$) for the pair $(c_{2k-1},c_{2k})$ of 
coordinates of a codeword $c\in C$.  
Finally let for $b=0$ or $1$ and 
$\epsilon=(\epsilon_1,\ldots,\epsilon_{d/2})\in \{+,-\}^{d/2}$
\begin{equation}\label{dgh-1}
\Gamma_{\epsilon}^{a}(b)=\bigoplus_{k=1}^{d/2}
D_1^2+\gamma_{\epsilon_k}^b(c(k))
\end{equation}
be a coset of $D_1^d$. 
\begin{lem}\label{ldecomp}
We have the following coset
decomposition of $\lc$ and $\tlc$ under the above $D_1^d$ sublattice:
\begin{eqnarray*}
\lc&=&\bigcup_{c\in C}
\bigcup_{ \epsilon\in\{+,-\}^{d/2}}
\Gamma_{\epsilon}^{0}(c),
\\
\tlc & = & \bigcup_{c\in C}\left(
\bigcup_{ \epsilon\in\{+,-\}^{d/2} \atop \prod \epsilon_i=+ }
\Gamma_{\epsilon}^{0}(c)
\cup
\bigcup_{ {\epsilon\in\{+,-\}^{d/2}} \atop  
\prod \epsilon_i=(-)^{d/8}}
\Gamma_{\epsilon}^{1}(c)
\right).
\end{eqnarray*}
\end{lem}

\pf The result follows from the definition of these two lattices.
\qed

\smallskip
We next interpret the decompositions in terms of codes over 
$\Z_4=\{0,1,2,3\}$ associated to $L_C$ and $\tilde L_C$. See~[CS2] for 
the relevant definitions for $\Z_4$-codes.

Let $L$ be a positive definite even lattice of rank $d$ which contains
a $D_1^d$ as a sublattice. We call such a sublattice a
{\em $D_1$-frame\/}. Note that $(D_1^*/D_1)^d$ is isomorphic
to $\Z_4^d$. Then $\Delta(L):=L/D_1^d\leq (D_1^*/D_1)^d$ is 
a code over $\Z_4$ and $\Delta(L)$ is self-annihilating if and only if $L$ is
self-dual. For the lattices $L_C$ and $\tilde L_C$ we give the following 
explicit description of the corresponding codes $\Delta$:

Let $\widehat{\phantom{X}}$ be the map from $\F_2^d$ to $\Z_4^d$ 
induced from $\hat{\phantom{.}}:
\F_2^2\cong D_2^*/D_2 \longrightarrow (D_1^*/D_1)^2\cong \Z_4^2$, 
$00\mapsto 00$, $11\mapsto 20$, $10\mapsto 11$ and $01\mapsto 31$.
Let $(\Sigma_2^n)_0$ be the subcode of the $\Z_4$-code 
$\Sigma_2^n=\{(00),(22)\}^n$ of length $2n$ consisting of codewords of weights 
divisible by $4$. 
Then we have
\begin{eqnarray}\label{gamma-def}
\Gamma &\!\!=\!\! &\Delta(L_C)=\widehat{C}+\Sigma_2^{d/2}, \\
\widetilde{\Gamma}&\!\!= \!\!&\Delta(\tilde L_C)\!=\!
\widehat{C}\!+\!(\Sigma_2^{d/2})_0\cup\widehat{C}\!+
\!(\Sigma_2^{d/2})_0\!+\!
\cases{ (1,0,\ldots,1,0,1,0),\! & if $d\equiv 0\pmod{16}$, \cr
        (1,0,\ldots,1,0,3,2),\! & if $d\equiv 8\pmod{16}$.}\nonumber 
\end{eqnarray}

\smallskip
An important invariant of a $\Z_4$-code $\Delta$ is the   
symmetrized weight enumerator, as defined in~[CS2].
\begin{de}\rm 
The {\em symmetrized weight enumerator\/}  
of a $\Z_4$-code $\Delta$ of length $d$ is given by
$${\rm swe}_{\Delta}(A,B,C)=\sum_{0\leq r,s\leq d}U_{r,s}\,A^{d-r-s}B^rC^s,$$
where $U_{r,s}$ is the number of codewords $\gamma\in \Delta$ having
at $r$ positions the value $\pm 1$ and at $s$ positions the value $2$.
\end{de}

To describe the symmetrized weight enumerator  
for our codes $\Gamma$ and $\widetilde{\Gamma}$ in terms of the marked 
binary code $C$, we introduce the analogous invariant for
marked binary codes. 
\begin{de} \rm
The {\em symmetrized marked weight enumerator\/} of a binary code $C$ of 
length $d$ with a marking ${\cal M}$ is the homogeneous polynomial
$${\rm smwe}_C(x,y,z)=
       \sum_{k=0}^{d/2}\sum_{l=0}^{[k/2]}W_{k,l}\,x^{d/2-k+l}y^{k-2l}z^l,$$
where $W_{k,l}$ is the number of codewords $c\in C$ of Hamming weight $k$
having the value $(c_{i_{2m-1}},c_{i_{2m}})=(1,1)$ 
for exactly $l$ of the $d/2$ pairs $(i_{2m-1},i_{2m})$ of the marking ${\cal M}$. 
\end{de}

{}From Lemma~\ref{ldecomp} we get:
\begin{cor}\label{ldecomp-cor}
The symmetrized weight enumerators of the $\Z_4$-codes $\Gamma$ and 
$\widetilde{\Gamma}$ are given by
\begin{eqnarray}
{\rm swe}_{\Gamma}(A,B,C) & = & {\rm smwe}_C(A^2+C^2,2B^2,2AC), \\
{\rm swe}_{\widetilde{\Gamma}}(A,B,C) & = & 
\frac{1}{2}\cdot{\rm smwe}_C(A^2+C^2,2B^2,2AC)  
+\frac{1}{2}\left((A^2-C^2)^{d/2}\right) \\ &&\qquad\qquad
+\frac{1}{2}\cdot 2^{d/2}\left((A+C)^{d/2}+(-1)^{d/8}(A-C)^{d/2}\right)
B^{d/2}. \nonumber\qquad \qed
\end{eqnarray}
\end{cor}

\medskip

Motivated by Lemmas~\ref{l3.1} and~\ref{l3.3},
define for $a \in \{0,1\}$, $\alpha\in\{+,-\}$ and 
$x\in\{0,\pm\ha,1\}$ the $16$ formal linear combinations
of $T_2$-modules $R_{\alpha}^{a}(x)$ by the following table:
$$\begin{array}{c||cccc}
  & 0 & 1 & \ha,-\ha \\ \hline\hline
R_{+}^{0} & M(0,0) & M(\ha,0)  &   \ha   M(\se,\se) \\
R_{-}^{0} & M(\ha,\ha) & M(0,\ha)  & \ha    M(\se,\se)  \\  \hline
R_{+}^{1} & \multicolumn{2}{c}{\frac{1}{\sqrt{2}}M(\se,\ha)} &
            \frac{1}{\sqrt{2}} M(\ha,\se)\\
R_{-}^{1} & \multicolumn{2}{c}{\frac{1}{\sqrt{2}}M(\se,0)} &   
             \frac{1}{\sqrt{2}}M(0,\se)
\end{array}$$
For $\mu\in\{+,-\}^{n}$ and an element $\gamma=(\gamma_k)\in\Z_4^{n}$ which
is identified with $\{0,1,\pm\frac{1}{2}\}$ we write 
shortly  
\begin{equation}\label{dgh0}
{\bf R}^{a}_{\mu}(\gamma)=\bigotimes_{k=1}^{n}R_{\mu_k}^{a}(\gamma_k).
\end{equation}

We see from Lemmas~\ref{l3.1} and~\ref{l3.3}  
that $R_{\pm}^0=V_{D_1+x}^{\pm}$
and $R_{\pm}^1=\frac{1}{\sqrt{2}}(V_{D_1}^{T_{\chi_{(-1)^{2x}}}})^{\pm}$.
The introduction of the extra factor $\frac{1}{\sqrt{2}}$ in the twisted
case enables us to write the Virasoro decompositions for 
the twisted sector $V_L^T$ in a neat way. 
The index $0$ in $R^0_{\pm}$ refers to the untwisted case while
the index $1$ in $R^1_{\pm}$ refers to the twisted case. 

Let $L$ be a self-dual even lattice of rank $d$ containing a $D_1$-frame.
So, $L$ is defined by the self-annihilating $\Z_4$-code $\Delta=
L/D_1^d\leq (D_1^*/D_1)^d$ of length $d$ which is now even in the sense 
that ${\rm swe}_{\Delta}(1,x^{\frac{1}{2}},x)$ is a polynomial in $x^2$ (cf.~[BS]).

\begin{thm}\label{d1-decomp}
The vertex operator algebras $V_L$ 
and $\widetilde{V}_L$ have the following decompositions as modules
for $T_{2d}$:
\begin{eqnarray*}
V_L &=& \bigoplus_{\gamma\in\Delta}
\bigoplus_{ {\mu\in\{+,-\}^{d}} }{\bf R}_{\mu}^{0}(\gamma),
\\ 
{\widetilde V}_{L} &= & 
\bigoplus_{\gamma\in\Delta}
\bigoplus_{ {\mu\in\{+,-\}^{d}} \atop \prod\mu_k=+}
{\bf R}_{\mu}^{0}(\gamma)
\oplus\bigoplus_{ \gamma\in\Delta}
\bigoplus_{ {\mu\in\{+,-\}^{d}}\atop \prod\mu_k=(-)^{d/8}} 
{\bf R}_{\mu}^{1}(\gamma).
\end{eqnarray*}
\end{thm}

In order to determine the decomposition for ${\widetilde V}_{L}$, we 
first study the decomposition of $V_L^T$.

Since $L$ is self-dual, $V_L$ has a unique irreducible $\theta$-twisted
module $V_L^T$~[D2]. In this case $T$ can be constructed in the following
way. Let $Q$ be a subgroup of $L$ containing the $D_1^{d}$ 
which is maximal such that $\<\alpha,\beta\>\in2\Z$ for $\a$, $\b\in Q$ 
(it exists since $L$ has ascending chain conditions on subgroups).
Let  $\hat Q$ be  the inverse image 
of $Q$ in $\hat L$. Note that $|L/Q|=2^{d/2}$. 
Then $\hat Q$ is a maximal abelian subgroup of $\hat L$ which
contains $\hat D_1^{d}\cong D_1^{d}\times \<\k\>$ and which contains $K$.
Let $\psi:\hat Q\to \<\pm 1\>$ be a character of $\hat Q$ such that 
$\psi(\k K)=-1$. Then $T$ can be realized as the induced $\hat L$-module
$T=\C[\hat L]\otimes_{\C[\hat Q]}\C_{\psi}$ where $\C_{\psi}$ is one dimensional
$\hat Q$-module defined by the character $\psi$.  For $a\in\hat L$, set 
$t(a)=a\otimes 1\in T$. 
It is easy to see that we can choose 
$a_i\in \hat D_1^{d}$ such that 
$\bar a_i=\alpha_i$ for all $i$ and $\psi(a_i)=1$. 
Then for $a\in \hat L$
\begin{equation}\label{dgh2}
a_it(a)=a_ia\otimes 1=(-1)^{\<\alpha_i,\bar a\>}aa_i\otimes 1=
(-1)^{\<\alpha_i,\bar a\>}a\otimes 1.
\end{equation}
Thus,  $\C\ t(a)$ is a one-dimensional representation for $D_1^{d}$,
with character $\chi$ given by $\chi(a_i)= (-1)^{\<\alpha_i,\bar a\>}$.
In fact, $\C\ t(a)$ is isomorphic to $T_{\chi}$ as  $D_1^{d}$-modules.
Let $\b_l\in L$ for $l=1$, $\ldots$, $2^{d/2}$ 
represent the distinct cosets of $Q$ in $L$.  
Choose $b_l\in\hat L$ with $\bar b_l=\b_l$ for all $l$. 
Then $\{t(b_l)\mid l=1,\,\ldots,\,2^{d/2}\}$
forms a basis of $T$ and each $t(b_l)$ spans  a one-dimensional module
for $D_1^{d}$. Denote the character of $D_1^{d}$ on $\C\ t(b_l)$ by $\chi_l$.
Then $M(1)\otimes t(a)$ is isomorphic to $V_{D_1^d}^{T_{\chi_l}}$ as 
$\theta$-twisted $V_{D_1^{d}}$-modules and as $T_{2d}$-modules. 
Thus, $V_L^T\cong \bigoplus_{l=1}^{2^{d/2}} V_{D_1^{d}}^{T_{\chi_l}}$.

\begin{prop}\label{t4.1} 
Let $\Theta\subseteq\Delta$ be a complete coset system for
the induced $\Z_4$-subcode $Q/D_1^d$ of $\Delta$. 
We have the Virasoro decomposition
\begin{eqnarray*}
V_{L}^T  & = & \bigoplus_{\gamma\in \Theta}
                V_{D_1^{d}}^{T_{\varphi_{\gamma}}}=
\bigoplus_{ \mu\in\{+,-\}^{d}}\bigoplus_{\gamma\in \Delta}
               {\bf R}_{\mu}^{1}(\gamma),
\end{eqnarray*}
where the character $\varphi_{\gamma}$ is determined by
$\varphi_{\gamma}(a_i)=(-1)^{2\gamma_i}$ if we identify 
$D_1^*/D_1\cong \Z_4$ with $\{0,1,\pm\frac{1}{2}\}$. 
\end{prop}

\pf The first equality has been proven in the previous discussion. 
In order to see the second equality note that by Lemma~\ref{l3.3} we have 
\begin{equation}\label{dgh1}
{\bf R}^{1}_{\mu}(\gamma)=2^{-d/2}\bigotimes_{k=1}^d
(V_{D_1}^{T_{\chi_{x_k}}})^{\mu_k}
\end{equation}
where $x_k=(-1)^{2\gamma_k}$. 
Observe that
$\<\a,\b\>\in 2\Z$ for any $\a$, $\b\in Q$.  Let $a$, $b\in \hat L$ such
that $\bar a+Q=\bar b +Q$. Then from (\ref{dgh2}), $M(1)\otimes t(a)$
and $M(1)\otimes t(b)$ are isomorphic $\theta$-twisted $V_{D_1^d}$-modules
and are isomorphic $T_{2d}$-modules. Thus 
${\bf R}_{\mu}^{1}(\gamma)= {\bf R}_{\mu}^{1}(\gamma')$
if $\gamma$ and $\gamma'$ are in the same coset of $Q/D_1^d$ in $\Delta$.
Since the coset of $Q/D_1^d$ in $\Delta$ has exactly $2^{d/2}$ elements, 
we immediately see from (\ref{dgh1}) that for $\gamma\in \Delta$
$$V_{D_1^d}^{T_{\varphi_{\gamma}}}=
 \bigoplus_{ \mu\in\{+,-\}^{d}}\bigoplus_{\sigma\in Q/D_1^{d}}{\bf R}_{\mu}^{1}(\gamma+\sigma).$$
This proves the second equality. \qed
\smallskip

{\bf Proof of Theorem~\ref{d1-decomp}:} For a subset $N$ of $L$ we denote 
by $\hat N$ the inverse image of $N$ in $\hat L$ and set 
$V[N]=M(1)\otimes\C\{\hat N\}$ 
Then $V_L=\bigoplus_{\gamma\in \Delta}V[D_1^d+\gamma]$ and $V[D_1^d+\gamma]$
is isomorphic to $V_{D_1^d+\gamma}$ (defined in Section~3) 
as $V_{D_1^d}$-modules.
The decomposition for $V_L$ follows immediately from
Corollary~\ref{c3.1} and Lemma~\ref{ldecomp}.

Now we study the decomposition of $V_L^+$. If $\gamma_j=\pm \frac{1}{2}$ for 
some $j$ then 
$$(V[D_1^d+\gamma]\oplus V[D_1^d-\gamma])^+$$ 
is isomorphic to $V[D_1^d+\gamma]$ as $T_{2d}$-modules and 
has the desired decomposition.
So we can assume that all $\gamma_j=0$, $1$. In Lemma~\ref{inv} below 
we will prove that the $\theta$ defined on $V_{D_1^d+\gamma}$ in Section~3 
coincides with
the $\theta $ on $V[D_1^d+\gamma]$.
We again use Corollary~\ref{c3.1} to see that 
$V[D_1^d+\gamma]^+$ has the desired decomposition.  	

For the twisted part, we use Proposition~\ref{t4.1} and
Corollary~\ref{c3.4} to obtain
\begin{eqnarray*}
(V_{L}^T)^+=\bigoplus_{\gamma\in\Delta}
\left({\bf R}_{\mu}^{1}(\gamma)\right)^+=\bigoplus_{ \gamma\in\Delta}
\bigoplus_{ {\mu\in\{+,-\}^{d}}\atop \prod\mu_k=(-)^{d/8}} 
{\bf R}_{\mu}^{1}(\gamma).\qquad\qed
\end{eqnarray*}

\medskip

\begin{lem}\label{inv} With the same notations as in the the proof
of Theorem~\ref{d1-decomp},  
the $\theta$ defined on $V_{D_1^d+\gamma}$ in Section 3 coincides with
the $\theta$ on $V[D_1^d+\gamma]$ if all $\gamma_j=0$, $1$.
\end{lem}

\pf Let $X$ be a sublattice of $L$
containing $D_1^d$ such that $\<x,y\>\in 2\Z$ for $x,y\in X$. Then 
the inverse image $\hat X$ of $X$ in $\hat L$ is an abelian group.
We can choose a section $e: X\to \hat X$ such that 
$e_{x}e_y=\kappa^{\<x,y\>/2}e_{x+y}$. Then 
$e_x^{-1}=\kappa^{\<x,x\>/2}e_{-x}$ for $x\in X$. Thus $\theta\iota(e_x)
=\iota(e_{-x})$. 

Take $X$ to be the sublattice generated by $D_1^d$ and $\gamma$.  
Then $V[D_1^d+\gamma]$ is generated by
$\iota(e_{\gamma})$ as $V_{D_1^d}$-module and $V[D_1^d+\gamma]^+$
is generated by $\iota(e_{\gamma})+\iota(e_{-\gamma})$ as
$V_{D_1^d}^+$-module. In fact  $V[D_1^d+\gamma]^+$ is an irreducible
$V_{D_1^d}^+$-module. Let $\mu: V_{D_1^d+\gamma}\mapsto V[D_1^d+\gamma]$
be an $V_{D_1^d}$-module isomorphism such that $\mu e^{\gamma}=\iota(e_{\gamma})$. 
We must prove that $\mu$ maps $V_{D_1^d+\gamma}^+$ to 
$(V[D_1^d+\gamma])^+$. As both $V_{D_1^d+\gamma}^+$ and 
$(V[D_1^d+\gamma])^+$ are irreducible $V_{D_1^d}^+$-modules, it is enough
to show that $\mu(e^{\gamma}+e^{-\gamma})=
\iota(e_{\gamma})+\iota(e_{-\gamma})$ or 
equivalently $\mu e^{-\gamma}=\iota(e_{-\gamma})$.

Let $J$ be the subset of $\{1,\ldots,d\}$ 
consisting of $j$ such that $\gamma_j=1$.  Note that 
$e^{-\gamma}$ is the coefficient of $z^{-2|J|}$ in 
$Y(e^{-2\gamma},z)e^{\gamma}$ and $\iota(e_{-\gamma})$  is the coefficient of 
$z^{-2|J|}$ in $Y(\iota(e_{-2\gamma}),z)\iota(e_{\gamma})$ 
as $(-1)^{\<2\gamma,\gamma\>/2}$ is even. Also note that $2\gamma\in D_1^d$. 
Thus $\mu e^{-\gamma} =\iota(e_{-\gamma})$. 
\qed

\medskip

Now we return our lattices $\lc$ and $\tlc$ associated to the
code $C$. We assume that $C$ is a {\em self-annihilating\/} 
(i.e., $C=C^\perp$) doubly-even binary code.
Then the $\Z_4$-codes $\Gamma$ and $\widetilde{\Gamma}$
are self-annihilating and even, or equivalently as mentioned above,
the lattices $\lc$ and $\tlc$ are self-dual and even. 

Combining Lemma~\ref{ldecomp} and Theorem~\ref{d1-decomp} we will see how
to read off the Virasoro decomposition directly from the marked code $C$.
For $a$,~$b\in \{0,1\}$, $\alpha$,~$\beta\in\{+,-\}$ and 
$x$,~$y \in\{0,1\}$, define the formal linear combinations of 
$T_4$-modules $N_{\alpha\beta}^{ab}((x,y))$ by
$$N_{\alpha\beta}^{ab}((x,y))=
\bigoplus_{{\alpha',\alpha''\in\{+,-\}}\atop{\alpha'\alpha''=\alpha}}
{\bf R}^a_{(\alpha',\alpha'')}(\gamma^b_{\beta}((x,y)))$$
where ${\bf R}^a_{(\alpha',\alpha'')}$ was defined in~(\ref{dgh0}). 
Explicitly,  we get the $64$ formal linear combinations as shown
in the following table:
{\scriptsize $$\begin{array}{c||cccc}
  & (0,0) & (1,1) & (0,1), (1,0)  \\ \hline\hline
N_{++}^{00} & M(0,0,0,0)\oplus M(\ha,\ha,\ha,\ha) & M(\ha,0,0,0) \oplus M(0,\ha,\ha,\ha) &
        \ha   M(\se,\se,\se,\se) \\
N_{-+}^{00} & M(0,0,\ha,\ha)\oplus M(\ha,\ha,0,0) & M(\ha,0,\ha,\ha) \oplus M(0,\ha,0,0) &
        \ha    M(\se,\se,\se,\se)  \\
N_{+-}^{00} & M(\ha,0,\ha,0)\oplus M(0,\ha,0,\ha) & M(0,0,\ha,0)\oplus M(\ha,\ha,0,\ha) &
        \ha M(\se,\se,\se,\se)  \\
N_{--}^{00} & M(\ha,0,0,\ha)\oplus M(0,\ha,\ha,0) & M(0,0,0,\ha) \oplus M(\ha,\ha,\ha,0) &
  \ha M(\se,\se,\se,\se)  \\ \hline
N_{++}^{01},N_{-+}^{01} & \multicolumn{2}{c}{
\ha M(\se,\se,0,0)\oplus \ha M(\se,\se,\ha,\ha)}  &  
\ha M(0,0,\se,\se) \oplus \ha M(\ha,\ha,\se,\se) \\
N_{+-}^{01},N_{--}^{01} &  \multicolumn{2}{c}{
\ha M(\se,\se,\ha,0)\oplus\ha M(\se,\se,0,\ha) } &  
\ha M(\ha,0,\se,\se) \oplus \ha M(0,\ha,\se,\se) \\ \hline
N_{++}^{10},N_{+-}^{10} & \multicolumn{2}{c}{
\ha M(\se,0,\se,0) \oplus \ha M(\se,\ha,\se,\ha)} & 
\ha M(0,\se,0,\se)\oplus \ha M(\ha,\se,\ha,\se) \\
N_{-+}^{10},N_{--}^{10} &  \multicolumn{2}{c}{ 
\ha M(\se,\ha,\se,0) \oplus \ha M(\se,0,\se,\ha)} &
\ha M(\ha,\se,0,\se)\oplus \ha  M(0,\se,\ha,\se)
\\ \hline
N_{++}^{11},N_{+-}^{11}&  \multicolumn{2}{c}{
\ha M(0,\se,\se,0)\oplus\ha  M(\ha,\se,\se,\ha) } & 
\ha M(\se,0,0,\se)\oplus \ha M(\se,\ha,\ha,\se) \\
N_{-+}^{11},N_{--}^{11}&  \multicolumn{2}{c}{
\ha  M(\ha,\se,\se,0)\oplus \ha M(0,\se,\se,\ha) } & 
\ha M(\se,\ha,0,\se)\oplus\ha  M(\se,0,\ha,\se) 
\end{array}$$ }
For $\mu$,~$\epsilon\in\{+,-\}^{n/2}$ and an element $c\in\F_2^n$ we write 
$${\bf N}^{ab}_{\mu,\epsilon}(c)=\bigotimes_{k=1}^{n/2}
N_{\mu_k\epsilon_k}^{ab}(c(k)),$$
where $c(k)=(c_{2k-1},c_{2k})$ as before. Let $\delta(c)$ be the
number of $k$ with $c(k)\in\{(0,1),(1,0)\}$. Recall that the lattice $D_1^d$ 
determines $2d$ commuting Virasoro algebras inside the four vertex operator 
algebras $\vlc$, $\vtlc$, $\tvlc$ and $\tvtlc$.
The following is the main theorem of this paper.

\begin{thm}\label{decomp}
For the Virasoro frame coming from the marked code $C$,  
we have the following decompositions 
\begin{eqnarray*}
& &\vlc = \bigoplus_{c\in C}
\bigoplus_{ {\mu,\epsilon\in\{+,-\}^{d/2}} }
{\bf N}_{\mu,\epsilon}^{00}(c),
\\ 
& &\vtlc =  \bigoplus_{c\in C\atop \delta(c)=0}
\bigoplus_{ {\mu,\epsilon\in\{+,-\}^{d/2}} \atop \prod\epsilon_k=+}
{\bf N}_{\mu,\epsilon}^{00}(c)
\oplus \ha\cdot
\bigoplus_{c\in C,\atop \delta(c)>0}
\bigoplus_{ {\mu,\epsilon\in\{+,-\}^{d/2}} } 
{\bf N}_{\mu,\epsilon}^{00}(c)
\\ & &\ \ \ \ \ \ \ \ \ \ \oplus
 \bigoplus_{ c\in C}
\bigoplus_{ {\mu,\epsilon\in\{+,-\}^{d/2}}
\atop \prod\epsilon_k=(-)^{d/8}} 
{\bf N}_{\mu,\epsilon}^{01}(c),
\\ 
& &\tvlc  =  \bigoplus_{c\in C\atop \delta(c)=0}
\bigoplus_{ {\mu,\epsilon\in\{+,-\}^{d/2}} \atop \prod\mu_k=+ } 
{\bf N}_{\mu,\epsilon}^{00}(c)
\oplus \ha\cdot
\bigoplus_{c\in C,\atop \delta(c)>0}
\bigoplus_{ {\mu,\epsilon\in\{+,-\}^{d/2}} }
{\bf N}_{\mu,\epsilon}^{00}(c)  
\\& &\ \ \ \ \ \ \ \ \ \ 
\oplus\bigoplus_{c\in C}
\bigoplus_{ {\mu,\epsilon\in\{+,-\}^{d/2}}\atop 
\prod\mu_k=(-)^{d/8} } 
{\bf N}_{\mu,\epsilon}^{10}(c),
\\ 
& &\tvtlc  = \bigoplus_{c\in C,\atop \delta(c)=0}
\bigoplus_{ {\mu,\epsilon\in\{+,-\}^{d/2}} \atop \prod\mu_k=\prod\epsilon_k=+ } 
{\bf N}_{\mu,\epsilon}^{00}(c)\oplus\
 \frac{1}{4}\cdot
\bigoplus_{c\in C,\atop \delta(c)>0}
\bigoplus_{ {\mu,\epsilon\in\{+,-\}^{d/2}} }
{\bf N}_{\mu,\epsilon}^{00}(c)
\\ & & \ \ \oplus \ha \cdot \bigoplus_{c\in C}\big[
\bigoplus_{ {\mu,\epsilon\in\{+,-\}^{d/2}}\atop 
\prod\epsilon_k=(-)^{d/8} } 
{\bf N}_{\mu,\epsilon}^{01}(c) \oplus 
\bigoplus_{ {\mu,\epsilon\in\{+,-\}^{d/2}}\atop 
\prod\mu_k=(-)^{d/8} } 
{\bf N}_{\mu,\epsilon}^{10}(c) \oplus 
\bigoplus_{ {\mu,\epsilon\in\{+,-\}^{d/2}}\atop 
\prod\epsilon_k=\prod\mu_k=(-)^{d/8} }
{\bf N}_{\mu,\epsilon}^{11}(c)\big].
\end{eqnarray*}
\end{thm}

\pf Recall~(\ref{dgh-1}). 
For a codeword $c\in C$ and $\epsilon\in \{+,-\}^{d/2}$ 
let $\gamma=\Gamma_{\epsilon}^b(c)$ and fix
$s\in\{+,-\}$. 
Using the definition of ${\bf N}^{ab}_{\mu,\epsilon}(c)$ we get
\begin{eqnarray*}
\bigoplus_{ {{\mu\in\{+,-\}^{d}}\atop \prod\mu_k=s} } 
      {\bf R}_{\mu}^a(\gamma) 
&= & \bigoplus_{ {{\mu',\mu''\in\{+,-\}^{d/2}}\atop \prod\mu'_k\mu''_k=s} }
\bigotimes_{i=1}^{d/2}
{\bf R}^a_{(\mu_i',\mu_i'')}((\gamma_{2i-1},\gamma_{2i}))
\\
& = & \bigoplus_{ {{\mu\in\{+,-\}^{d/2}}\atop \prod\mu_k=s} }
\bigotimes_{i=1}^{d/2}
\bigoplus_{ {\mu_i'\,\mu_i''\in\{+,-\}}\atop {\mu_i'\mu_i''=\mu_i}}
      {\bf R}^a_{(\mu_i',\mu_i'')}((\gamma_{2i-1},\gamma_{2i})) \\
& = & \bigoplus_{{ {\mu\in\{+,-\}^{d/2}}\atop \prod\mu_k=s }}
\bigotimes_{i=1}^{d/2}
 {N}^{ab}_{\mu_i\epsilon_i}(c(i))                   \\
&= & \bigoplus_{{ {\mu\in\{+,-\}^{d/2}}\atop \prod\mu_k=s }}
      {\bf N}^{ab}_{\mu,\epsilon}(c).
\end{eqnarray*}
By using the identification of $\Gamma_{\epsilon}^a(c)$ with codewords in
$\Gamma$ and $\tilde\Gamma$ the decomposition follows from
Lemma~\ref{ldecomp} and Theorem~\ref{d1-decomp} if we use the following two
observations:

First note that ${N}_{\mu_k\epsilon_k}^{00}(c(k))=
{N}_{\pm\mu_k\,\pm\epsilon_k}^{00}(c(k))$ for $c(k)\in\{(0,1),(1,0)\}$.
So for $\delta(c)>0$ we can for 
suppress the distinction between $\pm\epsilon_k$ 
(resp.~$\pm\mu_k$) in the decomposition and compensate it with 
one factor $\frac{1}{2}$. 

Second, the value of ${\bf N}^{01}_{\mu,\epsilon}(c)$ 
(resp.~${\bf N}^{10}_{\mu,\epsilon}(c)$ and ${\bf N}^{11}_{\mu,\epsilon}(c)$) 
depends for fixed $c$ only on $\epsilon$ (resp.~$\mu$). \qed

\smallskip

Now we discuss an action of $Sym_3$ (the permutation group
on three letters) defined in [FLM] and [DGM2] on $\vlc$ and $\tvtlc$
in terms of our decompositions. The resulting group of 
automorphisms is sometimes called the {\em triality group\/}.   

Recall from Chapters~11 and~12 of~[FLM] and Sections~7 and~8 
of~[DGM2] that the triality group is generated by distinct
involutions $\sigma$ and $\tau$. 
Also recall from Section~4 that $\a_1$, $\ldots$ , 
$\a_{d}$ form a $D_1$-frame in $L$. A straightforward computation
shows that $\sigma \omega_{4i-3}=\omega_{4i-3}$, 
$\sigma \omega_{4i}=\omega_{4i}$ and $\sigma$ interchanges 
$\omega_{4i-2}=\omega_{4i-1}$ for all $i=1$, $\ldots$, ${d\over 2}$.
Similarly, $\tau$ interchanges $\omega_{4i-3}=\omega_{4i-2}$ and
fixes $\omega_{4i-1}$ and $\omega_{4i}$. Thus the triality group 
is a subgroup of $G$ defined in (\ref{g1}) for both $\vlc$ and $\tvtlc$. 
Its image in $G/G_{\cal C}\leq Sym_{2d}$ is the above described
permutation of the elements of the VF $\{\omega_{\nu}\}$.

Additionally, the involution $\sigma$ defines an isomorphism
between $\vtlc$ and $\tvlc$. 

\smallskip

\begin{de}\label{dmz} \rm
Following~[DMZ], the {\em decomposition polynomial\/} of a FVOA
$V=\bigoplus m_{h_1,\ldots,h_r}\,M(h_1,\ldots,h_r)$
is defined as 
$$P_V(a,b,c)  =\sum_{i,j,k}A_{i,j,k}\,a^ib^jc^k$$
where $A_{i,j,k}$ is the number of $T_r$-modules 
$M(h_1,\ldots,h_r)$ in a $T_r$ composition series of $V$ for which the number of
coordinates in $(h_1,\ldots,h_r)$ equal to $0$, $\ha$, $\se$ is $i$, $j$, $k$,
respectively.
\end{de}
The polynomial is homogeneous of degree $r$ and, in general, depends 
on the chosen Virasoro frame $\{\omega_1,\ldots,\omega_r\}$ inside of $V$.

The following corollary is an immediate consequence of Theorem~\ref{decomp}.
\begin{cor}\label{decomp-pol} 
Using the symmetrized marked weight enumerator ${\rm smwe}_C(x,y,z)$ one has
\begin{eqnarray*}
P_{\vlc}(a,b,c)&=&
{\rm smwe}_C(a^4+6a^2b^2+b^4,2c^4,4a^3b+4ab^3), \\   
P_{\vtlc}(a,b,c)&=&
\frac{1}{2}\left(a^4-2a^2b^2+b^4\right)^{\frac{d}{2}} 
+\frac{1}{2}\,{\rm smwe}_C(a^4+6a^2b^2+b^4,2c^4,4a^3b+4ab^3) \\
& &\qquad +\frac{1}{2}\cdot 2^{d/2}\left((a+b)^{d}+
                    (-1)^{d/8}(a-b)^{d}\right)c^{d}, \\ 
P_{\tvlc}(a,b,c)&= &P_{\vtlc}(a,b,c), \\ 
P_{\tvtlc}(a,b,c)&=&
\frac{1}{4}\cdot 3\left(a^4-2a^2b^2+b^4\right)^{\frac{d}{2}} \\
& &\qquad+\frac{1}{4}\,{\rm smwe}_C(a^4+6a^2b^2+b^4,2c^4,4a^3b+4ab^3) \\
& &\qquad + 3\cdot\frac{1}{4}\cdot 2^{d/2}\left((a+b)^{d}+
(-1)^{d/8}(a-b)^{d}\right)c^{d}.
\qquad\qed
\end{eqnarray*}
\end{cor}
 
\begin{rem}{\rm 
{}From Theorem~\ref{decomp} we can deduce that 
$\tvtlc$ is a self-dual rational vertex operator algebra. The
proof for the special case of $V^{\natural}$ given in~[D3] works
in general since the Virasoro decompositions were the only
information needed. 
}\end{rem}
\section{Applications}

In this section we discuss some important applications for Theorem~\ref{decomp}.
The simplest example is for the Hamming code  $H_8$ of length $8$.
When $C$ is the Golay code ${\cal G}_{24}$ of length $24$ there 
is a special marking and we obtain a particular interesting decomposition of 
the moonshine  module 
$V^{\natural}=\widetilde{V}_{ \widetilde{L}_{{\cal G}_{24}}}$ 
under $48$ Virasoro algebras.

\medskip

\noindent{\em Example I: The Hamming code $H_8$, the root lattice $E_8$, and 
the lattice vertex operator algebra $V_{E_8}$}

\smallskip

The Hamming code $H_8$ is the unique self-annihilating doubly-even binary code of 
length $8$.  Its automorphism group is isomorphic to $AGL(\F_2^3)$.
The root lattice $E_8$ of the exceptional Lie group $E_8(\C)$ is the unique
even unimodular lattice of rank $8$. It has the Weyl group $W(E_8)$ as its
automorphism group. 
The lattice vertex operator algebra $V_{E_8}$, whose underlying
vector space is the irreducible level~$1$ highest weight representation of the 
affine Kac-Moody algebra $E_8^{(1)}$, is a
self-dual vertex operator algebra of rank $8$ whose automorphism group is 
the Lie group $E_8(\C)$.
One can show, under some additional conditions on the \voa, that
$V_{E_8}$ is the unique self-dual VOA of rank $8$ (cf.~[H1], Ch.~2).

The uniqueness of the lattice $E_8$ implies 
$E_8\cong L_{H_8}\cong \widetilde{L}_{H_8}$
and $V_{E_8}\cong V_{L_{H_8}} \cong V_{\widetilde{L}_{H_8}} 
\cong \widetilde{V}_{L_{H_8}} \cong \widetilde{V}_{\widetilde{L}_{H_8}}$
for the vertex operator algebras, since one has $\vtlc\cong\tvlc$ in general
(see~[DGM1], [DGM2] and the remark after Theorem~\ref{decomp}). 

\smallskip
We will determine up to automorphism all markings for the Hamming code, all
$D_1$-frames of the $E_8$ lattice, and five Virasoro frames inside $V_{E_8}$ 
and describe the corresponding decompositions. 
They are all coming from markings of the Hamming code.

\medskip
To fix notation we choose $\{(00001111),(00110011),(11000011),(01010101)\}$ as 
a set of base vectors for the Hamming code.

\begin{thm}\label{ham}
There are $3$ orbits of markings for the Hamming code $H_8$ under 
${ Aut}(H_8)$. Their main properties can be found in the next table.
The last column shows the symmetrized marked weight enumerator.

\smallskip
\noindent{\rm\footnotesize
\begin{tabular}{llccc}
orbit & orbit representatives & stabilizer & orbit size & ${\rm smwe}_{H_8}(x,y,z)$ 
                                          \\ \hline
$\alpha$ & $\left\{(1,2),(3,4),(5,6),(7,8)\right\}$ & $2^3\colon Sym_4$ & 
$7$ & $x^4+6\,x^2z^2+z^4+8\,y^4$   \\
$\beta$ & $\left\{(1,2),(3,4),(5,7),(6,8)\right\}$ & $2^2.Dih_8$ & 
$ 42$ & $x^4 +2\,x^2z^2+z^4 +8\,xzy^2 + 4\,y^4 $    \\
$\gamma$ & $\left\{(1,2),(3,5),(4,7),(6,8)\right\}$ & $Sym_4$ & 
$56$ & $x^4+z^4+12\,xzy^2 +2\,y^4$ 
\end{tabular}}
\end{thm}

The proof is an easy counting exercise (see Appendix~\ref{hamming-code}).

\smallskip 

We remark that every pair $(i,j)$ of the eight positions is the component of
exactly one of the seven markings of type $\alpha$:  
Every marking contains $4$ pairs, so we cover 
$7\cdot 4=28$ pairs. There are ${8\choose 2}=28$ different such pairs on 
which ${ Aut}(H_8)$ transitively acts.

\medskip
As explained in the last section before Lemma~\ref{ldecomp} in general, every 
marking of $H_8$ determines a $D_1^8$ sublattice inside 
$L_{H_8}\cong E_8$ and $\widetilde{L}_{H_8}\cong E_8$.

The following theorem shows that all possible $D_1$-frames in
$E_8$ are obtained in this way.
\begin{thm}[Conway-Sloane \hbox{[CS2]}]\label{e8-lattice}
There are $4$ orbits of $D_1^8$ sublattices inside $E_8$ under the action of
$W(E_8)$. Their main properties are shown in the next table. The column
``origin'' lists the corresponding (untwisted, resp.~twisted lattice) 
Hamming code marking and ${\rm swe}_{\Delta}(A,B,C)$ is the 
symmetrized weight enumerator of the decomposition code 
$\Delta=E_8/D_1^8\leq (D_1^*/D_1)^8\cong\Z_4^8$.

\smallskip
\noindent{\rm\footnotesize
\begin{tabular}{llccl}
orbit & origin &  $stabilizer$ & orbit size &
                          ${\rm swe}_{\Delta}(A,B,C)$ \\ \hline
${\cal K}_8$ &$\alpha$  & $\Z_2^7.\cdot Sym_8$ &  $135$ &
  $A^8+28A^2C^6+70A^4C^4+28A^6C^2+C^8+128B^8$      \\
${\cal K}'_8$ &$\beta $, $\widetilde{\alpha}$ & $2^7(4!)^2$ & $9450 $ &
$A^8+ C^8 + 12A^2C^2(A^4+C^4)+38 A^4C^4+$ \\
&&&& \qquad\qquad\qquad $ 64AC(A^2+C^2)B^4 + 64B^8$    \\
${\cal L}_8$ &$\gamma$, $\widetilde{\beta}$ & $2^8\cdot 4!$ &  $113400$ &
 $A^8+C^8+4A^2C^2(A^4+C^4)+22A^4C^4+$ \\
&&&& \qquad\qquad\qquad $96AC(A^2+C^2)B^4 + 32B^8$       \\
${\cal O}_8$ & \phantom{$\gamma$, } $\widetilde{\gamma}$ & $2.AGL(3,2)$ &$259200$ &
$ A^8+C^8+14A^4C^4+112AC(A^2+C^2)B^4+16B^8$       \\
\end{tabular}}
\smallskip
\end{thm}

\pf It was also explained in the last section that every $D_1^8$ sublattice
inside $E_8$ defines a $\Z_4$-code $\Delta\leq (D_1^{*}/D_1)^8\cong \Z_4^8$. 
Since $E_8$ is self-dual and even, $\Delta$ is self-annihilating
and even as a code over $\Z_4$. 
All self-annihilating $\Z_4$-codes of length $8$ are classified 
in [CS2], Theorem~2. Only ${\cal K}_8$, ${\cal K}'_8$, ${\cal L}_8$ and 
${\cal O}_8$ are even (see also [BS]). 
The order of ${ Aut}(\Delta)$ and ${\rm swe}_{\Delta}(A,B,C)$
are also described in~[CS2].
To show that these codes arise from the markings of the Hamming code as 
in the table we apply Corollary~\ref{ldecomp-cor}.
\qed

The remark after Theorem~\ref{ham}  about 
the Hamming code  
has an analogue here: every vector of squared 
length~$4$ inside $E_8$ is contained in exactly one $D_1$-frame belonging
to the orbit of type ${\cal K}_8$ since $135\cdot 16=2160$, the number
of vectors of squared length~$4$, and $W(E_8)$ acts transitively on such vectors.
These $135$ $D_1^8$-sublattices are in bijection with cosets of $2E_8$
in $E_8$ which have coset representatives of norm $4$.

\medskip
Every $D_1$-frame inside $E_8$ determines $16$ commuting Virasoro 
vertex operator algebras of rank~$\ha$ inside $V_{E_8}$ and  
$\widetilde{V}_{E_8}\cong V_{E_8}$. Altogether, one gets at least 
five different systems of commuting Virasoro subVOAs:

\begin{thm}
Let $\{\omega_1,\ldots,\omega_{16}\}$ be a Virasoro frame inside 
$V_{E_8}$. The possible decomposition polynomials are displayed
in the next table. They correspond by the untwisted or twisted lattice
construction to the $D_1$-frames inside $E_8$ as indicated in
the column origin. Furthermore, the first four cases belong to four
distinct orbits of Virasoro frames under the action of the Lie group $E_8(\C)$.
In the fifth case,  $\Omega$,   at least the $T_{16}$-module structure
is unique.

\smallskip{\small
\noindent\centerline{\begin{tabular}{lll}
case & origin & $P_{V_{E_8}}(a,b,c)$  \\ \hline
$\Gamma$ & ${\cal K}_8$ & $\frac{1}{2}\left[ (a+b)^{16}+(a-b)^{16} \right] +
  128\,c^{16}$     \\
$\Sigma$ & ${\cal K}'_8$, $\widetilde{{\cal K}}_8$ &  
 ${a^{16}} + {b^{16}}+
  56\,({a^{14}}\,{b^2} + \,{a^2}\,{b^{14}}) +
  924\,({a^{12}}\,{b^4} + \,{a^4}\,{b^{12}}) + $\\&&$\quad
 3976\,({a^{10}}\,{b^6} +\,{a^6}\,{b^{10}}) +
  6470\,{a^8}\,{b^8} + $\\&&$\quad
   \big(128\,({a^7}\,b+\,a\,{b^7})+ 896\,({a^5}\,{b^3}\,+{a^3}\,{b^b})\big){c^8} + 
   64\,{c^{16}}$ \\
$\Psi$ & ${\cal L}_8$, $\widetilde{{\cal K}}'_8$  & 
${a^{16}}+ {b^{16}} + 24\,({a^{14}}\,{b^2}+{a^2}\,{b^{14}} ) + 
   476\,({a^{12}}\,{b^4}+{a^4}\,{b^{12}} ) +  $\\&&$\quad
   1960\,({a^{10}}\,{b^6} + {a^6}\,{b^{10}})
  +3270\,{a^8}\,{b^8}  
     + $\\&&$\quad
  \big( 192\,({a^7}\,b + a\,{b^7})+
    1344\,({a^5}\,{b^3} + {a^3}\,{b^5}) \big) {c^8}
   + 32\,{c^{16}}$\\
$\Theta$ & ${\cal O}_8$, $\widetilde{{\cal L}}_8$ &
${a^{16}} +  {b^{16}}+ 
   8\,({a^{14}}\,{b^2} + {a^2}\,{b^{14}} )+
   252\,({a^{12}}\,{b^4} +{a^4}\,{b^{12}} )+ $\\
&&$\quad
   952\,({a^{10}}\,{b^6} + {a^6}\,{b^{10}} )
   1670\,{a^8}\,{b^8}  + $\\
&&$\quad
   \big( 224\,({a^7}\,{b} + {a}\,{b^7}) +
  1568\,({a^5}\,{b^3} + {a^3}\,{b^5}) \big) {c^8}
    + 16\,{c^{16}} $ \\
$\Omega$ & \phantom{${\cal O}_8$, } $\widetilde{{\cal O}}_8$ &  
${a^{16}} + {b^{16}}+ 140\,({a^{12}}\,{b^4} +{a^4}\,{b^{12}})+
448\,({a^{10}}\,{b^6} + {a^6}\,{b^{10}})+ 870\,{a^8}\,{b^8} + $\\&&$\quad
  \big( 240\,({a^7}\,b\,+a\,{b^7}) + 
   1680\,({a^5}\,{b^3}+ {a^3}\,{b^5})\big){c^8} + 8\,{c^{16}}$    \\
\end{tabular}}}
\smallskip
\end{thm}

\pf 
Using Corollary~\ref{decomp-pol} to compute $P_{V_{E_8}}(a,b,c)$ for 
the different Virasoro subVOAs $T_{16}$ coming from 
$V_{E_8}$, and $\widetilde{V}_{E_8}$ and a given $D_1^8$ sublattice in 
$E_8$ one checks that the polynomials
for $\Gamma$, $\Sigma$, $\Psi$, $\Theta$ and $\Omega$ correspond to the 
$D_1$-frames of $E_8$ as indicated. 

We show that there are no other possibilities for the decomposition polynomial
$P_{V_{E_8}}(a,b,c)$ and we will see directly that there 
is a unique ${Aut}(V_{E_8})$-orbit of $T_{16}$ subVOAs corresponding
to each of the cases $\Gamma$, $\Sigma$, $\Psi$ and $\Theta$. 

\smallskip
Assume a vertex operator subalgebra $T_{16}$ in $V_{E_8}$ is
given. First we determine the possible decomposition polynomials.

As described in the proof of Theorem~4.1.5 in~[H1], ${SL}_2(\Z)=
\langle S,T\rangle$ with $S=\left({\ 0\ 1}\atop {-1\ 0}
\right)$ and $T=\left({1\ 1}\atop {0\ 1}\right)$ acts on $\C[a,b,c]$ by 
$$\rho(S)=\left(\begin{array}{ccc} 
1/2 &1/2 & 1/\sqrt{2}\\1/2  &1/2 &-1/\sqrt{2} \\1/\sqrt{2} & -1/\sqrt{2}& 0
\end{array}\right), \qquad 
\rho(T)=e^{-2\pi i/48} \left(\begin{array}{ccc} 1& 0&0 \\ 0 & -1 & 0 \\
0 &0 & e^{2\pi i/16}\end{array}\right)$$
via substitution. Since $V_{E_8}$ is a self-dual VOA of rank $8$, the 
decomposition polynomial must be invariant under the action of $\rho(S)$ and
$\rho(T^3)$ (cf.~proof of Theorem~\ref{codeholo} or~Th.~2.1.2 and Th.~4.1.5 in~[H1]). 
They generate a matrix group $G=\langle \rho(S),\rho(T)^3\rangle$
of order $384$ as can easily be seen with the help of the program Gap~[Sgap].
The dimension of the space of invariant 
polynomials of degree $n$ is the multiplicity of the trivial representation
in the $n^{\rm th}$ symmetric power of $\rho$. This multiplicity is given by the
coefficient of $t^n$ in the expression 
$$\rho_G(t)=\frac{1}{|G|}\sum_{g\in G} \frac{1}{{\rm det}({\rm id} -g t)}\ .$$
For degree $16$ we obtain that the space of invariant
polynomials is two dimensional;
a possible base is given by $P_{V_{E_8}}^{\Gamma}(a,b,c)$ and 
$P_{V_{E_8}}^{\Omega}(a,b,c)$. The only polynomials $P(a,b,c)$ inside 
this space having positive coefficients and satisfying the 
necessary conditions $P(1,0,0)=1$ and $P(1,1,0)=|{\cal C}|=2^l$, 
$0\leq l\leq 15$ with integral $l$, are the five polynomials given 
in the theorem. 

Next we claim that the code ${\cal C}$ is uniquely determined from its weight 
enumerator $P(a,b,0)$: The weight enumerator of its annihilator code ${\cal C}^{\perp}$
is $a^{16}+(2^k-2)a^8b^8+b^{16}$, with $k=16-l=1$,~$\ldots$,~$5$. For $k=5$ the 
uniqueness of ${\cal C}^{\perp}$ and so of ${\cal C}$ is the uniqueness of 
the simplex code (see Appendix~\ref{moon-code}, Prop.~\ref{A.3.6}). 
For smaller $k$, it follows also from the proof of Prop.~\ref{A.3.6}. 

\smallskip 

The code ${\cal C}$ contains for $k=1$, $2$, $3$ and $4$ the subcode 
${\cal C}_0=\{(00),(11)\}^8$. By Corollary~\ref{c3.1} and the uniqueness 
statement of Proposition~\ref{C-voas} the corresponding subVOA must be 
$V_{D_1^8}$. Recall that $V_{E_8}=M(1)\otimes \C\{E_8\}$. 
The weight one subspace $(V_{E_8})_1$ is a Lie algebra under
$[u,v]=u_0v$ which is isomorphic to the Lie algebra of type $E_8$ 
and $(V_{D_1^8})_1$ is a Cartan subalgebra of  $(V_{E_8})_1$. 
{}From the construction of $V_{E_8}=M(1)\otimes \C\{E_8\}$ we have a 
canonical Cartan subalgebra $M(1)_1$ of  $(V_{E_8})_1$ which is identified with 
${\frak h}=\C\otimes_{\Z} E_8$. 
Since all Cartan subalgebras are conjugate under the adjoint action of the 
Lie group $E_8(\C)$ we can assume that $(V_{D_1^8})_1=M(1)_1\leq  (V_{E_8})_1$.

It is well-known that  $\C\{E_8\}=1 \otimes \C\{E_8\} = 
\{u\in V_{E_8}\mid h_nu=0,\, h\in \frak h,\, n>0\}$ 
which is the vacuum space for the Heisenberg algebra $\hat \frak h_{\Z}$ 
(see Section~3). Similarly,   
$\C\{D_1^8\}=\{u\in V_{D_1^8} \mid h_nu=0,\, h\in \frak h,\, n>0\}$. 
Thus,  $\C\{D_1^8\}$ is a subspace of $\C\{E_8\}$. 
Note that $\C\{E_8\}$ and $\C\{D_1^8\}$ are direct sums of weight spaces for the Cartan algebra $\frak h$
and the corresponding weight lattices are exactly $E_8$ and $D_1^8$. 
This determines a $D_1^8$ sublattice of $E_8$, unique up to the action of
$W(E_8)$.
That is, for $k=1$, $2$, $3$, $4$ the Virasoro frames 
$\{\omega_1,\ldots,\omega_{16}\}$ come from one of the four $D_1$-frames 
inside $E_8$ by the untwisted lattice-VOA construction.

\smallskip
It remains to show that for $k=5$, i.e.~in the case $\Omega$,
the obtained Virasoro decomposition is unique. As stated before the
code ${\cal C}^{\perp}$, which is equal to ${\cal D}$ 
by Theorem~\ref{codeholo}, is the simplex code and so ${\cal C}$ 
is the extended Hamming code of length $16$. The uniqueness of the 
code ${\cal C}$ implies by Theorem~\ref{t2.1}~(4) that $V_{E_8}^{0}$ is unique 
as a $T_{16}$-module. Let $I\in {\cal D}$ such that $|I|=8$.
Take an irreducible $T_{16}$-module $W$ in $V^I$. Using  the action 
of $\Ve$ on $W$ we see that all such $M(h_1,\ldots,h_{16})$ occur in $V^I$
where $h_i=\frac{1}{16}$ if $i\in I$ and 
$h_i\in \{0,\frac{1}{2}\}$ if $i\not\in I$
and the number of $i$ with $h_i=\frac{1}{2}$ is odd. 
So,  there are $2^7$ nonisomorphic
$T_{16}$-modules inside $V^I$. Since ${\cal D}$ has $30$ codewords of
weight $8$, we get at least $30\cdot 2^7$ such nonisomorphic $T_{16}$-modules. 
But $30\cdot 2^7$ is exactly the coefficient of $c^8$ in $P(1,1,c)$.
This shows that all these modules have 
multiplicity one. 
Finally the multiplicity of $M(\frac{1}{16},\ldots,\frac{1}{16})$ is $8$. 
Therefore the decomposition in the last case is unique.
\qed

\begin{rem}{\rm 
(1) We expect that also in the fifth case the \voa structure is unique,
i.e.~$\Omega$ corresponds to a unique ${E_8}(\C)$-orbit of Virasoro frames.

(2) A different proof would follow if the list of the $71$ unitary 
self-dual VOA {\em candidates\/} of rank $24$ given by Schellekens~[Sch] is 
complete:

\noindent{}The fusion algebras for $M(0)$ and the Kac-Moody VOA 
$V_{B_{1,1}}$ are isomorphic and one 
can identify the corresponding intertwiner spaces (cf.~[MoS], Appendix~D). 
{}From this, one can define for every VOA $V$ of rank $c$ 
containing $M(0)^{\otimes 2c}$
a VOA $W$ of rank $3c$ containing $V_{B_{1,1}}^{\otimes 2c}$. 
There are five candidates $W$ on Schellekens list containing 
$V_{{B_{1,1}}}^{\otimes 16}$, namely
$W_{D_{{24},1}}$, $W_{D_{{12},1}^2}$, 
$W_{D_{{6},1}^4}$, $W_{A_{{3},1}^8}$ and
$W_{A_{{1},2}^{16}}$. They correspond to 
$\Gamma$, $\Sigma$, $\Psi$, $\Theta$ and $\Omega$ 
in this order. Again the uniqueness of $W_{A_{{1},2}^{16}}$
is unknown. The decomposition of $W_{A_{{1},2}^{16}}$ as a 
$V_{{B_{1,1}}}^{\otimes 16}$-module obtained
in~[Sch] by a computer calculation follows from our 
analysis of the case $\Omega$.}
\end{rem}
\def\dd{$\swarrow\ \searrow$}

\medskip
The next table summarizes the relation between the markings for $H_8$, the
$D_1$-frames inside $E_8$ and the Virasoro frames 
$\{\omega_1,\ldots,\omega_{16}\}$ inside
$V_{E_8}$ as obtained in the last three theorems. The arrow 
$\swarrow$ (resp.~$\searrow$) denotes the untwisted (resp.~twisted) 
construction. For a detailed explanation of the second row of the table  
see~[H2]. Self-dual {\em Kleinian codes\/} are a generalization of the so 
called {\em type\/}~IV codes over $\F_4$. Especially,
the notation of a {\em marking\/} of a Kleinian code is defined in~[H2]. 
Finally, $\Xi_1$ is the $D_8^*/D_8$-code $\{0,s\}$ of length $1$ where
the $D_8$-coset $s\in D_8^*/D_8$ has minimal squared length~$2$.

\smallskip

{\small
\noindent\begin{tabular}{lc|ccccccccc}
type & object    &  \multicolumn{9}{c}{marking/frame} \\ \hline
$D_8^*/D_8$-code:   & $\Xi_1$      & & & & &A& & & & \\
                &              & & & & &\dd & & & & \\ 
Kleinian codes: & $\epsilon_2$ & & & &a& &b& & & \\
               &              & & & &\dd& &\dd & & & \\ 
binary codes:$\!\!$ & $H_8$        & & &$\alpha$& &$\beta$ & &$\gamma$ & & \\
               &              & & &\dd & &\dd & &\dd & & \\ 
lattices:      & $E_8$        & &${\cal K}_8$&&${\cal K}'_8$&&${\cal L}_8$&&${\cal O}_8$ \\
               &              & &\dd & &\dd & &\dd & &\dd & \\ 
VOAs:         & $V_{E_8}$    &$\Gamma$ & &$\Sigma$ & &$\Psi$ & &$\Theta$ &&$\Omega$  \\
\end{tabular}
}
\phantom{xxx}

\bigskip

\noindent{\em Example II:
The Golay code ${\cal G}_{24}$, the Leech lattice $\Lambda$, and the moonshine 
module $V^{\natural}$}

\smallskip

The moonshine module $V^{\natural}$ is the $\Z_2$-orbifold \voa
of $V_{\L}$ associated to the Leech lattice $\L$ which is itself the
twisted lattice $\widetilde{L}_{{\cal G}_{24}}$ coming from the Golay code
${\cal G}_{24}$ (cf.~[B], [FLM]). That is, 
$V^{\natural}=\widetilde{V}_{\widetilde{L}_{{\cal G}_{24}}}$. 

\medskip
To describe Virasoro decompositions of the moonshine module coming from 
markings of the Golay code, we must study these markings first.
For the decomposition polynomial $P_{V^\natural}(a,b,c)$ only, it is 
enough to compute the coefficients $W_{k,l}$ of the symmetrized marked 
weight enumerator. The possible values for $W_{8,l}$ (and so for $W_{16,l}$) 
for the Golay code were computed by Conder and McKay in~[CM]. 
They found $90$ possibilities. It is not clear if the numbers
$W_{12,l}$, which are also needed, can be determined from the $W_{8,l}$ alone.

The markings for the Golay code are classified by
the double cosets $\Z_2^{12}.Sym_{12}\backslash Sym_{24}/M_{24}$. (The
first subgroup is the stabilizer of a partition of the 24-set into
2-sets; the second is $M_{24}$, the automorphism group of ${\cal G}_{24}$.)
In fact there are $1858$ different classes of markings~[Be]. 

\smallskip 

The binary linear code ${\cal C}\leq \F_2^{48}$ as defined in 
Section~\ref{one} depends also on the chosen marking. Since for the moonshine
module we have ${\rm dim}\, V^{\natural}_1=0$ the minimal weight of ${\cal C}$
is at least four. The following easy result gives an restriction on
the dimension of ${\cal C}$. 
\begin{lem}\label{sphere-bound}
For every frame of $48$ Virasoro vertex operator algebras of rank~$\ha$
inside the moonshine module the dimension of ${\cal C}$ is
smaller than or equal to $41$. 
\end{lem}
\pf Deleting one coordinate of the codewords of a 
$k$-dimensional code ${\cal C}$ of minimal weight $4$ leads to a 
code of length $47$, dimension $k$ and minimal weight at least $3$. 
Minimal weight $3$ implies that the spheres of radius one around
the codewords of this code are all disjoint, i.e.~we have the sphere 
packing condition $2^k\cdot(1+47)\leq 2^{47}$ or $k\leq 41$. \qed

There is indeed a special marking ${\cal M}^*$ where ${\cal C}$ meets this bound. 
A good way to define it, is to describe the Golay code itself by a 
``double twist'' construction. Starting from the 
glue code $\Xi_3$ of the Niemeier lattice with root sublattice $D_8^3$ one 
gets first the hexacode ${\cal H}_6$, a code over the Kleinian fourgroup, 
and from the hexacode one obtains the Golay code ${\cal G}_{24}$:

\smallskip

As a code over $D_8^{*}/D_8=\Z_2\times \Z_2=\{0,1,s,\bar{s}\}$
(where $1$, $s$, $\bar{s}$ are  the $D_8$-cosets represented by
$(0^7,1)$,  $((\ha)^{7},\pm\ha)$, respectively) one has (cf.~[V])
$$
\Xi_3=\left\{(000),(s11),(1s1),(11s),
(0\bar{s}\bar{s}),(\bar{s}0\bar{s}),
(\bar{s}\bar{s}0),(sss) \right\}.
$$
The hexacode as a code over $D_4^{*}/D_4=\Z_2\times \Z_2=\{0,a,b,c\}$ 
(where $a=[(0,0,0,1)]$, $b=[(\ha,\ha,\ha,\ha)]$ and $c=[(\ha,\ha,\ha,-\ha)]$)
can be defined by
$$
{\cal H}_6=\widetilde{\Xi}_3:=\left(\widehat{\Xi}_3+(\delta_2^3)_0\right)
\cup \left(\widehat{\Xi}_3+(\delta_2^3)_0+(b0b0ca)\right).
$$
Here $\widehat{\phantom{\Xi}}$ is the map induced from $\hat{\phantom{i}}:
D_8^{*}/D_8 \longrightarrow (D_4^{*}/D_4)^2$, $0\mapsto 00$, 
$1\mapsto a0$, $s\mapsto bb$ and $\bar{s}\mapsto cb$, 
and $(\delta_2^n)_0$ is the subcode of 
the Kleinian code $\delta_2^n:=\{(00),(aa)\}^{n}$ of length~$2n$ 
consisting of codewords of weights divisible by $4$. 

\smallskip
In a similar way one gets 
$$ 
{\cal G}_{24}=\widetilde{\cal H}_6:=\left(\widehat{{\cal H}}_6+(d_4^6)_0\right) \cup
\left(\widehat{\cal H}_6+(d_4^6)_0+(1000\ 1000\ \ldots\ 1000\ 0111)\right),
$$
where $\widehat{\phantom{\Xi}}$ is the map induced from $\hat{\phantom{i}}:
D_4^{*}/D_4 \longrightarrow (D_2^{*}/D_2)^2\cong\F_2^4$, 
$0\mapsto 0000$, $a\mapsto 1100$, $b\mapsto 1010$, $c\mapsto 0110$,
and $(d_4^n)_0$ is the subcode of the binary code $d_4^n:=\{(0000),(1111)\}^n$ 
of length $4n$ consisting of codewords of weights divisible by $8$.
This is the usual MOG or hexacode construction of the Golay code
and is a special case of the twisted construction of binary codes from 
Kleinian codes (cf.~[H2], last section).

In this description of the Golay code we let ${\cal M}^*=\{(1,2),\ldots,(47,48)\}$ 
the special marking mentioned above. The marking used in~[DMZ] and~[H1] arose from 
the way the Golay code was written there as a cyclic code.   

The symmetrized marked weight enumerator for the marking ${\cal M}^*$ of the Golay
code is easily computed (using for example the above description)
and one gets
\begin{eqnarray}\label{smwe-golay}
{\rm smwe}_{{\cal G}_{24}}(x,y,z)& = &
{x^{12}}+{z^{12}} + 39\,({x^4}\,{z^8}+ {x^8}\,{z^4})+48\,{x^6}\,{z^6} \\ &&\nonumber
+\left( 96\,({x^6}\,{z^2}+\,{x^2}\,{z^6})+
192\,{x^4}\,{z^4} \right)\,{y^4}\\&&\nonumber  +
\left( 576\,(\,{x^5}\,z  + \,x\,{z^5})+ 1920\,{x^3}\,{z^3} 
\right)\,{y^6} \\&&\nonumber 
+\left( 48\,({x^4}+\,{z^4})  +  288\,{x^2}\,{z^2} \right) \,{y^8}
+128\,{y^{12}}.
\end{eqnarray}

Another property of the marking ${\cal M}^*$ is, that it has the largest 
stabilizer inside $M_{24}$ among all the different markings, namely
$2^6\colon [Sym_4\times Sym_3]$ of order $2^{10}3^2 = 9216$ 
(see Appendix~\ref{golay-code}), as was noted in~[CM]. 

\begin{rem}\label{golaymarkings}{\rm
Assume that a marking is represented by the standard 
partition $\{(1,2),(3,4),\ldots,(23,24)\}$.
The markings of the Golay code that arise from markings of the hexacode
in the sense of Kleinian codes (cf.~end of last subsection) are exactly the ones 
for which the code $(d_4^6)_0$ is a subcode of ${\cal G}_{24}'$,
a code equivalent to ${\cal G}_{24}$. } 
\end{rem}

\medskip
{}From Lemma~\ref{ldecomp}, we get the decomposition of the Leech lattice
$\Lambda\cong\widetilde{L}_{{\cal G}_{24}}$ under the $D_1$-frame belonging 
to the marking ${\cal M}^*$. For the symmetrized weight enumerator of the 
corresponding code $\widetilde{\Gamma}\leq \Z_4^{24}$ 
(see (\ref{gamma-def})). Corollary~\ref{ldecomp} gives:

\vspace{-2mm}
{\small
\begin{eqnarray*}
{\rm swe}_{\widetilde{\Gamma}}(A,B,C)&=& 
   {A^{24}} +  {{{\rm C}}^{24}}+
   23439\,({A^{16}}\,{{{\rm C}}^8} + {A^8}\,{{{\rm C}}^{16}}) +
   4032\,({A^6}\,{{{\rm C}}^{18}} + {A^{18}}\,{{{\rm C}}^6})  
\\ && +
   378\,({A^4}\,{{{\rm C}}^{20}} + \,{A^{20}}\,{{{\rm C}}^4}) +
   60480\,({A^{10}}\,{{{\rm C}}^{14}} + {A^{14}}\,{{{\rm C}}^{10}} )+
   85484\,{A^{12}}\,{{{\rm C}}^{12}}  
\\ && + \Big(
   3072\,({A^2}\,{{{\rm C}}^{14}} + {A^{14}}\,{{{\rm C}}^2}) + 
   43008\,({A^{12}}\,{{{\rm C}}^4} + {A^4}\,{{{\rm C}}^{12}})  
\\ && +
   193536\,({A^{10}}\,{{{\rm C}}^6} + {A^6}\,{{{\rm C}}^{10}}) + 
   307200\,{A^8}\,{{{\rm C}}^8}
   \Big) \,{B^8}
\\ && + \Big(
   86016\,({A^{11}}\,{\rm C} + A\,{{{\rm C}}^{11}}) + 
   1576960\,({A^9}\,{{{\rm C}}^3} +{A^3}\,{{{\rm C}}^9}) 
\\ && +
   5677056\,({A^7}\,{{{\rm C}}^5} + {A^5}\,{{{\rm C}}^7})  
   \Big)\,{B^{12}}
\\ && + \left(
   6144\,({A^8}\, + {{{\rm C}}^8}) + 
   172032\,({A^6}\,{{{\rm C}}^2} + \,{A^2}\,{{{\rm C}}^6}) + 
   430080\,{A^4}\,{{{\rm C}}^4} 
   \right)\,{B^{16}}\  
\\ && + 
   262144\,{B^{24}}.
\end{eqnarray*}
}
\vspace{-2mm}

As stated before, the markings for the Golay code are classified by
the double cosets $\Z_2^{12}.Sym_{12}\backslash Sym_{24}/M_{24}$.  

The classification of all $D_1$-frames in the Leech lattice
would seem to be more complicated. From equation~(\ref{gamma-def}), 
we see that in the case where the 
$D_1$-frame comes from a marking of the Golay code the corresponding
$\Z_4$-code $\widetilde{\Gamma}$ contains the subcode $(\Sigma_2^{12})_0$.
The following result gives the converse. Recall that the {\em Euclidean weight\/}
of a codeword is the minimal Euclidean squared norm of a coset representative in
$(D_1^*)^{24}$. 

\begin{lem}\label{leechframes}
Every self-annihilating
even $\Z_4$-code $\Delta$ of length $24$ and minimal Euclidean
weight $4$ containing the subcode $(\Sigma_2^{12})_0$
can be obtained from a marking of the Golay code as in equation~(\ref{gamma-def}).
\end{lem}

\pf Let $K=\bigoplus_{i=1}^{24}\Z\,a_i$ a lattice of type $A_1^{24}$ in 
$\R^{24}$, i.e.~the $a_i$ are pairwise orthogonal vectors of squared
length $2$. Set $L=\bigoplus_{i=1}^{24}\Z\,b_i$, with $b_{2i-1}=
a_{2i-1}+a_{2i}$ and $b_{2i}=a_{2i-1}-a_{2i}$ for $i=1$, $\ldots$, $12$,
i.e.~$L$ is a lattice of type $D_1^{24}$. Finally let $M=2K$.
On $K$, the group $\Z_2^{24}\colon Sym_{24}$ acts by monomial
matrices with entries $\pm 1$ with respect to the basis 
$\{ a_i  \mid  i=1, \ldots, 24 \}$. 
The lattice $L$ is fixed at least by the group of sign changes.
Clearly $K^*/K\cong\Z_2^{24}$, $L^*/L\cong\Z_4^{24}$, $M^*/M\cong\Z_8^{24}$,
with the induced action of  $\Z_2^{24}\colon Sym_{24}$ on $\Z_2^{24}$ 
and $\Z_8^{24}$ and of $\Z_2^{24}$ on $K$.

The code $\Delta\leq L^*/L$ determines a self-dual even  
lattice $\Lambda$ of rank $24$ and minimal length $4$.
(This must be the Leech lattice since it is the unique self-dual 
even rank $24$ lattice of minimal length $4$.)

To prove the lemma we have to find a doubly-even self-annihilating
binary code ${{\cal G}}_{24}'\leq K^*/K$ equivalent to the Golay code 
${\cal G}_{24}$ such that ${{\cal G}}_{24}'$ determines 
$\Delta=\widetilde{\Gamma}$ as in~(\ref{gamma-def}). 
(Instead of changing the marking ${\cal M}=\{(1,2),\ldots,(23,24)\}$, 
the choice which is determined by the relation between $K$ and $L$, 
we are permuting the code ${{\cal G}}_{24}$; these procedures are equivalent.) 

The lattice $\Lambda$ defines a self-annihilating even
$\Z_8$-code $\Omega=\Lambda/M\leq M^*/M$ of minimal Euclidean weight~$4$. 
If we start with the with our standart copy of the Golay code 
${{\cal G}}_{24}$ we get a lattice $\widetilde{\Lambda}$, 
a $\Z_4$-code $\widetilde{\Delta}\subset L^*/L$, and a $\Z_8$-code 
$\widetilde{\Omega}\subset M^*/M$.

Since $(\Sigma_2^{12})_0$ is contained in $\Delta$,  we see easily that the code 
$\Omega$ contains all $24 \choose 2$ vectors of type $(4^20^{22})$.
As a main step in the uniqueness proof of $\Lambda$ in [Co], it was shown that 
such a code is unique up to the action of $\Z_2^{24}\colon Sym_{24}$, i.e.~we 
have a $\pi$ in this group such that $\pi(\widetilde{\Lambda})/M=\Sigma$.
The copy ${{\cal G}}_{24}'=\pi({{\cal G}}_{24})$ of the Golay code gives the code
$\Delta$ in $L^*/L$. \qed

\medskip
\noindent{}Finally we come to the Virasoro decomposition of the moonshine 
module \hbox{$V^{\natural}=\widetilde{V}_{\widetilde{L}_{{\cal G}_{24}}}$.}

The following theorem gives an precise description of the codes ${\cal C}$ and 
${\cal D}$ as defined in section~\ref{one}.

\begin{thm}\label{48code} The code ${\cal C}$ 
associated to the special marking ${\cal M}^*$ of the Golay code
has length $48$ and dimension $41$. Its annihilator code 
${\cal C}^{\perp}=\{d\in\F_2^{48}\mid (d,c)=0\ \hbox{for all\ }
c\in {\cal C}\}$ is of dimension $7$ and equals the code
${\cal D}$ which has generator matrix
$$\left(\begin{array}{ccc}
1111 1111 1111 1111 & 0000 0000 0000 0000 & 0000 0000 0000 0000 \\
0000 0000 0000 0000 & 1111 1111 1111 1111 & 0000 0000 0000 0000 \\
0000 0000 0000 0000 & 0000 0000 0000 0000 & 1111 1111 1111 1111 \\ 
0000 0000 1111 1111 & 0000 0000 1111 1111 & 0000 0000 1111 1111 \\
0000 1111 0000 1111 & 0000 1111 0000 1111 & 0000 1111 0000 1111 \\
0011 0011 0011 0011 & 0011 0011 0011 0011 & 0011 0011 0011 0011 \\
0101 0101 0101 0101 & 0101 0101 0101 0101 & 0101 0101 0101 0101 
\end{array}\right)_{\textstyle .}$$
\end{thm}

\pf Recall the description of the Golay code given above.
The codes ${\cal H}_6$ and ${\cal G}_{24}$ are unions of two parts.  
The first part we call the untwisted part and the second is called the
twisted part. 

First we show that the above matrix is a parity check matrix for 
${\cal C}$. From Theorem~\ref{decomp} we see that a 
codeword $c\in {\cal G}_{24}$ 
gives us an irreducible
$T_{48}$-module $M(h_1,\ldots,h_{48})$ with all $h_i$ 
different from 
$\frac{1}{16}$ if and only if $c(k)\in \{(0,0),(1,1)\}$ for all
$k$. The codewords with this property are exactly the ones that are 
coming from the codeword $(000)\in\Xi_3$. This gives the first three rows
of the parity check matrix. The next two rows correspond to the selection
of the subcodes $(\delta_2^3)_0\subset \delta_2^3$ and 
$(d_4^6)_0\subset d_4^6$.
Let ${\cal B}_2^n$ the FVOA $(M(0,0)\oplus M(\ha,\ha))^{\otimes n}$
with binary code ${\cal C}({\cal B}_2^n)=\{(0,0),(1,1)\}^n$ of length $2n$. 
The subVOA $({\cal B}_2^n)_0$ is the FVOA belonging to the subcode of 
${\cal C}({\cal B}_2^n)$ consisting of codewords of weights divisible 
by $4$ (cf.~Proposition~\ref{C-voas}). 
Then the last two rows of the parity check matrix
correspond to the selection of the subcodes
$(\Sigma_2^{12})_0\subset \Sigma_2^{12}$
and ${\cal C}(({\cal B}_2^{24})_0)\subset {\cal C}( {\cal B}_2^{24} )$:
this amounts to the conditions $\prod \epsilon_k=+$ and 
$\prod \mu_k=+$. There are no further conditions.

To determine ${\cal D}$ note first that the inclusion 
${\cal D}\leq {\cal C}^{\perp}$ is Proposition~\ref{codes}~(3).
To see ${\cal C}^{\perp}\leq  {\cal D}$ observe that the codewords 
$\{(s11),(1s1),(11s)\}\subset\Xi_3$ 
correspond to the first three lines of the generator matrix, 
the twisted parts of ${\cal H}_6$ and ${\cal G}_{24}$
to the next two, and two of the last three summands of 
$V^{\natural}=\widetilde{V}_{\widetilde{L}_{{\cal G}_{24}}}$ 
in Theorem~\ref{decomp} correspond to the last two lines 
of the generator matrix. 

Alternatively,  one can compute ${\cal D}$ by using the self-duality of the 
moonshine VOA~[D3] and apply Theorem~\ref{codeholo}. \qed

The code ${\cal C}$ is also the lexicographic code of 
length $48$ and minimal weight $4$ (see [CS3],~Th.~6). 
As mentioned there, it is a ``shortened extended Hamming code" of 
length $64$ in the following sense: 
If we extend the generator matrix of ${\cal D}$ by
the block
$$\left(\begin{array}{c}
1111 1111 1111 1111 \\
1111 1111 1111 1111  \\
1111 1111 1111 1111  \\ 
0000 0000 1111 1111  \\
0000 1111 0000 1111  \\
0011 0011 0011 0011  \\
0101 0101 0101 0101  \\
\end{array}\right)_,$$
we obtain a parity check matrix for the extended Hamming code $H_{64}$ of 
length $64$. The vectors $c\in \F_2^{64}$ with $0$'s in the last~$16$
coordinates belong to $H_{64}$ if and only if
the vector of the first $48$ coordinates belongs to ${\cal C}$.

The automorphism group of this code is of type 
$2^{12}[GL(4,2)\times Sym_3]$ and has
order $495452160$ (see Appendix~\ref{moon-code} for a proof).

\smallskip
For future references we give the decomposition polynomial as obtained 
{}from Corollary~\ref{decomp-pol} in full. Remember that $a$, $b$ and $c$ count
the modules of conformal weight $0$, $\frac{1}{2}$, resp.~$\frac{1}{16}$ 
(see Definition~\ref{dmz}).
\begin{cor}
The complete decomposition polynomial for the moonshine module belonging to 
the special marking ${\cal M}^*$ is given by
\begin{eqnarray*}
& &P_{V^{\natural}}^{{\cal M}^*}(a,b,c)=
   {a^{48}}+{b^{48}} + 
   3300\,({a^{44}}\,{b^4}+{a^{4}}\,{b^{44}})+ 
   189504\,( {a^{42}}\,{b^6} + {a^{6}}\,{b^{42}})\\
& &\ \ \ \ \ + 
   5907810\,({a^{40}}\,{b^8} + {a^{8}}\,{b^{40}} )+
   102156864\,({a^{38}}\,{b^{10}} + {a^{10}}\,{b^{38}} )\\
& &\ \ \ \ \ +
   1088684372\,({a^{36}}\,{b^{12}} + {a^{12}}\,{b^{36}})+
   7535996160\,({a^{34}}\,{b^{14}} + {a^{14}}\,{b^{34}})\\
& &\ \ \ \ \ +
   35232581487\,({a^{32}}\,{b^{16}} + {a^{16}}\,{b^{32}})+
   114215080192\,({a^{30}}\,{b^{18}} + {a^{18}}\,{b^{30}})\\
& &\ \ \ \ \ +
   261496913352\,({a^{28}}\,{b^{20}} + {a^{20}}\,{b^{28}})+ 
   427898196864\,({a^{26}}\,{b^{22}} +{a^{22}}\,{b^{26}})\\
& &\ \ \ \ \ + 
   503871835740\,{a^{24}}\,{b^{24}}  \\
& &\ \ \ \ \ +\left(
   6144\,({a^{30}}\,{b^2}+{a^{2}}\,{b^{30}}) + 
   430080\,({a^{28}}\,{b^4}+{a^{6}}\,{b^{28}})\right.\\
& &\ \ \ \ \ +
   10881024\,({a^{26}}\,{b^6}+{a^{6}}\,{b^{26}}) + 
   126197760\,({a^{24}}\,{b^8}+{a^{8}}\,{b^{24}}) \\
& &\ \ \ \ \ + 
   774199296\,({a^{22}}\,{b^{10}}+{a^{10}}\,{b^{22}}) + 
   2709417984\,({a^{20}}\,{b^{12}}+{a^{12}}\,{b^{20}})\\
& &\ \ \ \ \ \left.+ 
   5657364480\,({a^{18}}\,{b^{14}}+{a^{14}}\,{b^{18}}) + 
   7212810240\,{a^{16}}\,{b^{16}}
   \right) c^{16}\\
& &\ \ \ \ \ + \left(
   184320\,({a^{23}}\,b+{a}\,{b^{23}} ) + 
   15544320\,{a^{21}}\,{b^3}+{a^{3}}\,{b^{21}} )\right.\\
& &\ \ \ \ \ + 
   326430720\,({a^{19}}\,{b^5}+{a^{5}}\,{b^{19}} ) + 
   2658078720\,({a^{17}}\,{b^7}+{a^{7}}\,{b^{17}} )\\
& &\ \ \ \ \  \left.+ 
   10041630720\,({a^{15}}\,{b^9}+{a^{9}}\,{b^{15}} ) + 
   19170385920\,({a^{13}}\,{b^{11}}+{a^{11}}\,{b^{13}} )
   \right){c^{24}}\\
& &\ \ \ \ \  + \left(
   3072\,({a^{16}})+{b^{16}})  + 
   368640\,({a^{14}}\,{b^2}+{a^2}\,{b^{14}}  )\right.\\
& &\ \ \ \ \ + 
   5591040\,({a^{12}}\,{b^4}+{a^4}\,{b^{12}} )+
   24600576\,({a^{10}}\,{b^6}+{a^6}\,{b^{10}} )\\
& &\ \ \ \ \ \left.+ 
   39536640\,({a^8}\,{b^8})   \right)\,{c^{32}} + 
   131072\,{c^{48}}.
\end{eqnarray*}
\end{cor}
It was shown in Chapter~4 of~[H1] that for a self-dual 
vertex operator algebra $V$ the decomposition
polynomial belongs to the ring $\C[a,b,c]^G$ of invariants for some
$3\times 3$-matrix group $G$ of order $1152$. The space of invariant
homogeneous polynomials of degree $48$ is $7$-dimensional and it can 
be checked that the above polynomial indeed belongs to this space
by using the explicit base given in~[H1].

We expect that the analog of Remark~\ref{golaymarkings} 
and Lemma~\ref{leechframes} holds: 
Every self-dual FVOA of central charge $24$ and minimal weight $2$ 
(i.e.~${\rm dim\,}V_1=0$) containing the subVOA $({\cal B}_2^{24})_0$
can be obtained from a $D_1$-frame of the Leech lattice as in
the second equation of Theorem~\ref{d1-decomp}. 
\vspace{1 cm}


\appendix
\newtheorem{thma}{Theorem}[subsection]
\newtheorem{propa}[thma]{Proposition}
\newtheorem{lema}[thma]{Lemma}
\newtheorem{rema}[thma]{Remark}
\newtheorem{cora}[thma]{Corollary}
\newtheorem{conja}[thma]{Conjecture}
\newtheorem{dea}[thma]{Definition}
\newtheorem{nota}[thma]{Notation}

{}{\Large\bf{}Appendix}
\section{Orbits on markings of a Hamming Code}\label{hamming-code}

\begin{notation}\rm   
Let $H$ be the unique binary code with parameters $[8,4,4]$, the
Hamming code (see  Appendix~\ref{moon-code}).  We take it to be the span of 
$(00001111)$, $(00110011)$, $(01010101)$, $(11111111)$ (see (\ref{A.3.4})).  
Let $A:= Aut(H) \cong AGL(3,2)$ (\ref{A.3.5}).  A {\em marking\/} is a 
partition of the index set into $2$-sets.   
\end{notation}

The number of markings is 
${8 \choose 2}{6 \choose 2}{4 \choose 2}{2 \choose 2}/4!=2520/24=105$.  
We show that there are 
three  orbits of $A$ on the set of markings and determine the 
stabilizers.  
This group is triply but not quadruply transitive on the eight indices.

\begin{notation}\rm
It helps to interpret the index set 
as $V \cong \Bbb F_2^3$ with the obvious action of $A$. 
So, $2$-sets correspond to affine subspaces of dimension 1.  
Take  a linear subspace $U \le
V$ of dimension 1. 
Let $T$ the translation subgroup of $A$ and let 
$L:=Stab_A(0) \cong GL(3,2)$.
Let $M$ be a marking, $S:= Stab_A(M)$.   
By double transitivity, we may assume $U\in M$.   Let $R \le T$ be
the group of order 2 corresponding to $U$.   
\end{notation} 

\smallskip

Case $\alpha$.  We assume that all four parts of the marking $M$ are cosets of
$U$.  Then $S=TStab_L(U) \cong 2^3\colon Sym_4$, a group of index~$7$ in $A$.  

\smallskip

Case $\beta$.  
We assume that  exactly two parts of the marking are cosets of 
$U$, say $U$ and $W$.  Let $P$ and $Q$ be the other two parts.  
Then $X:=U \cup W$ is a dimension 2 linear subspace of $V$ and $Y:=P \cup
Q$  is its complement.  Both $P$ and $Q$ are cosets of a common linear
$1$-dimensional subspace $U^* \ne U$ of $X$.  

Let $R^*$  be the fours group in $T$ which corresponds to $X$;  
$R^* > R$.  Then,  $R^*$ stabilizes both $\{U,W\}$ and $\{P,Q\}$, 
whence  $R^* \le S$; in fact, $R^*=T \cap S$.  

Since $A$ acts transitively on pairs of parallel affine 1-spaces, 
$S$ acts transitively on $\{X,Y\}$; let $S_0:=Stab_S(X)=Stab_S(Y)$.  
Then $S_0$ has index 2 in $S$ and acts transitively on $\{U,W\}$; let 
$S_1$ be the common stabilizer, index 2 in $S_0$.  
Also,  $S_0$ acts transitively on $\{P,Q\}$; 
let $S_2$ be the common stabilizer,
index 2 in $S_0$.   Then $S_1 \ne S_2$  
(since $R$ stabilizes $U$ and $W$ but
interchanges $P$ and $Q$),  and $S_4:=S_1 \cap S_2 \triangleleft S$ and
$S/S_4 \cong Dih_8$, a Sylow $2$-group of $Sym_4$ 
(via its action on the marking).   

It suffices to show that $|S_4|=4$.  Clearly, elements of $S_4$ have 
square $1$. The involution which is trivial on $X$ and interchanges 
the points within each of $P$ and $Q$ is in $L$.  
The same idea, with $U$, $W$ replaced by $P$, $Q$ gives
an involution which is in a conjugate of $L$, say in $L^g$, where $g \in A$
interchanges $X$ and $Y$.  Since these involutions are different, $|S_4| 
\ge 3$. If $1 \ne u \in  S_4$ has a fixed point, say $v \in V$, it may be 
interpreted as a linear transformation by taking $v$ as the origin; 
since $u$ is an involution, its
fixed point subspace has dimension $2$, and is a union of members of $M$, 
so is one of $X$ or $Y$; this means $u$ is one of the two involutions already
defined.  Therefore, $|S_4| \le 4$, whence equality.   
 
\smallskip

Case $\gamma$.  
We assume that all parts of the marking besides $U$ are not 
cosets of $U$.  
It follows that $S \cap T=1$, so $S$ embeds in $L$.  Clearly, $7$ does not
divide $|S|$, so $S$ embeds as a proper subgroup of order dividing $24$. 
Thus, the orbit here has length divisible by $8\cdot 7=56$.  By our  above 
count of the
number of markings, this must be the exact number.  We conclude that 
$S \cong Sym_4$, since the only subgroups of odd index in $GL(3,2)$ are
parabolic subgroups [Ca],~8.3.2.    


\section{Automorphisms of a marked Golay code}\label{golay-code}

\def\la{\langle}
\def\ra{\rangle}
\def\Euc{$\RR^{\ell}$ }

We settle the stabilizer in $M_{24}$ of the special marking $\cal M^*$
we obtained in our description of the Golay code and identify $\cal M^*$ with the
exceptional marking of Blackburn, Conder and McKay~[CM] with parameters
$(48,576,96,0,39)$. 

\smallskip
As noted in section~5 our construction of $\cal G$ is equivalent to 
the usual hexacode construction, as in [G2] (5.25). 
The marking ${\cal M^*}$ in this notation
is gotten from the usual sextet partition of 
the $24$-set $\Omega$

$$\begin{array}{cccccc}
0\ \ \ \ &\bullet  \ \ \ \ \ \bullet\ 
& \ & \bullet  \ \ \ \ \ \bullet\ 
& \ &\bullet  \ \ \ \ \ \bullet\cr
1\ \ \ \ &\bullet  \ \ \ \ \ \bullet\ 
& \ & \bullet  \ \ \ \ \ \bullet\ 
& \ &\bullet  \ \ \ \ \ \bullet\cr
\omega\ \ \ \ &\bullet  \ \ \ \ \ \bullet\ 
& \ & \bullet  \ \ \ \ \ \bullet\ 
& \ &\bullet  \ \ \ \ \ \bullet\cr
\bar \omega\ \ \ \ &\bullet  \ \ \ \ \ \bullet\ 
& \ & \bullet  \ \ \ \ \ \bullet\ 
& \ &\bullet  \ \ \ \ \ \bullet\cr
\end{array}
$$
by intersecting  the columns with the unions $Row_0 \cup Row_1$ and
$Row_\omega \cup Row_{\bar \omega}$:  

$$\begin{array}{cccccc}
0\ \ \ \ &\bullet  \ \ \ \ \ \bullet\ 
& \ & \bullet  \ \ \ \ \ \bullet\ 
& \ &\bullet  \ \ \ \ \ \bullet\cr
1\ \ \ \ &\bullet  \ \ \ \ \ \bullet\ 
& \ & \bullet  \ \ \ \ \ \bullet\ 
& \ &\bullet  \ \ \ \ \ \bullet\cr
\cr 
\omega\ \ \ \ &\bullet  \ \ \ \ \ \bullet\ 
& \ & \bullet  \ \ \ \ \ \bullet\ 
& \ &\bullet  \ \ \ \ \ \bullet\cr
\bar \omega\ \ \ \ &\bullet  \ \ \ \ \ \bullet\ 
& \ & \bullet  \ \ \ \ \ \bullet\ 
& \ &\bullet  \ \ \ \ \ \bullet\cr
\end{array}
$$
The set of twelve resulting $2$-sets form $\cal M$.  

In this appendix, we settle the stabilizer in $M_{24}$ of the
marking above. 

\smallskip

In~[G2], the action of the associated sextet group on $\Omega$ is 
described. 
The group has shape $H:=2^63\cdot Sym_6$ and may be thought of as
${\cal H}\colon { Aut}^*({\cal H})$, the affine hexacode group ($\cal H$
denotes the hexacode)~(5.25).  As in~[G2], Chapters~5 and~6, we use
the notation $K_i$  for  the $4$-set in $\Omega$ occurring as the $i^{\rm th}$ 
column above and $K_{ij\dots}$ denotes the union $K_i \cup K_j \cup \cdots$.  

\smallskip

The obvious subgroup of $H$ which preserves $\cal M$ is ${\cal H}:P$,
where $P=S \times \la t \ra \cong Sym_4 \times 2$, where $S$ is generated
by the groups of permutations  (1)  the four-group of row-respecting
column permutations which interchange columns within evenly many
coordinate blocks $K_{12}$, $K_{34}, K_{56}$;  (2) the copy of $Sym_3$ 
obtained by permuting the three coordinate blocks 
(respecting the order within the blocks);  (3)  the permutation $t$ 
is given by the diagram
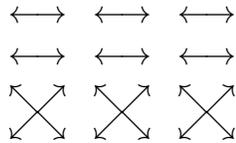
\vbox{  
$$\begin{array}{ccc}
\leftarrow\! \rightarrow&\leftarrow\! \rightarrow&\leftarrow\! \rightarrow\cr
\leftarrow \!\rightarrow&\leftarrow\! \rightarrow&\leftarrow\! \rightarrow\cr
{\nwarrow\! \nearrow} & {{\nwarrow\! \nearrow}} & {\nwarrow\! \nearrow}
\end{array}$$
\vspace{-6.5mm}
$$\begin{array}{ccc}
 {\swarrow\! \searrow} &
 {\swarrow\! \searrow} &
 {\swarrow\! \searrow} 
\end{array}
$$}
[G2]~(5.38),~UP2.   

\smallskip

The corresponding subgroup $Sym_4 \times 2$ of $Sym_6$ is maximal (since
it is the stabilizer of a $2$-set in a sextuply transitive action).  Since 
the  ``scalar'' transformation UP9~(5.38)~[G2] (fixes top row elementwise,
cycles rows $2$, $3$, $4$ downward) 
$$\begin{array}{cccccc}
\bullet & \bullet & \bullet & \bullet & \bullet & \bullet \cr
\downarrow & \downarrow &\downarrow &\downarrow &\downarrow
&\downarrow  \end{array}
$$
does not stabilize $\cal M$, it follows that $\cal H:P$ is the 
stabilizer of $\cal M$ in $H$.  

\begin{notation}\label{A.2.1}\rm   
$R:= Stab_G(\cal M)$,  $G:=M_{24}$.  
\end{notation}  

The next step is to determine $R$; we know that $R \cap H ={\cal H}\colon P$.  

We take a clue from the symmetrized marked 
weight enumerator ${\rm smwe}_{\cal G}(x,y,z)$ of the Golay code as given 
in~(\ref{smwe-golay}) and see that the parameters in the sense of~[CM] are 
$(48,576,96,0,39)$. The next result is an exercise.

\begin{lem}\label{A.2.2}  (i) The octads which contribute to 
contribute to $c_0=48$ are those with even parity and which are labeled by
a hexacode word of the form $(00xxxx)$, where $x=\omega$ or $\bar
\omega$ and where the zeroes occur in any of the three coordinate blocks; 
these octads are unions of $2$-sets 
 which  are subsets of columns labeled by $x$.  

(ii)  The octads which contribute to 
contribute to $c_4=39$ have even parity and 
are one of  $K_{ij}$ ($15$ of these) or are  
octads labeled by hexacode words $(001111)$, with the zeroes 
occurring in any of the three coordinate blocks; these octads are unions 
of parts
of $\cal M$ which occur in columns labeled by $1$ ($24$ such octads).  
\end{lem}

Clearly, $R$ permutes the sets of octads  (i) and (ii). 

\begin{lem}\label{A.2.3}   
The orbits of ${\cal H}:P$ on $X$, the set of octads in (ii),  are the
following: 

(a) $K_{12}$, $K_{34}$, $K_{56}$ (length 3); 

(b) $K_{ij}$, for all $\{ij\} \not = \{12\}, \{34\}, \{56\}$ 
(length $12$); 

(c) octads labeled by some $(001111)$ (length $24$).
\end{lem}

\begin{thm}\label{A.2.4}  
$R$ is a subgroup of index 7 in the stabilizer of the trio 
(a), whence $|R:R \cap H|=3$ and $R \cong 2^6\colon [Sym_4 \times Sym_3]$. 
\end{thm}

\pf
We consider the action of $R$ on $X$.  
The octads in $X$ which have only a $0$- or a $4$-set as intersection with all
members of $X$ are the three in (a).  So, $R$ preserves this trio and so 
is in the
trio group, $J$, of the form $2^6[GL(3,2) \times Sym_3]$.  
The group  ${\cal H}:P$ is a  subgroup of $J$ of index $21$.
We consider the possibility that $7$ divides $|R|$.  
Let $g \in R$ be an element of order $7$. Then $g$ fixes
at least $1$ of the remaining $36$ members of $X$.  An element of order $7$ in
$G$ fixes exactly three octads and clearly these are just the octads of
our trio (a), a contradiction.   So, $R$ has order $2^{10}3$ or 
$2^{10}3^2$.  
We eliminate the former by exhibiting a permutation in 
$R \setminus H$;  UP13 from (5.38)~[G2] does the job.
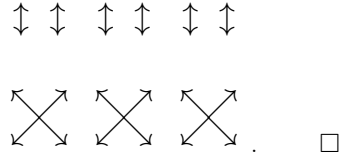
\vbox{
$$\begin{array}{ccc}
\updownarrow \  \ \updownarrow &
\updownarrow \  \ \updownarrow &
\updownarrow  \ \ \updownarrow  \cr
&&\cr 
\nwarrow\! \nearrow   &
\nwarrow\! \nearrow   &
\nwarrow\! \nearrow  
\end{array}_{\phantom{.\qquad\qed}}$$
\vspace{-6.6mm}
$$\begin{array}{ccc}
\swarrow\! \searrow &
\swarrow\! \searrow &
\swarrow\! \searrow 
\end{array}_{.\qquad\qed}
$$}

\medskip
Finally we can identify our marking ${\cal M}^*$ with the one in~[CM].
This is not completely obvious since the labelings chosen in [CM] are
different from the standard ones, e.g.~in [Atlas] or~[G2].

There are $|G|/|R|=26565$ markings equivalent to
${\cal M}^*$, but this is exactly the number of markings obtained
in~[CM] with parameters $(48,576,96,0,39)$, i.e.~there is only one
orbit of markings with this parameters.   


\section{Automorphism group of certain codes ${\cal C}$ and 
${\cal D}$ of length~$3\cdot2^d$}\label{moon-code} 

\def\FF{\F}
\def\la{\langle}
\def\ra{\rangle} 

We are studying binary codes ${\cal C} \le \FF_2^\Omega$, where
$|\Omega|=3\cdot 2^d$ and ${\cal D}:={\cal C}^\perp$ is spanned by the 
$d+3$ rows of the matrix
\begin{equation}\label{A.3.1x}
M := \left(\begin{array}{ccc}
1111 \dots\, \dots 1111  &  0000 \dots\, \dots 0000  &  0000 \dots\, \dots 0000 \cr 
0000 \dots\, \dots 0000  &  1111 \dots\, \ldots 1111  &  0000 \dots\,  \dots 0000 \cr   
0000 \dots\, \dots 0000  &  0000 \dots\, \dots 0000  &  1111 \dots\, \dots 1111 \cr  
00 \ldots 0011 \ldots 11  & 00 \ldots 0011 \ldots 11   & 00 \ldots 0011 \ldots 11  \cr 
\vdots  & \vdots   &  \vdots  \cr
0011 \ldots\, \dots 0011 &  0011 \ldots\, \dots 0011  & 0011 \ldots\, \dots 0011  \cr
0101 \ldots\, \dots 0101 & 0101 \ldots\, \dots 0101 & 0101 \ldots\, \dots 0101 \cr 
\end{array}\right)_{\textstyle .} 
\end{equation}

Our problem is to find $F:=Aut({\cal C})=Aut({\cal D}) \le Sym_\Omega$.  

\begin{notation}\label{A.3.2}  We partition $\Omega$ into three coordinate blocks
$\Gamma_1:=\{1,2,\dots,2^d \}$,  
$\Gamma_3:=\{ 2^d+1, 2^d+2,\dots ,2\cdot 2^d\}$
$\Gamma_3:=\{2\cdot 2^d+1,2\cdot 2^d+2,\dots ,3\cdot 2^d\}$.  
\end{notation}

Here is our main result; it was referred to after Theorem~\ref{48code},
for $d=4.$   

\begin{thm}\label{A.3.3}   
$F \cong 2^{3d}[GL(d,2) \times Sym_3]$, where the $2^{3d}$ may
be interpreted as a tensor product of a $d$ and $3$ dimensional
module for the factors of $GL(d,2) \times Sym_3$.  
\end{thm}

The two main parts of the proof consist of showing that $F$ preserves 
the partition $\{\Gamma_i\}$ then getting the
automorphism groups of the related length $2^d$ codes.

\begin{thm}\label{A.3.4}  There is a unique binary code  of length
$2^d$, dimension $d$ in which all nonzero weights are $2^{d-1}$.  It may be
viewed as a subspace of $\FF_2^V$, where $V:=\FF_2^d$, consisting of $0$, 
the all ones vector, and all affine $(d-1)$-dimensional subspaces of $V$.  
\end{thm}

Existence follows from checking that the set of vectors described in
(\ref{A.3.4}) forms a code with the right properties.    

\begin{lem}\label{A.3.5}      
Assume that $J$ is a code of length $2^d$,
dimension $d$ in which all nonzero weights are $2^{d-1}$.
If $x_1$, $\dots$, $x_m$ are linearly independent
elements of $J$, then $|x_1 \cap \dots \cap x_m|=2^{d-m}$.
\end{lem}  

\pf   The statement is easy to check for $m \le 2$, so we may assume that 
$m\ge 3$. Let $Y := x_3 \cap \dots \cap x_m$ and set $A := x_1 \cap Y$, 
$B := x_2 \cap Y$,  $ C:= (x_1+x_2) \cap Y = A+B$.  Then, by induction
$|A|=|B|=|C|=2^{d-m+1}$.  Since   $A \cap B=|x_1 \cap \dots \cap
x_m|$ and $|A|+|B|-2|A \cap B|=|C|$, we complete the induction step.  
\qed

\begin{prop}\label{A.3.6} 
Assume that $J$ is a code of length $2^d$, dimension $d$ in
which all nonzero weights are $2^{d-1}$.   Then any basis of $J$ is 
equivalent (by coordinate permutations) to the following set of 
vectors of length $2^d$:
$$\begin{array}{c}
 0000 \dots 0000  \ 1111 \dots 1111    \cr 
 0000 \dots 1111 \ 0000 \dots 1111   \cr 
\vdots   \cr  
 0011 \dots 0011 \ 0011 \dots 0011   \cr 
 0101 \dots 0101 \ 0101 \dots 0101   \cr 
\end{array}_{\textstyle .}  
$$
\end{prop}

\pf Use the preceding lemma, permuting coordinates at each stage to right 
justify all the ones and leave earlier basis vectors unchanged.  
\qed 

At once, Theorem~(\ref{A.3.4})  follows and we get a 
closely related result.

\begin{thm}\label{A.3.7}   
(i) The code $J$ spanned by the  $d$ vectors in (\ref{A.3.6}) has
automorphism group $GL(d,2)$.  

(ii) The code spanned by $J$ and the all ones vector has automorphism 
group $AGL(d,2)$; the normal translation subgroup is the group of 
automorphisms which are trivial modulo the span of the all ones vector.
\end{thm}  

\pf
(i)  By the uniqueness, $J$ may be identified with a subspace of the 
power set
$\FF_2^V$ on $V:=\FF_2^d$ consisting of all complements of $d-1$
dimensional linear subspaces; the leftmost index corresponds to $0 \in 
V$. 
Clearly, $GL(V)$ acts as automorphisms.  A  $3$-set $T$ of nonzero vectors sums
to zero in $V$ iff there is an intersection $I=x_1 \cap \dots \cap x_{d-2}$
(cardinality $2^{d-(d-2)}=4$) so that $T= \{ x+y  \mid x \ne y, x,y \in 
I  \}$.  This
characterization of $T$ implies that any automorphism  of $J$ preserves 
addition on $V$, so is in our copy of $GL(V)$.   

(ii)  Let $z=(111 \dots 1) \in \FF_2^V$ be the all ones vector and let
$K:=J+\FF_2 z$.   Adding $z$ to a vector replaces it by its complement (as
usual in binary codes, we identify a vector with its support).   Clearly, 
$z$ is
stable under $G:= Aut(K)$ and, since $J \cong K/\FF_2$ we have a map of
$G$ onto $Aut(J) \cong GL(J) \cong GL(V)$.  
Note that $Aut(K)$ is faithful on $K$; 
this follows from  (\ref{A.3.5}).   The kernel is contained in a subgroup of
$Aut(K)$ naturally identified with  a subgroup of $Hom(K/\FF_2 z,\FF_2 z)
\cong 2^d$.  It suffices to show that the translations in $AGL(d,2)$ act as
automorphisms of $K$.  But this is clear, for  elements of $K \setminus 
\FF_2z$ may be thought of as affine $(d-1)$- dimensional subspaces.   \qed

\begin{notation}\label{A.3.8}\rm
Let ${\cal R}$ be the span of the first three rows of $M$ and 
${\cal S}$ the span of the last $d$. 
Note that the projections of ${\cal S}$ or ${\cal D}$  to any
$\Gamma_i$ block is a code as described by (\ref{A.3.6}) or (\ref{A.3.7}).   
We observe that every element of ${\cal R}$ has cardinality 
$0$, $2^d$, $2\cdot 2^d$ or $3\cdot 2^d$ and that every element of 
${\cal D} \setminus {\cal R}$ has cardinality $3\cdot 2^{d-1}$.  
To check this, just verify it for elements of ${\cal S}$, a triply thickened
[G2]~(3.19) extended Hamming code, and note that the effect of 
adding an element of ${\cal R}$  to an  element  $d \in {\cal D}$ is, 
for each $i$, to take the $i^{\rm th}$ projection of $d$ to itself or its 
complement with respect to $\Gamma_i$.  So, $F$ fixes ${\cal R}$.  
\end{notation}

\begin{lem}\label{A.3.9}   $F:=Aut({\cal D})$ permutes the partition 
$\Gamma_i$,  $i=1$, $2$, $3$, as $Sym_3$.  
\end{lem}

\pf
Since $F$ preserves ${\cal R}$, we deduce that $F$ preserves the partition by
examining the three minimal weight elements of ${\cal R}$. 
On the other hand, any blockwise permutation fixes the set of rows of $M$
(permutes the first three, fixes the rest).  \qed

\begin{notation}\label{A.3.10} 
Let $H$ be the subgroup of $F$ which fixes each $\Gamma_i$;   
the code ${\cal S}$ (\ref{A.3.8}) is a triply thickened version of 
the $d$-dimensional length $2^d$ code associated to $GL(d,2)$, as 
in (\ref{A.3.6}).  It is clear that the natural
action of a group  $F_0 \cong GL(d,2) \times Sym_3$  (first factor $F_1$ 
acting diagonally and the second $F_2$ as block permutations) is in $F$ 
and stabilizes ${\cal S}$.  
Note that the second factor acts trivially on ${\cal S}$. \qed

\end{notation}

{\bf Proof of Theorem~\ref{A.3.3}:}

There is a group $T_i$ acting as translations on $\Gamma_i$, identified with
$\FF_2^d$ as in (\ref{A.3.4}), and trivially on $\Gamma_j$, for $j \ne i$; we  
choose these identifications to be compatible with the action of $F_2$.  
The direct product $T:=T_1 \times  T_2 \times T_3$ is in $F$.  

Since $H$ fixes ${\cal R}$, we consider the action of $H$ on 
${\cal D}/{\cal R}$.  The kernel of this action corresponds naturally to
a subgroup $Hom({\cal D}/{\cal R},{\cal R})$, 
order $2^{3d}$,  and may be interpreted as an element of $T$ as in the above
paragraph.   Since  $T \le H$  and $|T|=2^{3d}$, this kernel is  $T$.  
Since $F_1 \le H$ induces $GL({\cal D}/{\cal R})$
on ${\cal D}/{\cal R}$,  $H=TF_1$ and $F=TF_0$.   
\qed 

\smallskip

\begin{rem}{\rm A study of the automorphism groups of the Reed-Muller codes, 
a class which include the simplex codes and Hamming codes, can be also found 
in~[AK], Chapter~5, for example.}
\end{rem}

\section{Lifting minus the identity}\label{append-lift}

\def\e{\varepsilon}
\def\ee{$\varepsilon$}
\def\hl{$\hat L$ }
\def\til{^{\widetilde {}} \ }

\begin{de} 
{\rm Let $L$ be an even integral lattice. 
A {\em lift of $-1$\/} is an automorphism $\theta$ of the lattice VOA
$V_L$ such that for all $x \in L$,  there is a scalar $c_x$ so that $\theta
\colon e^x \mapsto  c_x e^{-x}$. (Here, $e^x$ means $1 \otimes e^x$, 
where $1$ is the constant polynomial.) }
\end{de}

As usual, there is an epsilon function in the description of
the lattice VOA $V_L$, $\varepsilon : L \times L \rightarrow \C^\times$, 
which is bimultiplicative and satisfies 
$\varepsilon (x,y) \varepsilon (y,x)^{-1} = (-1)^{(x,y)}$.  

\begin{lem}\label{a4.1} Let $x$, $y \in L$.  For some integer $k$ and scalar 
$c$,  $e^x_ke^y=ce^{x+y}$  ($a_kb$ means the value of the $k^{\rm th}$ binary 
composition on $a$, $b$).  In  fact, we take  $k=-1-(x,y)$ and $c = \e (x,y)$,
which is always nonzero.    
\end{lem}

\pf This is obvious from the form of the vertex operator representing $e^x$. 
\qed

\begin{lem}\label{a4.2} If the set $S=-S$ spans $L$, then the set of all 
$e^x$, for $x \in S$, generates the associated lattice VOA $V_L$.  
\end{lem}

\pf  By Lemma~\ref{a4.1}, we may assume that $S=L$. 
Let $V'$ be the subVOA so generated. Note that for any $\alpha\in L$,
$\alpha(-1)=e^{\alpha}_{(\alpha,\alpha)-2}e^{-\alpha}$. Thus $V'$ contains
all $e^{\alpha}$ and $\alpha(-1)$ for $\alpha\in S$. 
It is clear that $V_L$ is irreducible under the component operators
of $Y(e^{\alpha},z)$ and  $Y(\alpha(-1),z)$ for $\alpha\in S$, 
hence $V'$ contains all $p \otimes e^{\alpha}$, where $p$ is a polynomial
expression in $\alpha (n)$, for $n < 0$.  
It follows immediately that $V'=V_L$. \qed

\begin{notation} 
{\rm  Let $M$ be the set of lifts of $-1$  and $T$ the rank $\ell$ torus
of automorphisms of $V_L$ associated to $L$.  There is an identification 
$T= \C^\ell / L^*$ so that $t=v+L^*  \in T$ sends $e^x$ to $e^{2\pi i(v,x)}e^x$.}
\end{notation}

\begin{lem}\label{a4.3} Let $A$ be an abelian group, $\<u\>$ 
be a group of order 2 which acts on $A$ by letting $u$ invert
 every element of $A.$ Set $B:=A\langle u \rangle,$ the semidirect 
product. Every element of the coset $Au$ is an involution, and two such 
involutions $cu$ and $du$ are conjugate in $B$ (equivalently, by an 
element of $A$) iff $cd^{-1}$ is the square of an element of $A$.  
This last condition follows if $A$ is divisible, e.g.~a torus.
\end{lem}  

\begin{thm}
$M$ forms an orbit under conjugation by $T$ in ${ Aut}(V_L)$.  
\end{thm}

\pf  Let $x_1$, $\dots$, $x_\ell$ form a basis of $L$.   Given 
an element of $M$, we may compose it with an element  $r \in T$ to assume it 
satisfies  $e^{\pm  x_i} \mapsto e^{\mp x_i} $, for all $i$.  
The conditions  $e^{\pm  x_i} \mapsto  e^{\mp x_i} $ 
characterize an automorphism, since these $2\ell$ 
elements generate the VOA, by Lemma~\ref{a4.2}.  
This composition is the same as conjugation  by $s \in T$
such that $s^2=r$ or $r^{-1}$.  
So, we are done if we prove that $e^x \mapsto e^{-x}$ for all $x \in L$.  But
this is clear from Lemma~\ref{a4.1} since 
$\varepsilon (-x, -y) = \varepsilon (x,y)$. 

\begin{cor} Given two lifts of $-1$ on $V_L$, their fixed point subVOAs 
are isomorphic.  In fact, these subVOAs are in the same orbit of $Aut(V_L)$.  
\end{cor}


\medskip
\noindent Department of Mathematics, University of California, 
Santa Cruz, CA 95064 USA (C.D.)

\noindent Department of Mathematics, University of Michigan, 
Ann Arbor, MI 48109-1109 USA (R.L.G.)

\noindent Mathematisches Institut, Universit\"at Freiburg, Eckerstr. 1, 
D-79104 Freiburg, Germany (G.H.)

\end{document}